\documentclass[11pt]{article}

\usepackage{pdfsync}
\usepackage{placeins,graphicx,diagbox }
\usepackage{bm}
\usepackage{array}
\usepackage{amsmath,amssymb}  
\usepackage{pgfplots}

\usepackage{amsmath,algorithm} 
\usepackage{amssymb}  
\usepackage{color}
\usepackage{cite}

\newtheorem{assumption}{Assumption}
\newtheorem{corollary}{Corollary}
\usepackage{algorithmic}

\usepackage{graphicx,epstopdf}

\usepackage{amsopn}

\usepackage{xspace}
\usepackage{bold-extra}
\usepackage{placeins,diagbox }

\usepackage{wrapfig} 
\usepackage{tikz,pgfplots,mathtools,float}
\usepackage{booktabs,caption}
\usepackage[flushleft]{threeparttable}
\usepackage{subfigure}
\usepackage{multirow}
\usepackage{amssymb,latexsym}  
\usepackage{extarrows}
\usepackage{dsfont} 
\usepackage{color}
\usepackage{cite}

\def\Real{\mathbb{R}}
\def\argmin{\mathop{\rm argmin}}

\newcommand{\blue}[1]{{\color{black}#1}}

\allowdisplaybreaks
\usepackage{graphicx,bm}
\usepackage{amsfonts,amsmath,fullpage,bbm}
\usepackage{amssymb,multirow,verbatim}
\usepackage{acronym,wrapfig,plain,mathrsfs,enumerate,relsize,color}
\usepackage{subfig}
\usepackage{algorithmic}
\usepackage{pifont}

\usepackage{wrapfig} 
\usepackage{times} 
\usepackage{amsmath,mathtools,float}
\usepackage{tikz,pgfplots}
\usepackage{booktabs,caption}
\usepackage[flushleft]{threeparttable}
\usetikzlibrary{arrows,shapes,snakes,shadows,positioning,automata,patterns}
\usetikzlibrary{trees,decorations.pathmorphing,decorations.markings}
\usepackage{amsmath}
\usepackage{multirow}
\usepackage{amssymb,latexsym}  
\usepackage{dsfont} 

\usepackage{amsmath,algorithm} 
\usepackage{amssymb}  
\usepackage{color}
\usepackage{cite}

\def\Real{\mathbb{R}}

\def\argmin{\mathop{\rm argmin}}

\newtheorem{theorem}{Theorem}
\newtheorem{lemma}{Lemma}

\newtheorem{proposition}{Proposition}

\newtheorem{remark}{Remark}

\usepackage{verbatim}
\usepackage{ulem}
\newcommand{\vus}[1]{{\color{black}#1}}

\newcommand{\uss}[1]{{\color{black}#1}}

\newcommand{\red}[1]{{\color{black}#1}}

\def\be{\begin{enumerate}}
\def\ee{\end{enumerate}}

\def\argmin{\mathop{\rm argmin}}

 \newcommand{\remove}[1]{}

\def\Real{\mathbb{R}}

\def\argmin{\mathop{\rm argmin}}

\newcommand{\bxi}{\boldsymbol{\xi}}

\usepackage{tikz}
\usetikzlibrary{decorations.pathreplacing,calc}

\date{}

\begin{document}
\title{Variance-Reduced  Accelerated First-order Methods: Central Limit Theorems and Confidence Statements}

\author{  Jinlong  Lei  and  Uday  V. Shanbhag \thanks{ Email:  leijinlong@tongji.edu.cn (Jinlong Lei), udaybag@psu.edu
(Uday  V. Shanbhag). } \thanks{  Lei is  the Department of Control Science and Engineering, Tongji University; She is also with the Shanghai Institute of Intelligent Science and Technology, Tongji University, Shanghai 200092, China. Shanbhag is with
the Department of Industrial and Manufacturing Engineering, Pennsylvania State
University, University Park, PA 16802, USA}}

\maketitle

\begin{abstract} In this paper, we consider a  strongly convex stochastic
    optimization problem and propose  three classes of variable sample-size
    stochastic first-order methods:  (i) the standard stochastic gradient
    descent method; (ii) its accelerated variant; and  (iii) the stochastic
    heavy ball method.  In each scheme,   the exact gradients are approximated
    by averaging across an increasing batch size of sampled gradients.  We
    prove that  when the sample-size    increases geometrically,  the generated
    estimates converge  in mean to the optimal solution at a geometric rate for
    schemes (i) -- (iii).  Based on this result, we provide central limit
    statements whereby it is shown that the rescaled estimation errors converge
    in distribution to a normal distribution with the associated  covariance
    matrix   dependent   on the Hessian matrix, covariance  of the gradient
    noise, and the steplength. If the sample-size   increases  at a polynomial
    rate, we show that the estimation errors decay at a corresponding
    polynomial rate and  establish the  associated central limit theorems
    (CLTs). Under  certain conditions, we discuss how both the algorithms
    and the associated limit theorems may be extended to constrained and
    nonsmooth regimes.   Finally, we provide an avenue to construct confidence
    regions for the optimal solution based on the established CLTs  and test
    the  theoretical  findings on a stochastic parameter estimation problem.
\end{abstract}

  \section{Introduction}
  In this paper, we consider the strongly convex optimization problem \eqref{Problem1}:
\begin{align}\label{Problem1} \min_{x\in \mathbb{R}^m} f(x) \triangleq
  \mathbb{E} [F(x,\bxi)],
\end{align}
where  $\bxi:  \Omega \to \Real^d$  is a random variable  {defined on the  probability space} $({ \Omega}, {\cal F},\mathbb{P})$, $F: \mathbb{R}^m \times \mathbb{R}^d \to \mathbb{R}$, and   the expectation is taken over the distribution of the random vector. Stochastic optimization problems  have been  extensively studied,
 given the   wide applicability  of such models  in almost all areas of
science and engineering, ranging from communication and queueing systems  to
finance (cf.~\cite{birge97introduction,shapiro2009lectures}). However, in most
situations,  this expectation   and  its derivative are unavailable in closed
form requiring the development of sampling approaches. Sample average approximation (SAA)~\cite{shapiro2009lectures} and   stochastic approximation (SA)  represent two  commonly used approaches  for contending with stochastic programs.
Our focus is on SA  schemes,  first considered by Robbins and Monro
\cite{robbins1951stochastic} for seeking roots of a regression
function with noisy observations.  The   standard  SA algorithm
$x_{k+1}=x_k-\alpha_k \nabla F(x_k,\xi_{k}),$ also known as  the
stochastic gradient  descent (SGD) algorithm,   updates the estimate
$x_{k+1}$ based on  a sampled gradient  $\nabla F(x_k,\xi_{k})$.  The  convergence analysis usually requires
suitable properties on the gradient map (such as Lipschitzian
requirements) and the steplength sequence (such as non-summable but
square summable). The almost sure convergence of $x_k$ to
$x^*$, the unique optimal solution of \eqref{Problem1},  was
established in
\cite{blum1954multidimensional,bottou1998online,chen2006stochastic}
on the basis of the  Robbins-Siegmund theorem
\cite{robbins1971convergence}   while ordinary differential equation (ODE) techniques  have also been
employed for claiming similar statements   in
\cite{hasminskii1973stochastic,kushner2003stochastic,borkar2009stochastic}.
 In addition, in~\cite{chee2018convergence}, the authors developed a statistical diagnostic test to detect the phase transition, during which the  iterative procedure converges towards a region of interest in the context of SGD with constant learning rate.
Tight bounds on the  rate of convergence can be obtained by establishing
the asymptotic distribution for the iterates (cf.
\cite{chen2006stochastic,hasminskii1973stochastic,kushner2003stochastic,borkar2009stochastic,fabian1968asymptotic}).
An instructive review of results up to around 2010 was provided in
\cite{pasupathy2011stochastic},  while  asymptotic normality of the suitably  scaled  iterates for   SA
with decreasing step-sizes \vus{has been proven in~}\cite{fabian1968asymptotic,hasminskii1973stochastic,kushner2003stochastic}. To be specific,  the rescaled   error process
$\sqrt{ \alpha_k}(x_k-x^*)$ asymptotically converges  in distribution to a
normal distribution  with zero mean  {and} with covariance depending on the
Hessian matrix, the covariance of the gradient noise, and the steplength.  Prior work on CLTs  for  standard stochastic approximation for smooth convex optimization can be traced to the seminal averaging paper by Polyak and Juditsky~\cite{polyak1992acceleration}.  The asymptotic normality was further investigated in \cite{chen2006stochastic}
for SA with expanding truncations, an avenue that does not require Lipschitz continuity   of the gradients.    The sequence of iterates
generated by  the  constant steplength SA scheme   has been   shown to be a homogeneous Markov chain with a unique stationary distribution;  see
\cite[Chapter 9]{kushner2003stochastic},  \cite[Chapter 17]{meyn2012markov}, and \cite{dieuleveut2020bridging}.

CLTs for SA schemes are significant from
the standpoint of algorithm design as well as inference.

\noindent  {\em (i)  Algorithm design.} Indeed, the optimal
selection of the steplength depends on the Hessian matrix at  $x^*$, see e.g.,
\cite[Chapter 3.4]{chen2006stochastic}. Since the  Hessian at $x^*$ is
unavailable, the optimal value of $\alpha_k$ can only be estimated, which has  led to the development of   adaptive SA methods (cf. \cite{venter1967extension,wei1987multivariate,spall2000adaptive}).
Motivated by the heavy dependence of SA schemes on steplength choices,
 the authors in \cite{yousefian10stochastic} developed a self-tuned rule that
adapts the steplength sequence to problem parameters. In a similar vein,
there has been an effort to develop optimal constant steplengths.
Specifically,    it was shown {in \cite{nemirovski2009robust}}  that with   suitably selected
 constant  step-sizes, the  expected  function values   at the averaged
iterate  converge to the optimum with rate  $\mathcal{O}(1/\sqrt{k})$ in
merely convex case and  rate $\mathcal{O}(1/k)$  in  strongly convex case,
matching the lower bounds~\cite{nemirovskii1983problem}.
 It is further  shown in   \cite{mandt2017stochastic} that the constant SGD  simulates a Markov chain
with a stationary distribution, which   might   be  used  to adjust the tuning parameters of constant SGD
 so as to improve the  convergence rate,     e.g.,  SGD for fitting  the  generalized linear
models \cite{toulis2014statistical}   and approximating the  Bayesian posterior
inference \cite{mandt2017stochastic}.

\noindent  {\em (ii) Confidence statements.} Furthermore, CLTs for SA
schemes might  allow for  the possibility  of constructing  confidence
regions for the optimal solution (e.g.~\cite{bickel2013asymptotic,chen2020statistical, hsieh2002recent}).   In particular, {Hsieh and Glynn~}\cite{hsieh2002recent}   designed an approach to rigorously characterize confidence regions without explicitly estimating the covariance
of the limiting normal distribution,  while  {Chen et al.}~\cite{chen2020statistical} proposed a plug-in estimator  and  a batch-means estimator for the asymptotic covariance of the average iterate from SGD.  In addition,  {Su and Zhu}~\cite{su2018uncertainty} advocated  performing SGD updates for a while and then  {splitting} the single thread into several threads; {notably, they proceeded to} construct $t$-based confidence interval  that can attain asymptotically exact coverage probability.
 The   {monograph~}\cite{kushner2003stochastic}  has extensively investigated CLTs for  SA schemes  under both constant and decreasing steplengths.
In the context of SAA schemes, there has been some recent work on developing confidence statements for  {stochastic optimization~\cite{shapiro2009lectures} and} stochastic variational inequality problems~\cite{lu2014new,LammLB17}.

Unfortunately,  SA schemes with  diminishing steps cannot  recover  the  deterministic convergence rates
 seen in exact gradient methods  while constant
steplength SA schemes are only characterized by convergence
guarantees to a  neighborhood of the optimal solution.
 Variance-reduction schemes employing an increasing batch-size of
sampled gradients (instead of the unavailable true gradient) appear
to have been first alluded to
in~\cite{bertsekas1989parallel,hdmello03variable,deng09variable,pasupathy10choosing} and
analyzed in smooth and strongly convex~\cite{friedlander12hybrid,
byrd12sample, shanbhag2015budget, schmidt12convergence}, smooth
convex regimes~\cite{ghadimi2016accelerated}, nonsmooth (but
smoothable) convex~\cite{jalilzadeh2018optimal},
nonconvex~\cite{lei18asynchronous},  and
game-theoretic~\cite{lei18game} regimes. Notably,  the linear rates in mean-squared error were  derived
for strongly convex smooth~\cite{shanbhag2015budget} and a
subclass of nonsmooth objectives~\cite{jalilzadeh2018optimal}, while
a   rate $\mathcal{O}(1/k^2)$ and $\mathcal{O}(1/k)$ was obtained
for expected sub-optimality  in convex
smooth~\cite{ghadimi2016accelerated, jofre19variance} and
nonsmooth~\cite{jalilzadeh2018optimal}, respectively.  In each instance, the
schemes achieve the corresponding optimal deterministic rates of
convergence under suitable growth in
sample-size sequences while in almost all cases, the optimal sample-complexity bounds were obtained.  An excellent
discussion relating sampling rates and the canonical rate (i.e. the Monte-Carlo rate) in stochastic
gradient-based schemes has been provided by Pasupathy et
al.~\cite{pasupathy18sampling}.\\

 {\bf Gaps and motivation.}  This paper is motivated by
the following gaps. (i) Despite a surge of  {interest} in
variance-reduced schemes, both  with and without acceleration, no limit
theorems are available for claiming asymptotic normality of the
  scaled sequence in either unaccelerated or Nesterov accelerated regimes coupled
with variance reduction. In particular, we have little understanding
regarding whether such avenues have detrimental impacts in terms of
the limiting behavior. Such statements are also unavailable for the
related heavy-ball method. (ii) The geometric growth in sample-size
required for achieving linear rates of convergence may prove onerous
in some settings.  Are CLTs available in settings where the
sample-size grows at a polynomial rate? (iii) Finally, given a CLT,
there is little by way of availability of rigorous confidence
statements, barring the work by Hsieh and
Glynn~\cite{hsieh2002recent}. Can such statements be developed for
the proposed class of variance-reduced first-order methods?\\

 {\bf Justification and relationship to other variance-reduced schemes.}

\noindent (i) {\em Terminology and applicability.} The moniker
``variance-reduced'' has been loosely used in stochastic optimization to
capture approaches that admit deterministic convergence rates. For instance,
for minimizing  finite-sum expected-residual minimization  problems in
machine learning, techniques such as stochastic variance reduced gradient (SVRG)~\cite{johnson2013accelerating} and
a fast incremental gradient method~\cite{defazio2014saga} achieve deterministic rates of convergence. In
general probability spaces, using increasing batch-sizes of gradients represent
progressively accurate approximations of the true gradient, as opposed to noisy
sampled variants employed in single sample schemes~\cite{friedlander12hybrid}.
The resulting schemes, often referred to as {\em mini-batch} SA schemes,
 may achieve deterministic rates of convergence under suitable growth rates
in sample-sizes.

\noindent (ii) {\em Weaker assumptions, and stronger statements.} The proposed
variance-reduced framework has several key benefits often unavailable in
single-sample regimes:  1) Under suitable assumptions, such schemes often
achieve optimal deterministic rates in terms of iteration complexity
 (e.g. in constrained regimes, complexity of an iteration is essentially that of computing a  projection onto a convex set)  while achieving near-optimal sample complexity. This assumes profound importance when
iterations are expensive as in high-dimensional and nonlinear settings;   2) In addition, the benefits of acceleration are less clear in stochastic single-sample regimes, absent such approaches (to the best of our knowledge).  While the present treatment is in an unconstrained regime, we also discuss how extensions to constrained and nonsmooth regimes may be developed.

\noindent (iii) {\em Sampling requirements.} Variance-reduced
schemes have obvious benefits when sampling is relatively cheap compared to computational effort of an iteration.  Though such schemes are often characterized by near-optimal
sample complexity,   one might question how to contend with batch-sizes (denoted by $N_k$) tending to $+\infty$. This
issue is somewhat of a red herring since most SA schemes are meant to provide
$\epsilon$-approximations; (see the detailed discussions   in Remark \ref{rem-geo}).     
With the   ubiquity of multi-core architecture, such
requirements are not terribly onerous.  {In addition, we also establish  the   polynomial rates and associated sample complexities for polynomially increasing sample size.  They are particularly important when sampling is expensive and one would like to modulate the sampling rate to ensure that the scheme is practical.    }

 \noindent (iv)   {\em  CLTs  and  confidence statements.}  Apart from  rate statements, we additionally establish  CLT results for  the proposed variance-reduced framework.
 In particular,  we provide a clean characterization of the limiting distribution and show the   dependence  of the  Hessian matrix,  the condition number,  the noise  covariance, etc   in the prescription of the CLT.
The confidence statements are established as well;  we should note that we are unaware  of such statements in the context of accelerated and heavy-ball settings in a variance-reduced regime.   	 \\

  {\bf Outline and contributions.}  To address these gaps, we present CLTs and confidence statements for   first-order
stochastic variance-reduced algorithms for resolving \eqref{Problem1}, including the classical
SGD~\cite{robbins1951stochastic}, the stochastic variants of  the
Nesterov's  accelerated  method~\cite{nesterov1983method},  and the
heavy ball method~\cite{polyak1964some}.  We provide statements
when batch-sizes increases at either a geometric or a polynomial
rate. Our main contributions are summarized next.

\noindent   {{\bf (I) CLTs for variable sample-size gradient
methods.}  In Section \ref{Sec:Alg:SGD}, we recall the
variance-reduced (VR) stochastic gradient algorithm
with  a constant steplength  (see Algorithm
\ref{Alg_1}), where the  gradient is  estimated by  the average
of  an  increasing  batch of sampled gradients. In Section
\ref{sec:CLT-geo}, when the  batch-size increases at a
geometric rate,  we observe that the mean-squared error
diminishes at a geometric   rate (see Proposition
\ref{prp1}) and   provide  a preliminary Lemma
\ref{lem-CLT2}  for  establishing  a  CLT for a
noised-corrupted     linear  recursion. Based on  this Lemma
and    the linear approximation  of the gradient function
at  $x^*$,  we  may derive a  CLT (see Theorem
\ref{thm-CLT1}) in this setting. We proceed to show}  that
the  covariance   of the limiting normal distribution
depends on the  Hessian matrix at the solution,  the
condition number,  the covariance of   gradient noise, etc
and the steplength.  Additionally, we show in    Section
\ref{sec:CLT-poly} that when the batch-size  is increased at
a polynomial rate, the sequence of iterates  converges
at  a corresponding   polynomial  rate (see Proposition
\ref{prp3}).  Then based on the CLT shown in Lemma
\ref{lem-CLT3}  for the time-varying linear recursion, the
CLT for Algorithm \ref{Alg_1}  with polynomially increasing
batch-size is established in  Theorem \ref{thm-CLT3}.

\noindent {\bf (II)  CLTs for VR-accelerated gradient method.} In Section \ref{Sec:Alg:ACC},  we consider a VR-accelerated
gradient algorithm with constant steplengths (see Algorithm \ref{Alg_2}). Then by
leveraging the  geometric rate of convergence (Proposition \ref{prp2}),
we establish amongst the first CLTs in accelerated regimes
(Theorem \ref{thm-CLT1-Alg2})  when the batch-size increases at a
geometric rate.   It is well-known that the accelerated variant
has better constants in the iteration complexity in deterministic regimes and this is also seen in stochastic settings.
 Akin to earlier,  when the  batch-size increases at a polynomial rate, a polynomial convergence   rate and the corresponding CLT are   in  Proposition \ref{prp4}   and Theorem \ref{thm-CLT2-Alg2}.

\noindent {\bf III. CLTs for VR-heavy-ball schemes.} In
Section \ref{Sec:Alg:HB},   we design  a VR-heavy ball scheme with constant steplengths (see Algorithm
\ref{Alg_3}).  \red{In contrast with  Algorithms \ref{Alg_1} and \ref{Alg_2},  the heavy ball method is analyzed for both  quadratic and non-quadratic objective functions. When the batch-size increases geometrically,  the  geometric    rate  of convergence are respectively  shown in Proposition \ref{prp-alg3-Q} and  \ref{prp-alg3} for    quadratic and non-quadratic objective functions, while the corresponding    CLT is established in Theorem \ref{thm-CLT1-Alg3} and Theorem \ref{thm-CLT1-Alg3-NQ}. Similarly, for  polynomially  increasing  batch-size, a
polynomial  rate  of convergence is shown in Proposition   \ref{prp2-alg3-Q} and \ref{prp2-alg3} for   quadratic and non-quadratic objective functions,   whereas the corresponding    CLTs  are respectively
established in   Theorem \ref{thm-CLT2-Alg3} and Theorem \ref{thm8}.}

\noindent {\bf IV. Confidence statements.} In Section
\ref{sec:CI}, inspired by  \cite{hsieh2002recent},  we
 provide rigorous  confidence regions  for the optimal solution
and  function value.  Then in Section \ref {sec:simu},
we  implement  some simulations on a  parameter estimation  problem
in the stochastic  environment to validate the theoretical findings.\\

   {\bf Notations.} Let $\mathbf{I}_m$ denote the identity matrix of dimension $m $ and $\mathbf{0}_m \in \mathbb{R}^{m \times m}$  denote   the   matrix with all entries equal zero.
  Let $\{X_k\}$ be a sequence of   random variables.   $X_k  \xrightarrow [k
  \rightarrow \infty]{d}   N(0,\mathbf{S}) $ denotes that $ X_k $ converges  in
  distribution  to a  normal distribution  $N(0, \mathbf{S})$   with mean zero
  and covariance $\mathbf{S}, $ and $X_k  \xrightarrow [k \rightarrow
  \infty]{P} X $ denotes that  $ X_k $ converges  in probability  to $X.$
 By utilizing  the definitions of   the ``small-$o$ in probability'' and ``big-$O$ in probability'' notations $o_P(\cdot)$ and $O_P(\cdot)$
  provided  in  \cite[Chapter 2.2]{van2000asymptotic}, the expression   ${e_k}  { \, = \, } o_P(1)$ implies that a sequence of random variables $\{e_k\}$  converges to
  zero in probability. Similarly, the  expression  $e_k  { \, = \, } O_P(1)$ implies
  that  a sequence of random variables   $\{e_k\}$ is bounded in probability. For a square matrix $ \mathbf{P}$, we denote  by   $ \rho(  \mathbf{P})$ and $\|  \mathbf{P}\|$ its   spectral radius and matrix two-norm, respectively.

\section{First-Order Variable Sample-size  Stochastic   Algorithms.}\label{Sec:Alg}
Since the exact gradient $\nabla f(x)$  is  expectation-valued and unavailable
in a closed form,   we assume  that  there exists a {\em stochastic first-order
oracle}   such that  for any $x$ and $\xi$, a  sampled gradient $\nabla F (x,\xi)$ is returned  and is
assumed to be an unbiased estimator of  $\nabla f(x)$.  In
this section, we present three  first-order  stochastic  algorithms to
find the optimal solution to \eqref{Problem1}. Throughout the
paper, time is slotted at $k=0,1,2,\dots$ and an iterate at time $k$ is denoted by $x_k\in \mathbb{R}^m.$

\subsection{Variance-reduced Gradient Method.}\label{Sec:Alg:SGD}
 We  present  a variable sample-size stochastic gradient algorithm  (Algorithm \ref{Alg_1}) to solve \eqref{Problem1}, where at iteration $k $,   the unavailable exact gradient $\nabla f(x_k)$ is estimated via the average of  an increasing batch-size of  sampled gradients.
 \begin{algorithm}[H]  \caption {Variance reduced SGD}\label{Alg_1}
Given an arbitrary  initial value $x_0\in \mathbb{R}^m, $  a positive constant $\alpha>0,$ and a positive  integer sequence $\{N_k\}_{k\geq 0}.$ Then iterate the following equation  for $k\geq 0.$
\begin{align}\label{Alg1} x_{k+1}=x_k - \alpha \tfrac{ \sum_{j=1}^{N_k} \nabla F(x_k,\xi_{j,k})}{N_k},
\end{align}
 where  $ \alpha>0$ is the constant steplength, $N_k$  is the number of sampled gradients used   at   time  $k$, and $  \xi_{j,k}, j=1, \cdots, N_k,$ denote the  {independent and identically distributed (i.i.d.)}
  realizations of \red{$\bxi$}.
\end{algorithm}

If  the gradient observation   noise $w_{k,N_k}$ is defined as
\begin{align} \label{def-noise}
w_{k,N_k}\triangleq \tfrac{ \sum_{j=1}^{N_k} \nabla F (x_k,\xi_{j,k})}{N_k}- \nabla f(x_k),
\end{align}
 then the  update \eqref{Alg1} can be rewritten as
\begin{align}\label{ref_alg1}
x_{k+1}=x_k - \alpha  (\nabla f(x_k)+w_{k,N_k}).
\end{align}
  For any $k$, define $\mathcal{F}_k$ as   $\mathcal{F}_k \triangleq \sigma\{x_0, \xi_{  j,t}, 1\leq j \leq N_{t} ,     0\leq t\leq k-1  \}.$
 Then   $x_k$ is adapted to $\mathcal{F}_k$ by   Algorithm  \ref{Alg_1}.
 We impose the following  conditions on the  objective function, the conditional expectation and the second
moments of the sampled gradients   produced by the stochastic first-order oracle.

  \begin{assumption}\label{ass-fun}
      (i)  $f$ is continuously differentiable on $\mathbb{R}^m$ with a Lipschitz continuous gradient, i.e., there exists a constant  $L>0$ such that $ \| \nabla f(x) -\nabla f(x')  \|  \leq L \| x-x'\| $ for any $x,x'\in \mathbb{R}^m $.
   (ii) $f$ is $\eta$-strongly convex, i.e.,
    $ ( \nabla f(x) -\nabla f(x') )^T (x-x')  \geq \eta \| x-x'\|^2$ for any $ x,x'\in \mathbb{R}^m .$ \\
\noindent (iii)
There exists  a constant $\nu>0$ such that for any $k\geq 0 $ and $j=1,\cdots, N_k,$
      $ \mathbb{E}\left[\nabla F (x_k, \bxi)\,\mid \, \mathcal{F}_k\right]=\nabla f (x_k)$ almost surely and
      $ \mathbb{E}\left[\| \nabla F (x_k, \bxi)-\nabla f (x_k) \|^2\, \mid \, \mathcal{F}_k\right] \leq \nu^2$
 almost surely (a.s.).
    \end{assumption}

Since  $x_k$ is adapted to $\mathcal{F}_k$ and the   samples $  \xi_{j,k}, j=1, \cdots, N_k,$  are    independent,   we obtain from \eqref{def-noise} and Assumption \ref{ass-fun}(iii) that  for any $k\geq 0$,
\begin{equation}\label{ass-bd-noise}
\begin{split}
    \mathbb{E}[w_{k,N_k}\, \mid \, \mathcal{F}_k]&=\tfrac{ \sum_{j=1}^{N_k}   \mathbb{E}\big[ \nabla F(x_k, \bxi ) - \nabla f(x_k) \, \mid \, \mathcal{F}_k\big]  }{ N_k}=0, \mbox{ and }
    \\  \mathbb{E}[\|w_{k,N_k}\|^2\, \mid \, \mathcal{F}_k]&=  \tfrac{\sum_{j=1} ^{N_k} \mathbb{E}   [  \|  \nabla F(x_k, \bxi ) - \nabla f(x_k) \|^2 \mid  \mathcal{F}_k    ]}{N_k^2}  \leq \tfrac{\nu^2}{N_k}.
\end{split}
\end{equation}
hold almost surely.

Since $f$ is strongly convex,  it has a unique optimal solution denoted by  $x^*$. Then by the optimality condition, $\nabla f(x^*)=0.$ We now introduce an inequality \cite[Eqn. (2.1.24)]{nesterov2013introductory} on  $f$  satisfying   Assumptions \ref{ass-fun}(i) and \ref{ass-fun}(ii):
\begin{align}\label{inequ-f}
(x-y)^T(\nabla f(x) -\nabla f(y)) \geq \tfrac{\eta L\|x-y\|^2}{\eta+L} +\tfrac{\|\nabla f(x)- \nabla f(y)\|^2}{\eta+L},\quad \forall x,y\in \mathbb{R}^m.
\end{align}
We then establish a  recursion on   the mean-squared estimation  error,
which can be proved by making a simple modification to  the proof of \cite[Theorem 2.1.15]{nesterov2013introductory}.
\begin{lemma}\label{Lemma_1}
    Let  Assumption \ref{ass-fun} hold. Consider    Algorithm  \ref{Alg_1}  with $\alpha \in $( $0, 2\over \eta+L$].
Then 
\begin{align}\label{recursion1}
\mathbb{E}[\| x_{k+1}-x^*\|^2  ] &   \leq
\left(1-\tfrac{2 \alpha \eta L}{\eta+L} \right)  \mathbb{E}[\| x_k - x^*\|^2]+\tfrac{ \alpha^2   \nu^2}{N_k},\quad \forall k\geq 0 .
\end{align}
\end{lemma}
{\bf Proof.} From \eqref{ref_alg1} it follows that
\begin{equation}\label{ms-x}
\begin{split}
\| x_{k+1}-x^*\|^2&
 =\| x_k - x^*\|^2 - 2 \alpha (x_k -x^*)^T  ( \nabla f(x_k)+w_{k ,N_k})  \\&+  \alpha^2  \left( \|\nabla f(x_k) \|^2 +2\nabla f(x_k)^T w_{k,N_k}+\|w_{k,N_k}\|^2 \right).
\end{split}
\end{equation}
Since $\nabla f(x^*)=0,$ by using  \eqref{inequ-f}, we obtain that
\begin{align}\label{inequ-fx} (x_k -x^*)^T   \nabla f(x_k) \geq \tfrac{\eta L }{ \eta+L} \|x_k-x^*\|^2+\tfrac{1}{ \eta+L}  \|\nabla f(x_k) \|^2. \end{align}
Since  $x_k$ is adapted to $\mathcal{F}_k$, by  taking expectations conditioned on $\mathcal{F}_k$ on both sides of  \eqref{ms-x}, and using   \eqref{ass-bd-noise} from Assumption \ref{ass-fun}(iii),  we obtain that
\begin{align*} \mathbb{E}[\| x_{k+1}-x^*\|^2 |\mathcal{F}_k] & \leq\| x_k - x^*\|^2 - 2 \alpha (x_k -x^*)^T   \nabla f(x_k) +   \alpha^2  \|\nabla f(x_k) \|^2  +  \tfrac{\alpha^2   \nu^2}{ N_k}  \\&
\overset{\eqref{inequ-fx}}{\leq} \left(1-\tfrac{2 \alpha \eta L }{ \eta+L} \right) \| x_k - x^*\|^2 - \alpha \left(\tfrac{2}{ \eta+L}-\alpha \right)  \|\nabla f(x_k) \|^2 + \tfrac{ \alpha^2   \nu^2 }{ N_k  }  ,\quad a.s. ~. \end{align*}
Then by  choosing  $\alpha \in (0, {2\over \eta+L}]$ and  taking  unconditional expectations,  we achieve  \eqref{recursion1}.
  \hfill $\Box$

 \subsection{Accelerated Gradient Method.}\label{Sec:Alg:ACC}
 Nesterov's  accelerated gradient   descent method \uss{generates a sequence
 that converges to the solution at}  a rate $\mathcal{O}(q^k)$ where $q
 \triangleq 1-\sqrt{\eta  / L}$ for $\eta$-strongly convex and $L$-smooth
 functions  \cite{nesterov2013introductory}, and  with a rate $
 \mathcal{O}(1/k^2 )$ for merely convex functions \cite{nesterov1983method}.
 Nesterov proved that this is the \uss{best} possible rate for any
 \uss{first-order} method. As such, we  combine the Nestrov's  accelerated
 method with  Algorithm  \ref{Alg_1}  and   propose an accelerated variable
 sample-size stochastic gradient descent algorithm (Algorithm \ref{Alg_2})  so
 as to improve the  rate of  convergence.  Such a scheme has been employed for
 smooth strongly convex~\cite{jalilzadeh2018optimal}, smooth
 convex~\cite{ghadimi2016accelerated}, and nonsmooth (but smoothable)
 convex~\cite{jalilzadeh2018optimal} stochastic optimization problems with
 associated rates of convergence given by $\mathcal{O}(q^k)$,
 $\mathcal{O}(1/k^2)$, and $\mathcal{O}(1/k)$, respectively.   The present
 paper  takes a crucial \uss{step} towards developing  CLTs and confidence
 regions for the optimal solution.
  \begin{algorithm}[H]
\caption{Variance-reduced Accelerated SGD} \label{Alg_2}
Given   arbitrary  initial values  $x_0=y_0\in \mathbb{R}^m, $   positive constants $\alpha,\beta, $ and a positive integer sequence $\{N_k\}_{k\geq 0}.$
Then iterate the following equations for $k\geq 0.$
\begin{subequations}
  \begin{align} y_{k+1}& =x_k -\alpha     (\nabla f(x_k)+w_{k,N_k}) ,\label{Alg21}
  \\ x_{k+1}&=y_{k+1}+\beta  (y_{k+1}-y_k),\label{Alg22}
\end{align}
\end{subequations}
where $w_{k,N_k} $ is defined as in \eqref{def-noise}, $\alpha>0$ and $\beta>0 $ are constant steplengths.
\end{algorithm}

Next,  we establish an upper bound on  the  expected  sub-optimality gap   of the iterates  generated by Algorithm  \ref{Alg_2}. Its  proof   is similar to that  in   \cite{nesterov2013introductory}; hence it is omitted here but included in the   supplementary material for completeness. This is an important preliminary result to be used in the rate analysis of Algorithm  \ref{Alg_2}.
 \begin{lemma}\label{lem-acc}    Let Algorithm  \ref{Alg_2}
be applied to the problem \eqref{Problem1}. Suppose Assumption \ref{ass-fun} holds and   $   \alpha \in ( 0, {1\over L}]$. Set  $\beta = {1 -\gamma \over 1+ \gamma} $  with $\gamma = \sqrt{ \alpha \eta } $. Then  for any $k\geq 1,$
\begin{align}\label{lem-fy}
    \mathbb{E}[f(y_{k }) -f^*] \leq     \tfrac{(\eta+L)(1-\gamma)^{k }}{2}     \mathbb{E}[\| x_0-x^*\|^2]  +  \nu^2 \left( \alpha +\tfrac{(1-\gamma)\gamma}{2\eta}\right)  \sum_{i=0}^{k-1}  \tfrac{(1-\gamma )^i}{N_{k-1-i}}.
\end{align}
\end{lemma}
\subsection{The Heavy-Ball Method.}   \label{Sec:Alg:HB}
The classical heavy-ball method of Polyak~\cite{polyak1964some}
takes the   form $ x_{k+1}  =x_k -\alpha \nabla f(x_k)
+\beta(x_k-x_{k-1})$  with a steplength $\alpha>0$ and a momentum parameter
$\beta>0.$  {For the class of  quadratic  and strongly convex functions,   the author in~\cite{polyak1964some} derives optimal step-size parameters and achieves  the
linear convergence rate   $(\sqrt{L}-\sqrt{\eta})/ (\sqrt{L}+\sqrt{\eta}) $.
 However  for the class of non-quadratic but $L$-smooth and $\eta$-strongly convex
functions,  the authors in  \cite{ghadimi2015global} provide the global linear rate of
convergence;  however, this rate  does not lead to an acceleration  compared with  unaccelerated
gradient descent method  (cf. \cite[Lemma 2.5]{goujaud2023provable}).}
Stochastic variants of  the  heavy-ball method have been employed widely in
practice (cf.
\cite{sutskever2013importance,szegedy2015going,krizhevsky2012imagenet}
for applications to machine learning). Recent efforts
\cite{yang2016unified,gadat2018stochastic}  have analyzed the
rate  in stochastic settings and \cite{yang2016unified}  derived
a sublinear rate $ \mathcal{O}(  1/\sqrt{k} )$  for general Lipschitz
continuous convex objectives  with bounded variance, while a rate $ \mathcal{O}(
1/k^{\beta} )$ with $\beta\in(0,1)$ was provided   in
\cite{gadat2018stochastic}   for the  case of the strongly  convex quadratic
functions.  Convergence properties remain unaddressed  in variance-reduced regimes. Therefore, we consider a variable sample-size  variant of the stochastic heavy ball method assuming constant steplengths, in contrast with   the diminishing steplengths utilized in \cite{yang2016unified,gadat2018stochastic} and prove  the global linear convergence  of the iterates.  In  the developed Algorithm \ref{Alg_3},   we   add a momentum term $\beta(x_k-x_{k-1}) $ to the variable sample-size  stochastic gradient step \eqref{Alg1} and obtain the update \eqref{Alg3}.

\begin{algorithm} [H] \caption {Variance-reduced heavy-ball SGD}\label{Alg_3}
Given an arbitrary  initial value $x_0\in \mathbb{R}^m, $  two positive constants $\alpha,\beta>0,$ and a positive  integer sequence $\{N_k\}_{k\geq 0}.$  Set $x_{-1}=x_0$. Then iterate the following equation  for $k\geq 0.$
\begin{align}\label{Alg3} x_{k+1}& =x_k -\alpha     (\nabla f(x_k)+w_{k,N_k}) +\beta(x_k-x_{k-1}),
\end{align}
where $w_{k,N_k} $ is defined as in \eqref{def-noise}, $\alpha>0$ and $\beta>0 $ are constant steplengths.
\end{algorithm}
\vspace{-0.1in}

In the following lemma, we first give an upper bound on the expected mean-squared  error of  $\{x_k\}$ for quadratic cost functions based on \cite{polyak1987introduction}.  The proof can be found in Appendix \ref{App-Lem-HB-Q}

\begin{lemma} \label{Lem-HB-Q}
  Suppose that Assumption \ref{ass-fun} holds and   $f$ is a quadratic function with   $  \nabla^2 f(x) \equiv \mathbf{H} $ for any $x\in \mathbb{R}^m.$  Consider Algorithm  \ref{Alg_3}, where     $   \alpha \in (0,4/L)$ and  $\beta \triangleq  \max\{|1-\sqrt{\alpha \eta}|^2,|1-\sqrt{\alpha L}|^2\}<1.$  Then for any $\iota \in (0,1-\sqrt{\beta}), $  there exist a  constant   $ c(\iota )$   such that the following holds for any $k\geq 0.$
 \begin{align}\label{lem-hb-q}
 \mathbb{E} \left [ \left \| \begin{pmatrix}
  x_{k+1}-x^*
 \\ x_k-x^* \end{pmatrix} \right\|^2    \right]\leq 2(c(\iota))^2(\sqrt{\beta}+\iota)^{2(k+1)}    \left [ \|
     x_{0}-x^*   \|^2 \right]+\alpha^2\nu^2 (c(\iota))^2 \sum_{t=0}^{k} \tfrac{(\sqrt{\beta}+\iota)^{2(k-t)}}{N_t}.
\end{align}
\end{lemma}

Additionally, we   follow the idea of  \cite{ghadimi2015global}  in  establishing  the convergence  rate of Algorithm  \ref{Alg_3} for  $L$-smooth and $\eta$-strongly convex functions that are non-quadratic. The proof is  similar to that of  \cite{ghadimi2015global}, hence it is omitted here but included in the supplementary material for completeness.
\begin{lemma} \label{Lem-HB}
  Suppose that Assumption~\ref{ass-fun} holds. Consider Algorithm  \ref{Alg_3}, where     $   \beta \in (0,1)$ and  $\alpha\in (0, \tfrac{2( 1-\beta) }{ L+\eta})$. Then  for any $k\geq 0.$
 \begin{align}\label{lem-hb} \mathbb{E}[\|x_k-x^*\|^2]\leq \tfrac{1}{2 \hat{m} \eta}\mathbb{E} \left [\|  x_{0}-x^*\|^2+ \hat{m} ( f(x_0)-f(x^*)) \right] q^{k+1}+ \tfrac{\alpha^2 \nu^2}{2(1-\beta)^2 \hat{m} \eta} \sum_{i=0}^k q^i/N_{k-i} .
\end{align}
\end{lemma}
where $\hat{m} \triangleq \tfrac{2\alpha  }{ 1-\beta}\left(\tfrac{\beta-\alpha \eta }{ 1-\beta}+
 \tfrac{2\eta }{\eta+L} \right)$ and $q=\max\{q_1,q_2\}$ with
\begin{align}\label{def-q1q2}
q_1\triangleq {\beta \over \beta+ \eta\left(  \tfrac{2(1-\beta)}{ L+\eta}-\alpha\right)}<1 {\rm~and~}
q_2\triangleq \max \left\{0, 1- \tfrac{2\alpha \eta L}{(1-\beta)(L+\eta)+2\beta L}\right\}.
 \end{align}

Based on the above two lemmas, for both the quadratic and non-quadratic   strongly convex functions,  we will separately derive both  a geometric rate and asymptotic normality statements in Section \ref{sec:CLT-geo} for
Algorithm \ref{Alg_3}  under geometrically  increasing batch-sizes,
whereas the corresponding  results    of Algorithm \ref{Alg_3} under
polynomially increasing batch-sizes will be  established  in
Section \ref{sec:CLT-poly}.

 \subsection{Pathway for addressing constrained and nonsmooth regimes.}
We believe that the presented avenues hold promise for contending with constrained and nonsmooth regimes.

\noindent {\em (i) Constrained problems.} Consider the constrained problem
\begin{align} \label{copt}
    \min_{x \in \mathcal{X}} \ f(x) \triangleq \mathbb{E}[F(x,\red{\bxi})],
\end{align}
where $\mathcal{X}$ is a closed and convex  set with a nonempty interior. {
    Algorithm \ref{Alg_1} can be  extended to a constrained regime with an additional  projection onto  the convex set $\mathcal X$.
Under suitable conditions,  Lemma \ref{Lemma_1} holds and $x_k$   converges almost surely to the optimal solution $x^*$. For the case when the optimal solution $x^*$ lies in the interior of $\mathcal X$,
the sequence $\{x_k\}_{k \geq K}$ will lie in the interior of the constrained set in an almost sure sense for sufficiently large $K$.
Then by proceeding in a similar fashion as in Section~\ref{sec:CLT-geo}, CLTs may be developed for constrained problems. If however, solutions are on the boundary of the set,
a possible resolution may  lie in computing CLTs of approximate minimizers of an unconstrained reformulation via penalization and barrier methods~\cite{nocedal99numerical}.

\noindent {\em (ii) Nonsmooth problems.} When the $f$ is not necessarily smooth, then one avenue lies in employing smoothing approaches. Under some conditions~\cite{beck17fom}, one can construct an $(a,b)$-smoothed convex approximation of $f$, denoted by $f_{\eta}$ where $f_{\eta}(x) \leq f(x) \leq f_{\eta}(x)+\eta b$ and $\|\nabla f_{\eta}(x) - \nabla f_{\eta}(y) \| \leq \tfrac{a}{\eta}\|x-y\|$ for all $x, y$. One may then integrate iterative smoothing within the  variance-reduced accelerated gradient scheme as follows.
    \begin{align*} y_{k+1}& =x_k -\alpha_k     (\nabla f_{\eta_k} (x_k)+ w_{k,N_k}) ,
  \\ x_{k+1}&=y_{k+1}+\beta_k  (y_{k+1}-y_k).
\end{align*}
Under suitable assumptions on $\alpha_k, \beta_k, \eta_k, N_k$, the sequence $\{y_k\}$ converges to a unique solution of the original problem. Under an assumption that $f$ is smooth in a neighborhood of $x^*$, one can again develop CLTs in this setting. If however, $f$ is not necessarily smooth at $x^*$, then one avenue might lie in developing $C^2$ smoothing-based approximations and provide CLTs for $\epsilon$-solutions.

A comprehensive examination of (i) and (ii) is beyond the scope of this paper but we believe the smooth and unconstrained analysis presented in this paper  is a crucial building block.

 \section{Central Limit Theorems  under Geometrically Increasing Batch-size.}  \label{sec:CLT-geo}
 In this section, we  establish CLTs for   Algorithms  \ref{Alg_1},   \ref{Alg_2},  and   \ref{Alg_3} when the number of sampled gradients,  denoted by $N_k$, increases at a geometric rate.

 \subsection{ Rate and Oracle Complexities.} \label{sec:geo_rate}
 Based on Lemmas \ref{Lemma_1}--\ref{Lem-HB}, we can  establish the geometric rate of convergence  along with  the iteration and oracle complexity  guarantees for  Algorithms~\ref{Alg_1}-\ref{Alg_3}. Related results regarding  linear convergence for stochastic gradient methods can be found in~\cite{byrd12sample,friedlander12hybrid,shanbhag2015budget,jofre19variance,pasupathy18sampling},
  while a linear rate of  the accelerated variants has been provided in~\cite{schmidt12convergence,jalilzadeh2018optimal}. The proofs of Propositions \ref{prp1}-\ref{prp-alg3-Q} are similar to those in our prior work  \cite[Theorem 4.2 and Corollary 4.7]{lei18asynchronous}.

 \begin{proposition}[Rate and Oracle Complexity   for Algorithm  \ref{Alg_1}]  \label{prp1}
 Let Assumption \ref{ass-fun} hold  and $   \alpha \in (0, {2\over \eta+L}]$. Consider   Algorithm  \ref{Alg_1}, where $N_{k } \triangleq \lceil \rho_1^{-(k+1)} \rceil $ for some $ \rho_1 \in (q,1) $  with $q \triangleq 1-{2 \alpha \eta L \over \eta+L}$.
   Then
  \begin{align}\label{recursion2}
\mathbb{E}[\| x_{k}-x^*\|^2  ]  \leq\rho_1^{k} \left(\mathbb{E}[ \| x_0 - x^*\|^2]+\tfrac{\alpha^2 \nu^2}{1- q/\rho_1}\right) ,\quad \forall k\geq 1.
\end{align}
  The  iteration and oracle  complexity for computing an $\epsilon-$solution, defined as $\mathbb{E}[\|x_k-x^*\|^2] \leq \epsilon$,  are    $\mathcal{O}\left (  \kappa    \ln(1/\epsilon)   \right)$ and   $\mathcal{O}\left ({ \kappa    / \epsilon} \right)$, where  $\kappa \triangleq {L \over \eta}$ denotes  the condition number.
 \end{proposition}
  {\bf Proof.} By  substituting $N_k=\lceil \rho_1^{-(k+1)} \rceil $  into \eqref{recursion1},     using $q= 1-{2 \alpha \eta L \over \eta+L} $ and $\rho_1 \in (q,1)$,  one obtains
 \begin{align*}
     \mathbb{E}[\| x_{k+1}-x^*\|^2  ]   & \leq q  \mathbb{E}[\| x_k - x^*\|^2]+\alpha^2\nu^2 \rho_1^{k+1}
 \leq q^{k+1} \mathbb{E}[ \| x_0 - x^*\|^2]+\alpha^2\nu^2 \sum_{t=0}^{k} q^t \rho_1^{k+1-t}    \\
& = q^{k+1}\mathbb{E}[ \| x_0 - x^*\|^2] +\alpha^2\nu^2  \rho_1^{k+1}\sum_{t=0}^{k} ( q/\rho_1)
  \leq \left(\mathbb{E}[ \| x_0 - x^*\|^2]+\tfrac{\alpha^2 \nu^2}{ 1- q/\rho_1}\right) \rho_1^{k+1}.
\end{align*}
Hence, we derive   \eqref{recursion2}.  Suppose  we set $\alpha\triangleq{2\over \eta+L}$ and  $\rho_1 \triangleq \left({\kappa   \over \kappa+1}\right)^2 >q=\left({\kappa-1 \over \kappa+1}\right)^2$.   From  \eqref{recursion2}, it follows that
$\mathbb{E}[\| x_{k }-x^*\|^2  ] \leq \epsilon $ for any $ k\geq K(\epsilon) $, where
\begin{align}\label{iter-bd} K(\epsilon) \triangleq \tfrac{  \ln \left(\mathbb{E}[ \| x_0 - x^*\|^2]+\tfrac{\alpha^2 \nu^2}{ 1- q/\rho_1}\right)  +\ln \left ({1/\epsilon} \right) }{ 2 \ln  \left (1+{1\over \kappa} \right)} .\end{align}
    It is noticed that $ \ln (1+{1/ \kappa})   \ge (1-\tfrac{\kappa}{\kappa+1}) = \tfrac{1}{(\kappa+1)} $.  Consequently,  the number of iterations   required  to obtain an  $\epsilon $-optimal solution in a mean-squared sense,  i.e. $\mathbb{E}[\| x  -x^*\|^2  ] \leq \epsilon$,  is $  \mathcal{O} \big(\kappa \ln \left ({1/\epsilon} \right)\big)$.
Finally, based on the iteration complexity bound \eqref{iter-bd},   we achieve the following   {oracle complexity}  bound, measured by the number of sampled gradients,  for obtaining an  $\epsilon $-optimal solution:
\begin{align*}
\sum_{k=0}^{K(\epsilon) -1 } N_k  & \leq \sum_{k=1}^{K(\epsilon)  } \rho_1^{-k}\leq \int_{1}^{K(\epsilon)+1} \rho_1^{-t} dt
  \leq { \rho_1^{-(K(\epsilon)+1)} \over \ln(1/\rho_1)} ={\mathbb{E}[ \| x_0 - x^*\|^2]+{\alpha^2 \nu^2 \over  1- q/\rho_1} \over  2 \epsilon \rho_1 \ln(1+1/\kappa)}  =\mathcal{O}\left ({ \kappa    \over \epsilon} \right) .  \qquad \blacksquare
\end{align*}

 In the following lemma, we show that the result of Proposition \ref{prp1}  holds
as well  when   the  noise condition Assumption  \ref{ass-fun}(iii) is
   replaced by  some state-dependent noise condition.
    \begin{lemma}\label{Lemma1s}
Consider  Algorithm  \ref{Alg_1},
where $N_{k } \triangleq \lceil \rho_1^{-(k+1)} \rceil$ for some $ \rho_1 \in (q,1) $  with $q \triangleq 1-{2 \alpha \eta L \over \eta+L}$.
Let $\alpha \in (0, {2\over \eta+L}]$ and Assumptions \ref{ass-fun}(i)-(ii) hold.
       In addition, suppose  that  there exist  constants $\nu_1,\nu_2>0$ such that
  \begin{equation}\label{ass-bd-noise2}
\begin{split}
 \mathbb{E}[w_{k,N_k}| \mathcal{F}_k]&= 0
{\rm ~ and ~}  \mathbb{E}[\|w_{k,N_k}\|^2| \mathcal{F}_k]   \leq \tfrac{\nu_1^2+\nu_2^2\|x_k\|^2}{N_k},~a.s.,~ \forall k\geq 0.
\end{split}
  \end{equation} Then there  exists a  constant $c>0$ such that the following  holds.
\begin{align}\label{ms-rate}
\mathbb{E}[\| x_{k+1}-x^*\|^2  ] &   \leq\rho_1^{k} \left(\mathbb{E}[ \| x_0 - x^*\|^2]+\tfrac{\alpha^2 c}{1- q/\rho_1}\right) ,\quad \forall k\geq 0.
\end{align}
        Further, the  iteration and oracle  complexity for computing an $\epsilon$-solution   are  bounded by  $\mathcal{O}\left (  \kappa    \ln(1/\epsilon)   \right)$ and   $\mathcal{O}\left ({ \kappa    / \epsilon} \right)$,
respectively.
\end{lemma}
{\bf Proof.}
Since  $x_k$ is adapted to $\mathcal{F}_k$, by  taking expectations conditioned on $\mathcal{F}_k$ on both sides of    \eqref{ms-x},   using $\alpha \in (0, {2\over \eta+L}]$,  \eqref{inequ-fx} and \eqref{ass-bd-noise2},  we obtain the following.
\begin{align*}& \mathbb{E}[\| x_{k+1}-x^*\|^2 |\mathcal{F}_k]
 \leq  \left(1-\tfrac{2 \alpha \eta L }{ \eta+L} \right) \| x_k - x^*\|^2  + \tfrac{\alpha^2 (\nu_1^2+\nu_2^2\|x_k\|^2)}{N_k}   ,\quad a.s. ~. \end{align*}
 By  taking    unconditional expectations,  using $N_{k } \triangleq \lceil \rho_1^{-(k+1)} \rceil$ and
 $\|x_k\|^2\leq 2(\|x_k-x^*\|^2+\|x^*\|^2)$,  we obtain that
 \begin{align}\label{bd-ms}& \mathbb{E}[\| x_{k+1}-x^*\|^2  ]
 \leq  \left(1-\tfrac{2 \alpha \eta L }{ \eta+L} +2\alpha^2 \nu_2^2 \rho_1^{k+1} \right)\mathbb{E}[ \| x_k - x^*\|^2]  + \alpha^2 \rho_1^{k+1}  (\nu_1^2+2\nu_2^2 \|x^*\|^2 ) . \end{align}

 Next, we show  that $\mathbb{E}[ \| x_k - x^*\|^2] $ is uniformly bounded by some constants.

  {\bf Case 1:} If $2\alpha^2 \nu_2^2 \leq \rho_1-q  ,$  then $1-\tfrac{2 \alpha \eta L}{ \eta+L} +2\alpha^2 \nu_2^2 \rho_1^{k+1} \leq q +2\alpha^2 \nu_2^2 \rho_1 \triangleq \tilde{q}< q +2\alpha^2 \nu_2^2 \leq \rho_1.$  Hence, we conclude that  $ \tilde{q}<\rho_1. $
 Similarly to the derivation of \eqref{recursion2}, we obtain
 from \eqref{bd-ms} that for any $k\geq 0,$
 \begin{align*}
\mathbb{E}[\| x_{k}-x^*\|^2  ]
 & \leq   \tilde{q}\mathbb{E}[ \| x_k - x^*\|^2] +   \alpha^2 \rho_1^k (\nu_1^2+2\nu_2^2 \|x^*\|^2 )\leq\rho_1^{k} \left(\mathbb{E}[ \| x_0 - x^*\|^2]+\tfrac{\alpha^2(\nu_1^2+2\nu_2^2 \|x^*\|^2 )}{1- \tilde{q}/\rho_1}\right)
 \\&\leq \mathbb{E}[ \| x_0 - x^*\|^2]+\tfrac{\alpha^2(\nu_1^2+2\nu_2^2 \|x^*\|^2 )}{1- \tilde{q}/\rho_1} \triangleq c_1.
\end{align*}

  {\bf Case 2:} If $2\alpha^2 \nu_2^2 > \rho_1-q  ,$   we define $\tilde{k}\triangleq \left \lceil { \ln\left(\tfrac{  2 \alpha^2 \nu_2^2 }{ \rho_1-q } \right) \over \ln(\rho_1^{-1})} \right \rceil $. Then for any $k\geq \tilde{k},$
 $ \rho_1^{-(k+1)}>\tfrac{  2 \alpha^2 \nu_2^2 }{ \rho_1-q }  $ and hence
 $2\alpha^2 \nu_2^2 \rho_1^{k+1} < \rho_1-q  .$ Note from
 \eqref{bd-ms} that   for any $k\leq \tilde{k},$
  \begin{align*} \mathbb{E}[\| x_{k}-x^*\|^2  ]
      & \leq  \left(1 +2\alpha^2 \nu_2^2 \rho_1  \right) \mathbb{E}[ \| x_{k-1} - x^*\|^2]  + \alpha^2(\nu_1^2+2\nu_2^2 \|x^*\|^2 )
 \\& \leq  (1 +2\alpha^2 \nu_2^2 \rho_1 )^{\tilde{k}} \Big( \mathbb{E}[ \| x_0 - x^*\|^2] + \alpha^2(\nu_1^2+2\nu_2^2 \|x^*\|^2 ) { (1 +2\alpha^2 \nu_2^2 \rho_1 )^{\tilde{k}}-1\over(1 +2\alpha^2 \nu_2^2 \rho_1)  -1}\Big)\triangleq c_2. \end{align*}
 By defining $\hat{q}\triangleq q +2\alpha^2 \nu_2^2 \rho_1^{\tilde{k}+1},$
 we have  $1-{2 \alpha \eta L \over \eta+L} +2\alpha^2 \nu_2^2 \rho_1^{k+1}  \leq \hat{q}<\rho_1$  for any $k\geq \hat{k}.$
Then it follows from  \eqref{bd-ms}  that for any $k\geq \tilde{k},$
 \begin{align*}
 \mathbb{E}[\| x_{k+1}-x^*\|^2  ]
     & \leq    \hat{q}\mathbb{E}[ \| x_k - x^*\|^2] +   \alpha^2 \rho_1^{k+1} (\nu_1^2+2\nu_2^2 \|x^*\|^2 )
 \\&\leq\hat{q}^{k+1-\tilde{k}} \mathbb{E}[ \| x_{\tilde{k}} - x^*\|^2]+ \alpha^2 (\nu_1^2+2\nu_2^2 \|x^*\|^2 ) \sum_{t=0}^{k-\tilde{k}}\rho_1^{k+1-t} \hat{q}^t
  \leq c_2+\tfrac{\alpha^2(\nu_1^2+2\nu_2^2 \|x^*\|^2 )}{1- \hat{q}/\rho_1} .
\end{align*}
Thus, for the case  $2\alpha^2 \nu_2^2 > \rho_1-q  ,$ we have that
 $\mathbb{E}[\| x_k-x^*\|^2  ] \leq c_2+\tfrac{\alpha^2(\nu_1^2+2\nu_2^2 \|x^*\|^2 )}{1- \tilde{q}/\rho_1} $ for any $ k\geq 1$.

By combing  $2\alpha^2 \nu_2^2 \leq  \rho_1-q  $ ({\bf Case 1})  and    $2\alpha^2 \nu_2^2 > \rho_1-q   $ ({\bf Case 2}), we conclude that  there exits some constant $c_3>0$ such  that $\mathbb{E}[\| x_k-x^*\|^2  ] \leq  c_3 $ for any $ k\geq 1$. This combines
with \eqref{bd-ms} produces
 \begin{align*} \mathbb{E}[\| x_{k+1}-x^*\|^2  ]
 &\leq  \left(1-\tfrac{2 \alpha \eta L }{\eta+L}  \right)\mathbb{E}[ \| x_k - x^*\|^2]  + \alpha^2 \rho_1^{k+1}  (\underbrace{\nu_1^2+2\nu_2^2 \|x^*\|^2 +2\nu_2^2c_3}_{\triangleq c}) \\ &= q\mathbb{E}[ \| x_k - x^*\|^2]  + \alpha^2 c \rho_1^{k+1}   . \end{align*}
 The rest of the proof is the same as that of Proposition \ref{prp1}.
  \hfill $\blacksquare$

 Next, we provide similar statements for Algorithm~\ref{Alg_2},  and Algorithm \ref{Alg_3} (for both quadratic and non-quadratic
 objective functions),
 for which the proofs are omitted here but included   in the supplementary material  for completeness.
 \begin{proposition}[Rate and Oracle Complexity for  Algorithm  \ref{Alg_2}]  \label{prp2} Let Assumption \ref{ass-fun} hold. Consider Algorithm \ref{Alg_2},  where  $   \alpha \in ( 0, {1\over L}]$, $\gamma\triangleq \sqrt{\alpha \eta}$,   $\beta \triangleq {1 -\gamma \over 1+ \gamma},$ and $N_{k } \triangleq \lceil \rho_2^{-(k+1)} \rceil $  with   $\rho_2 \in (1-\gamma,1)$.  Then for any $k\geq 0 $,
\begin{align}\label{lem-fy2}
\mathbb{E}[f(y_{k }) ]-f^* \leq      \rho_2^k  \left( \tfrac{\eta+L}{2}    \mathbb{E}[\| x_0-x^*\|^2] +  \tfrac{\rho_2 \nu^2}{\rho_2-(1-\gamma)}\left( \alpha +\tfrac{(1-\gamma)\gamma}{2\eta}\right)  \right) ,
\end{align}
{and $\mathbb{E}[\|x_k-x^*\|^2] \leq c\rho_2^k  $ for some constant $c>0$. In addition,}
 the  iteration and  oracle  complexity for obtaining  an $\epsilon-$solution  are         $\mathcal{O}\left (  \sqrt{ \kappa } \ln(1/   \epsilon) \right)$ and   $\mathcal{O}\left ({ \sqrt{ \kappa }  / \epsilon} \right)$, respectively.
\end{proposition}

\begin{proposition}[Rate and Oracle Complexity for Algorithm  \ref{Alg_3} on Quadratic Functions]  \label{prp-alg3-Q}
  Suppose that Assumption \ref{ass-fun} holds and   $f(\cdot)$ is a quadratic function with   $  \nabla^2 f(x) \equiv \mathbf{H} $ for any $x\in \mathbb{R}^m.$    Consider  Algorithm  \ref{Alg_3},  where  $   \alpha \triangleq { 4 \over (\sqrt{\eta}+\sqrt{L})^2}$,
$\beta \triangleq   \left(  { \sqrt{\kappa} -1 \over \sqrt{\kappa} +1}\right)^2,$  and  $N_{k } \triangleq \lceil \rho_3^{-(k+1)} \rceil $  with   $\rho_3 \in \left(\beta,1\right)$.  Then for  any $\iota \in (0,\sqrt{\rho_3}-\sqrt{\beta})$, there exist a constant  $ c(\iota)$  such that
 \begin{align}\label{prp-hb-q}
 \mathbb{E} \left [ \left \| \begin{pmatrix}
  x_{k+1}-x^*
\\ x_k-x^* \end{pmatrix} \right\|^2   \right]\leq     c(\iota )^2 \left( 2 \mathbb{E}[ \| x_0 - x^*\|^2]+ {\alpha^2\nu^2\over 1-(\sqrt{\beta}+\iota )^{2}/\rho_3} \right)\rho_3^{k+1}.
\end{align}
  In addition, the  iteration and  oracle  complexity required for obtaining  an $\epsilon-$solution in a mean-squared sense are      respectively   $\mathcal{O}\left (  \sqrt{ \kappa } \ln(1/   \epsilon) \right)$ and   $\mathcal{O}\left ({ \sqrt{ \kappa }  / \epsilon} \right)$.
\end{proposition}

\begin{proposition}[Rate for Algorithm   \ref{Alg_3} on Non-Quadratic Functions]  \label{prp-alg3}
Let  Assumption \ref{ass-fun} hold. Consider Algorithm  \ref{Alg_3}, where     $   \beta \in (0,1)$, $\alpha\in (0, \tfrac{2( 1-\beta) }{ L+\eta})$, and $N_{k } \triangleq \lceil \rho_4^{-(k+1)} \rceil $  with   $\rho_4 \in \left(q,1\right)$. Here  $q=\max\{q_1,q_2\}$ with
$q_1$ and $q_2$ defined in \eqref{def-q1q2}.  Then  for any $k\geq 0.$
 \begin{align}\label{prp-hb}  \mathbb{E}[\|x_k-x^*\|^2]\leq \tfrac{1}{2 \hat{m} \eta}
 \left( \mathbb{E} \left [\|  x_{0}-x^*\|^2+ \hat{m} ( f(x_0)-f(x^*)) \right] + \tfrac{\alpha^2 \nu^2}{(1-\beta)^2  (1-q/\rho_4)}\right) \rho_4^{k+1}  ,
\end{align}
where $\hat{m} \triangleq \tfrac{2\alpha  }{ 1-\beta}\left(\tfrac{\beta-\alpha \eta }{ 1-\beta}+
 \tfrac{2\eta }{\eta+L} \right)$.  
\end{proposition}

\begin{remark}
From the  complexity  statements in Propositions \ref{prp1} and \ref{prp2}, we see that the  dependence on the condition  number is improved from $\kappa$  (in Algorithm \ref{Alg_1})
to $\sqrt{\kappa}$  by the accelerated Algorithm \ref{Alg_2}.
    The rate and   oracle  complexity for the heavy-ball method  (Algorithm \ref{Alg_3}) of quadratic function  are
    similar  to that  observed in  Algorithm \ref{Alg_2},  while Algorithm \ref{Alg_3}   for non-quadratic functions might not lead an acceleration  (in the sense of the dependence on $\kappa$ in the constant)  of Algorithm \ref{Alg_1} (see also Remark \ref{rem9}).  In addition, Propositions \ref{prp1}-\ref{prp-alg3-Q} imply that a smaller $\kappa$ leads to a smaller constant in the oracle complexity.
\end{remark}

\begin{remark}[Convergence  in Probability]\label{rem-con-pro}
Because  mean-squared convergence implies  convergence  in probability,    on the basis of   Propositions \ref{prp1}-\ref{prp-alg3}, we may conclude that the sequences $\{x_k\}$   and $\{y_k\}$  generated by  Algorithms   \ref{Alg_1}-\ref{Alg_3}  converge in probability to the optimal solution $ x^* $.
\end{remark}

 \begin{remark}\label{rem-geo}
We discuss the choices employed in Algorithm 2 for different settings of $\epsilon$ and $\kappa$   in the following table.  In particularly,  we observe that the sample size at the last iterate is {\em relatively} modest.
\medskip
\begin{center}
\begin{tabular}{c|c|c|c|c|c}
\hline
$\kappa$ & $\epsilon$  & $K$ & $\beta = \tfrac{\sqrt{\kappa}-1}{\sqrt{\kappa}+1} $ &  $\rho \in (\beta,1)$ &  $N_K = \lceil \rho^{-K}\rceil$  \\
 \hline
 $10$ & $10^{-3}$  &  $\sqrt{\kappa} \ln(1/\epsilon)\approx 22$  &  $\approx 0.52$  & $\rho = 0.55$  & $\approx  5.1e5 $ \\ \hline
 $ 10$ & $10^{-3}$  &  $\sqrt{\kappa} \ln(1/\epsilon)\approx 22$  &  $\approx 0.52$  & $\rho = 0.85$ & $\approx  36 $ \\
\hline
 {$ 10^3$} & $10^{-3}$  & {$\sqrt{\kappa}  \ln(1/\epsilon)\approx 219$}  & $\approx 0.939$ & $\rho = 0.97$ & $\approx   {789}$\\ \hline
 {$ 10^3$} & $10^{-4}$  &  {$\sqrt{\kappa}  \ln(1/\epsilon)\approx 292$}  & $\approx 0.939$ & $\rho = 0.97$ & $\approx {7.3e3}$\\ \hline
\end{tabular}
\end{center}
\medskip

It can be seen  that for  well-conditioned problems ($\kappa$ close to $1$), the sampling rate at the last iteration reaches  $ 10^5$   when $\rho$ is chosen closer to $\beta$ and $\epsilon=10^{-3}$.
 However, it is possible to  choose $\beta$ closer to $1$ to obtain far more reasonable
sampling rates. In more practical settings with higher  $\kappa$ (e.g., $\kappa = 10^3$),  it is seen that sampling
rates are  {relatively modest}. In addition, by leveraging the parallelizable structure,
current multi-core and multi-processing techniques can contend with the
challenge of computing gradient estimators with large  sample sizes.  Finally, it
is noticed  that most large-scale machine learning problems defined
on finite sample spaces with cardinalities of $10^6$ to $10^9$ or even more.
Application of variance-reduced gradient methods (such as SVRG) necessitate
taking full gradients intermittently, implying that computing a gradient
estimate with this sample-size is well within the reach of current
computational constraints.    \hfill $\blacksquare$
\end{remark}

\subsection{Preliminary Lemmas.}
  Before establishing CLTs for Algorithms \ref{Alg_1}-\ref{Alg_3},  we first introduce a preliminary   CLT on doubly-indexed random variables \cite[Lemma 3.3.1]{chen2006stochastic}. We state it as Lemma \ref{lem-CLT},
  whose proof is found in~\cite[Chapter 12]{hall2014martingale}.
  \begin{lemma}\label{lem-CLT}
  Let $\bxi_{kt},t= 1,\cdots, k$ be $m$-dimensional random vectors. Define
 \begin{align} \label{def-RS}
 & \mathbf{S}_{kt}  \triangleq \mathbb{E}[\bxi_{kt}\bxi_{kt}^T], \quad  \mathbf{R}_{kt}\triangleq \mathbb{E}[\bxi_{kt}\bxi_{kt}^T| \bxi_{k1},\cdots, \bxi_{k,t-1}],
  {\rm ~and ~} \mathbf{S}_k\triangleq \sum_{t=1}^k  \mathbf{S}_{kt} \\
& \mbox{\rm Assume that \ }   \mathbb{E}[\bxi_{kt}| \bxi_{k1},\cdots, \bxi_{k,t-1}]=0, ~ \sup_{k\geq 1} \sum_{t=1}^k \mathbb{E}[\| \bxi_{kt}\|^2]<\infty, \label{cond1}\\
   &  \lim_{k\to \infty}  \mathbf{S}_k =\mathbf{S},~ \lim_{k \to \infty} \sum_{t=1}^k \mathbb{E}[\|  \mathbf{S}_{kt}- \mathbf{R}_{kt}\|]=0 , \label{cond2}\\
 & {\rm and~ }    \lim_{k\to \infty}  \sum_{t=1}^k \mathbb{E}[\| \bxi_{kt}\|^2I_{[\|\bxi_{kt}\|\geq \epsilon]}]=0, \quad \forall \epsilon>0.\label{cond3}
  \end{align}
  Then   $ \sum_{t=1}^k \bxi_{kt}   \xlongrightarrow [k \rightarrow \infty]{d}   N(0,\mathbf{S}).$
  \end{lemma}

To establish  the CLTs for  Algorithms~\ref{Alg_1}--\ref{Alg_3},
   we further require  the following conditions.
 \begin{assumption}\label{ass-CLT}
 (i)  $ \nabla f:\mathbb{R}^m \to \mathbb{R}^m$  is differentiable at $x^*$ with Hessian   matrix $ \mathbf{H} \in \mathbb{R}^{m\times m}$, and
 \begin{align}\label{Ass2}
     \nabla f(x) =  \mathbf{H}(x-x^*) +\mathbf{D}(x)(x-x^*) ~  {\rm with~}   \|\mathbf{D}(x)\| \leq R_D \|x-x^*\|   {\rm ~for~ some~} R_D>0. 
 \end{align}

 (ii)  The noise  sequence $\{w_{k,N_k}\} $ further satisfies
 \begin{align}
&\lim_{k\to \infty} \Big( N_k \mathbb{E}\left[ w_{k,N_k}w_{k,N_k}^T | \mathcal{F}_{k } \right]\Big) =\lim_{k\to \infty} \Big( N_k \mathbb{E}\left[ w_{k,N_k}w_{k,N_k}^T  \right]\Big) =\mathbf{S}_0,\quad a.s., \label{Ass3}
 \\ &\label{Ass32} {\rm and~}
\lim_{r\to \infty} \sup_k \mathbb{E}\left[ \left \|\sqrt{N_{k}} w_{k,N_{k}} \right \|^2 I_{\left[\|\sqrt{N_{k}} w_{k,N_{k}}\|> r\right]}\right]=0,
 \end{align}
 where  $I_{[a>b]}=1$ if $a>b$, and $I_{[a>b]}=0$, otherwise.
  \end{assumption}

{
    \begin{remark}\label{new-rem} We now  provide two  easily understood sufficient conditions to guarantee that Assumption  \ref{ass-CLT}(i)  {holds}.
    {In (A), we assume that the gradient $\nabla f$ is twice differentiable with bounded second derivatives, i.e. $f$ is assumed to be thrice differentiable while in (B), we require that $f$ is twice differentiable with Lipschitz continuous Hessians. }Define {$F$ as} $F(x)\triangleq \left(
                         \begin{array}{c}
                           F_1(x) \\
                           \vdots \\
                            F_m(x) \\
                         \end{array}
                       \right)=\nabla f(x)$, where $F_j: \mathbb{R}^m  \to \mathbb{R}$ for each $j=1,\cdots,m$.\\

        \noindent (A)  Suppose   $F$ is twice continuously differentiable and that $ \nabla^2 F_j(x) $ is bounded in $x\in \mathbb{R}^m$  for $j = 1, \cdots, m$.
Then for any $j$th ($j=1,\cdots,m$) component of  $ F (x) $,
        by the  second-order mean-value theorem,  we know there exists $\tilde{x}^j \in \mathbb{R}^m$ on the line segment connecting $x$ and  $x^*$ such that
\begin{align}\label{rec-Fj}
F_j(x) &=F_j(x^*)+ \nabla F_j(x^*)^T(x-x^*)+ {1\over 2}(x-x^*)^T \nabla^2 F_j(\tilde{x}^j )  (x-x^*) \notag
    \\ & =F_j(x^*)+   \nabla F_j(x^*)  (x-x^*)  + \underbrace{\left( {1\over 2}\nabla^2 F_j(\tilde{x}^j )(x-x^*) \right)^T}_{\, \triangleq \, \mathbf{D}_j(x )} (x-x^*).
\end{align}
        Since $\|\nabla^2 F_j(x)\| \leq R_{D_j}$ for all $x$, we have that
        $\|\mathbf{D}_j(x)\| \leq R_{D_j} \|x-x^*\|$. Thus,   we conclude  that
        $   \mathbf{D}(x)  \triangleq \,  \begin{pmatrix}  \mathbf{D}_1(x) \\ \vdots \\ \mathbf{D}_n(x) \end{pmatrix} $ satisfies  $
            \left\| \mathbf{D}(x) \right\| \,= \sqrt{\sum_{j=1}^m\| \mathbf{D}_j(x) \|^2 } \leq \, R_D \|x-x^*\| $, where  $R_D=\sqrt{\sum_{j=1}^m R_{D_j}^2 }.$
Since $F(x^*)=0$ and $ \nabla F(x^*)=\mathbf{H}$ is a symmetric matrix, we derive \eqref{Ass2} by stacking \eqref{rec-Fj} for $j = 1, \cdots, m$. \\

 (B). Suppose $f$ is twice continuously differentiable with Lipschitz continuous Hessians with constant $L_H$. By the integral mean-value theorem for the gradient, we have that
        \begin{align*}
            \nabla f(x) & = \nabla f(x^*) + \nabla^2 f(x^*) (x-x^*) + \int_{0}^1 \left(\nabla^2 f(x^*+ \zeta (x-x^*))-\nabla^2 f(x^*)\right) (x-x^*) d\zeta \\
                    & =  \mathbf{H} (x-x^*) + \mathbf{D}(x)(x-x^*), \mbox{ where } \mathbf{D}(x) \triangleq \int_{0}^1 \left(\nabla^2 f(x^*+\ \zeta (x-x^*))-\nabla^2 f(x^*)\right)  d\zeta.
        \end{align*}
        We now derive a bound on $\mathbf{D}(x)$ as follows.
    \begin{align*}
        \left\| \mathbf{D}(x) \right\| & = \left\| \int_{0}^1 \left(\nabla^2 f(x^*+\zeta(x-x^*))-\nabla^2 f(x^*)\right)  d\zeta \right\| \\
                    & \leq  \, \max_{0 \leq \zeta \leq 1} \left\| \left(\nabla^2 f(x^*+\zeta(x-x^*))-\nabla^2 f(x^*)\right)\right\|
                     \leq L_H \|x-x^*\|.
    \end{align*}
    Again, we see that \eqref{Ass2} holds.
     \hfill $\blacksquare$

\end{remark}}

\begin{remark}\label{rem3} Since  the Lindeberg condition   condition   \eqref{Ass32} is  less easily verified,  we   provide a   sufficient condition for the ease of understanding
and verification in practice.  Suppose there exists a constant $\delta>0$ and   a finite value $b>0$ such that
 $\mathbb{E}\big[   \|\sqrt{N_k} w_{k,N_k}   \|^{2+\delta}\big]  \leq b $ for any $ k\geq 0.$  Therefore,
\begin{align*}
& \quad \mathbb{E}\left[ \left \|\sqrt{N_{k }} w_{k,N_k} \right \|^2 I_{\left[\|\sqrt{N_{k }} w_{k,N_k}\|> r\right]}\right] \\
 & \leq  \tfrac{1}{ r^{\delta}}\mathbb{E}\left[ \left \|\sqrt{N_{k }} w_{k,N_k} \right \|^{2+\delta} I_{\left[\|\sqrt{N_k } w_{k,N_k}\|> r\right]}\right]
 \leq  \tfrac{1}{ r^{\delta}}\mathbb{E}\left[ \left \|\sqrt{N_{k }} w_{k,N_k} \right \|^{2+\delta}  \right] \leq   \tfrac{b}{ r^{\delta}}.
\end{align*}
 This implies that $\sup_k \mathbb{E}\left[ \left \|\sqrt{N_{k }} w_{k,N_k} \right \|^2 I_{ [\|\sqrt{N_{k }} w_{k,N_k}\|> r ]}\right] \leq  \tfrac{b}{ r^{\delta}} $.
 Hence \eqref{Ass32} {\rm~holds}. \hfill $\blacksquare$
\end{remark}

In the following, we  establish  the  central limit theorem of  a    linear  recursion, for which  the proof is provided  in Appendix \ref{app-lem-CLT2}. This result will be applied in the proof of Theorems \ref{thm-CLT1}-\ref{thm-CLT1-Alg3-NQ}.

  \begin{lemma}\label{lem-CLT2}  Suppose  that  $ \mathbf{P}$ is a square matrix   with spectral radius,  denoted by $ \rho(  \mathbf{P})$,  strictly smaller than 1 (i.e., $\rho(  \mathbf{P}) <1$),    $N_k=\lceil \rho^{-(k+1)} \rceil$ with $\rho\in(0,1)$ for any $ k\geq 0, $ and that  $\{ w_{k,N_k} \} $ satisfies Assumption \ref{ass-CLT}(ii).  Let  $\{e_k\}$ be generated by
\begin{align}\label{recursion-e}
e_{k+1}& =    \mathbf{P}  e_k -\alpha \rho^{-(k+1)/2} \mathbf{G}  w_{k,N_k}+\zeta_{k+1},
\end{align}
      where $\mathbb{E}[\| e_0\|^2]<\infty$ and $\zeta_k    \xlongrightarrow [k \rightarrow \infty]{P} 0 $. Then
 \begin{align*}
  &\alpha^{-1}   e_k\xlongrightarrow [k \rightarrow \infty]{d}   N(0,\bm{\Sigma})
\quad {\rm with} \quad   \bm{\Sigma}  \triangleq \sum_{t=0}^{\infty}  \mathbf{P} ^{t}  \mathbf{G} \mathbf{S}_0    \mathbf{G}^T \left( \mathbf{P}^t \right)^T  .\end{align*}
\end{lemma}

  \subsection{Central Limit Theorems.}\label{sec:geo_CLT}
Based on  Proposition \ref{prp1} and  Lemma \ref{lem-CLT2}, using Assumption \ref{ass-CLT}, we present  the central limit theorem for  Algorithm \ref{Alg_1} with geometrically increasing  batch-sizes.

 \begin{theorem}[CLT of Algorithm \ref{Alg_1} with Geometrically increasing $N_k$]\label{thm-CLT1} Let Assumptions \ref{ass-fun}  and \ref{ass-CLT} hold. Consider Algorithm \ref{Alg_1},  where  $   \alpha \in (0, {2\over \eta+L}]$.    Set $q\triangleq 1-{2 \alpha \eta L \over \eta+L} $, $N_{k } \triangleq \lceil \rho_1^{-(k+1)} \rceil $ with $\rho_1 \in (q,1)$,
 and  $\mathbf{P}_1  \triangleq  \rho_1^{-1/2}(\mathbf{I}_m-\alpha  \mathbf{H})$.  Then  $\rho(  \mathbf{P}_1) <1$ and
 \begin{align}  \label{Geo-CLT}
  &\alpha^{-1} \rho_1^{-k/2}( x_{k}-x^*)  \xlongrightarrow [k \rightarrow \infty]{d}   N(0,\bm{\Sigma}_1),
~{\rm where}~  \bm{\Sigma}_1 \triangleq \sum_{t=0}^{\infty}  \mathbf{P}_1 ^{t} \mathbf{S}_0    \mathbf{P}_1 ^{t}  .\end{align}
 \end{theorem}
 {\bf Proof. }   By  using \eqref{Ass2}, we can rewrite \eqref{ref_alg1}  as
 \begin{align}\label{recursion3}
x_{k+1}-x^*& =x_k-x^* -\alpha  \mathbf{H}(x_k-x^*)-\alpha  (\nabla f(x_k)- \mathbf{H}(x_k-x^*)+w_{k,N_k}) \notag \\
& = (\mathbf{I}_m-\alpha  \mathbf{H})(x_k-x^*) -\alpha  ( \mathbf{D}(x_k)(x_k-x^*) +w_{k,N_k}).
\end{align}
 Define $e_{k} \triangleq \rho_1^{-k/2}(x_{k}-x^*)$. Then by multiplying  both sides of  \eqref{recursion3} by  $\rho_1^{-(k+1)/2}$   and by using the definition $  \mathbf{P}_1  \triangleq  \rho_1^{-1/2}(\mathbf{I}_m-\alpha  \mathbf{H}) $, we  obtain
 \begin{align}\label{recursion4}
e_{k+1}& =   \rho_1^{-1/2}(\mathbf{I}_m-\alpha  \mathbf{H})\rho_1^{-k/2}(x_{k}-x^*) -\alpha \rho_1^{-(k+1)/2} (\mathbf{D} (x_k) (x_k-x^*)+w_{k,N_k}) \notag
\\& =    \mathbf{P}_1   e_k -\alpha \rho_1^{-(k+1)/2} w_{k,N_k}+\zeta_{k+1} ~{\rm with~}\zeta_{k+1} \triangleq - \alpha \rho_1^{-1/2} \mathbf{D}(x_k)e_{k} .
\end{align}

In the following, we prove that $ \zeta_{k }   \xlongrightarrow [k \rightarrow \infty]{P} 0 .$ From \eqref{recursion2},  it follows  that for any $k\geq 0,$
\begin{align}\label{bd-ek2}
 {\operatorname {var} (e_k) } \leq \mathbb{E}[\|e_{k}\|^2]=\mathbb{E}[\| x_{k}-x^*\|^2 /\rho_1^k ]  \leq \mathbb{E}[\|x_0-x^*\|^2]+{\alpha^2 \nu^2 \over  1-q/\rho_1} \triangleq v_e^2 {\rm~with~}v_e>0 .\end{align}
Chebyshev's inequality asserts that if $X$ is a random variable with mean $\mu$  and variance $\sigma^2$, then for any real number $h>0$,
 \begin{align}\label{bd-OP} \mathbb{P}( \|X- \mu \|   \leq h \sigma) \geq 1-h^{-2}.\end{align}
     By setting $h=   \chi^{-1/2}$  for any $\chi\in (0,1)$ and applying \eqref{bd-OP} \red{where  $X = e_k$}, we have
\begin{align}\label{pro-A1} \mathbb{P}\big( \|e_k-\mathbb{E}[e_{k} ]\|   \leq \chi^{-1/2} {\operatorname {var} (e_k) }\big) \geq 1-\chi.\end{align}
Define the events
$\mathcal{A}_1\triangleq \big \{ \|e_k-\mathbb{E}[e_{k} ]\|   \leq \chi^{-1/2} {\operatorname {var} (e_k)} \big\} $,
$\mathcal{A}_2\triangleq \big \{\|e_k-\mathbb{E}[e_{k} ]\|   \leq \chi^{-1/2} v_e^2 \big\} ,$
$\mathcal{A}_3  \triangleq \big\{\|e_k\| \leq \| \mathbb{E}[e_{k} ]\| + \chi^{-1/2} v_e^2 \big\} ,   $ $\mathcal{A}_4 \triangleq
  \big\{\|e_k\| \leq \mathbb{E}[\| e_{k} \|] + \chi^{-1/2} v_e^2 \big\} $, and $\mathcal{A}_5 \triangleq
  \big\{\|e_k\| \leq v_e + \chi^{-1/2} v_e^2 \big\} $. Note by \eqref{bd-ek2}  that $\mathcal{A}_1\subseteq \mathcal{A}_2.$
  We observe that $\mathcal{A}_2\subseteq \mathcal{A}_3 $ by the   inequality $\|x_1-x_2\| \geq \|x_1\|-\| x_2\|$.
  Since $\|\bullet\|$ is convex in $x,$  by the Jensen's inequality we have  $\|\mathbb{E}[X]\| \leq\mathbb{E}[\|X\|]$, and   hence  $\mathcal{A}_3\subseteq \mathcal{A}_4. $
Since $\|\bullet\|^2$ is a convex function, by using \eqref{bd-ek2}  and Jensen's inequality,
we have  $\mathbb{E}[\|e_{k} \|]\leq \sqrt{\mathbb{E}[\|e_{k} \|^2]} \leq  v_e$ and  hence  $\mathcal{A}_4\subseteq \mathcal{A}_5. $
Thus, we have $\mathcal{A}_1\subseteq \mathcal{A}_5.$ This together with \eqref{pro-A1}  implies  that for any $k\geq 0,$
  \[\mathbb{P}(\|e_k\| \leq v_e + \chi^{-1/2} v_e^2) =\mathbb{P}(\mathcal{A}_5) \geq  \mathbb{P}(\mathcal{A}_1)= \mathbb{P}\big( \|e_k-\mathbb{E}[e_{k} ]\|   \leq \chi^{-1/2} {\operatorname {var} (e_k) }\big) \geq 1-\chi.\]
That is to say that for any $\chi\in (0,1)$,
\begin{align}\label{bd-absek} \mathbb{P}(\|e_k\| > v_e + \chi^{-1/2} v_e^2)  <\chi ,~\forall k\geq 0. \end{align}
 Therefore, we conclude that  $ \| e_k \|$ is   bounded  in probability (i.e., $\| e_k\| =O_P(1)$).

Note by Proposition \ref{prp1} and the Markov's inequality that   $ \| x_k -x^*\|   \xlongrightarrow [k \rightarrow \infty]{P}  0$, i.e., $\| x_k -x^*\| =o_P(1)$.
 Then by invoking the bound that  $ \| \mathbf{D}(x) \| \leq R_D \| x- x^*\|  $,
we conclude that   $ \| \mathbf{D}(x_k)\|=  o_P(1)  .$ 
Recall  from  \cite[p.12]{van2000asymptotic} that  symbols $o_P(\cdot)$ and $O_P(\cdot)$  satisfy
  \begin{align} \label{bigO} O_P(1) o_P(1)=o_P(1) . \end{align}
      This  together with $ \| \mathbf{D}(x_k)\|=  o_P(1)   $ and $\| e_k\| =O_P(1)$  implies $  \| \zeta_{k+1} \| \leq \alpha \rho_1^{-1/2} \|  \mathbf{D}(x_k) \| \| e_k\|=  O_P(1)  o_P(1)  
 \overset{ \eqref{bigO}}{ =}  o_P(1)  .$
 Thus, we conclude that $ \zeta_{k }   \xlongrightarrow [k \rightarrow \infty]{P} 0 .$

Since $\mathbf{H}$ is the Hessian matrix of    $f(x)$ at $x=x^*,$ by Assumptions \ref{ass-fun}(i) and  \ref{ass-fun}(ii),
we conclude  that $\mathbf{H}$ has eigenvalues $\lambda_1, \cdots, \lambda_m$ that satisfy $0<\eta \leq \lambda_m\leq \lambda_{m-1} \leq \dots\leq \lambda_2\leq \lambda_1\leq L.$  Then  $\|  \mathbf{I}_m-\alpha  \mathbf{H}  \|_2 \leq   \max\{ |1-\alpha \eta|,  |1-\alpha L| \} .$ We first show that  {$ |1-\alpha \eta| \geq  |1-\alpha L| $}. It can be seen that  $|1-\alpha \eta| \geq  |1-\alpha L|$  when $\alpha\in (0,1/L]$.
While for any $\alpha\in [1/L,{2\over \eta +L}]$, $| 1-\alpha \eta|=1-\alpha \eta,  |1-\alpha L|=\alpha L -1 $, and hence $|1-\alpha \eta| - |1-\alpha L| =2-\alpha (\eta+L) \geq 0$. Then for any     $\alpha \in (0,{2\over \eta +L}] $,
$$\|  \mathbf{I}_m-\alpha  \mathbf{H}  \|  \leq   | 1-\alpha \eta | =1-\alpha \eta \leq \sqrt{1-\tfrac{2 \alpha \eta L}{\eta+L}}  .
  $$
 The last inequality holds because  $ \alpha \leq {2\over \eta+L}\Rightarrow {2L \over \eta+L}+\alpha\eta-2\leq 0\Rightarrow {2\alpha\eta L \over \eta+L}+(\alpha\eta)^2-2\alpha\eta\leq 0 \Rightarrow (1-\alpha \eta)^2 \leq  1-{2 \alpha \eta L \over \eta+L}.  $  Since $q=1-{2 \alpha \eta L \over \eta+L}$ and   $\rho_1 \in (q,1)$,  we have that
\begin{align}\label{bd-P}
\|  \mathbf{P}_1 \| =\rho_1^{-1/2}\| (\mathbf{I}_m-\alpha  \mathbf{H}) \|  \leq     (q/\rho_1)^{1/2}<1.
\end{align}
 Since the spectral radius  of a symmetric matrix $ \mathbf{P}_1$ equals its two norm,  i.e., $ \rho(\mathbf{P}_1)=\| \mathbf{P}_1\|<1,$
by  invoking  Lemma \ref{lem-CLT2}  and  setting $
\mathbf{G}=\mathbf{I}_m$ and $ \mathbf{P}= \mathbf{P}_1$
(symmetric),  we conclude that the sequence $\alpha^{-1}e_k$
generated by  \eqref{recursion4}  converges in distribution to a
normally distributed random variable with  zero mean and  covariance
$\bm{\Sigma}_1$ defined as in \eqref{Geo-CLT}. Therefore, the result follows by recalling that $ e_k=\rho_1^{-k/2}( x_{k}-x^*)$.   \hfill $\blacksquare$

 Based on Lemma 4,  by employing similar proof arguments as in Theorem \ref{thm-CLT1}, we conclude that the central limit result established in Theorem \ref{thm-CLT1}    holds as well under Assumptions \ref{ass-fun}(i)-(ii),
 Assumption \ref{ass-CLT},  and the state-dependent  noise condition \eqref{ass-bd-noise2}.

 \begin{remark}
Suppose we set  $\alpha=\tfrac{2}{\eta+L}$ and $\rho_1 \triangleq \left(\tfrac{\kappa}{\kappa+1}\right)^2 $in  Algorithm \ref{Alg_1}.    If $V \sim  N(0,\mathbf{I}_m)$, it follows from \eqref{Geo-CLT}  that
 $\tfrac{ \eta+L}{2} \left(  \tfrac{ \kappa+1}{\kappa }\right)^k( x_{k}-x^*)  \xlongrightarrow [k \rightarrow \infty]{d}   N(0,\bm{\Sigma}_1) $. Namely,
\begin{align*}  x_k  \overset {\tiny D}{\approx } x^* +\tfrac{2}{\eta+L }\left(1- \tfrac{ 1}{\kappa+1 }\right)^k\bm{\Sigma}_1^{1/2} V \mbox{~for large~}k.
\end{align*}
This result implies that the sequence $ \{x_{k}\} $ converges in distribution to the optimal solution $x^* $ at    rate 	   $(\tfrac{\kappa}{\kappa+1 })^k $, and  $\{x_k\}$ is asymptotically normally distributed  for large $k$.   This provides the possibility of assessing   confidence regions of the estimate  from the normal distribution.   In addition, the estimation error for large $k$ depends on the     structure of the studied problem  (including $\eta, L$, and the Hessian matrix $\mathbf{H}$) and probability distribution  of the gradient   noise,   measured through  the coupling   matrix $\bm{\Sigma}_1.$  Thus, the problem's difficulty is largely characterized by  the covariance matrix $\bm{\Sigma}_1.$
\end{remark}

Based on the CLT established in Theorem \ref{thm-CLT1}, we
proceed to  use the Delta method~\cite{shapiro2009lectures} to    derive the
asymptotic distribution  of the gradient
and  the objective function based on  the estimation  sequence  $\{x(k)\}$.

\begin{corollary}\label{cor1}
 Consider Algorithm \ref{Alg_1} and  suppose  all    conditions of Theorem \ref{thm-CLT1} hold. Then

\noindent (i) $\alpha^{-1} \rho_1^{-k/2}  \nabla f (x_{k})    \xlongrightarrow [k \rightarrow \infty]{d}   N(0,  \mathbf{H} \bm{\Sigma}_1 \mathbf{H} ).$

\noindent  (ii) $\alpha^{-2} \rho_1^{-k}  (f (x_{k}) -f(x^*))   \xlongrightarrow [k \rightarrow \infty]{d}   {1\over 2 }N(0,  \bm{\Sigma}_1  )^T \mathbf{H}N(0,  \bm{\Sigma}_1  ).$
\end{corollary}
{\bf Proof.} (i)  Note  by Assumption \ref{ass-CLT}(i) that  the gradient
mapping  $ \nabla f(x):\mathbb{R}^m \to \mathbb{R}^m$  is differentiable at
$x^*$ with Hessian   matrix $ \mathbf{H}$.  By  using \eqref{Geo-CLT} and the
Delta theorem~\cite[Eqn. (7.182)]{shapiro2009lectures}, we have that
$$\alpha^{-1} \rho_1^{-k/2} \big( \nabla f (x_{k})- \nabla f (x^*)\big)
\xlongrightarrow [k \rightarrow \infty]{d}   \mathbf{H}  N(0,  \bm{\Sigma}_1
).  $$
Then  by  the fact that $\nabla f (x^*)=0 $  and $ \mathbf{H} $ is symmetric,  the assertion (i) holds.

 \noindent (ii) By  using \eqref{Geo-CLT}, $\nabla f (x^*)=0,$  and the second-order Delta theorem~\cite[Theorem 7.70]{shapiro2009lectures},   one obtains
$$\alpha^{-2} \rho_1^{-k}  \big(f (x_{k}) -f(x^*)\big)   \xlongrightarrow [k \rightarrow \infty]{d}   {1\over 2 } f_{u}^{''}(x^*)  {~\rm with ~}u=N(0,  \bm{\Sigma}_1  ), $$
  where $f_{u}^{''}(x^*) $  denotes  the second order directional  derivative at   $x^*$ along the direction $u$.   Since   $f:\mathbb{R}^m \to \mathbb{R} $ is twice continuously differentiable, $f_{u}^{''}(x^*)=u^T \mathbf{H} u$.
Then   result (ii) follows from $u=N(0,  \bm{\Sigma}_1  )$.
\hfill $\blacksquare$

 Akin  to Theorem \ref{thm-CLT1},  based on Proposition \ref{prp2}  and
Lemma \ref{lem-CLT2},   we    can also establish  the CLT for
Algorithm \ref{Alg_2} with geometrically increasing batch-sizes.

 \begin{theorem}[CLT for Algorithm \ref{Alg_2} with  Geometrically Increasing $N_k$]\label{thm-CLT1-Alg2} Suppose Assumptions \ref{ass-fun}  and \ref{ass-CLT} hold. Consider Algorithm \ref{Alg_2}, where  $   \alpha \in (0, {1\over  L}]$,  $\gamma \triangleq \sqrt{ \alpha \eta } $,   $\beta \triangleq {1 -\gamma \over 1+ \gamma}$, and  $N_{k }\triangleq \lceil \rho_2^{-(k+1)} \rceil $  with   $\rho_2 \in (1-\gamma,1)$. Then
 \begin{align*}
  &\alpha^{-1} \rho_2^{-k/2} \begin{pmatrix}
y_{k }-x^* \\
y_{k-1}-x^* \end{pmatrix} \xlongrightarrow [k \rightarrow \infty]{d}   N(0,\bm{\Sigma}_2) {~\rm with~}\bm{\Sigma}_2 \triangleq \sum_{t=0}^{\infty}  \mathbf{P}_2 ^{t} \begin{pmatrix} \mathbf{S}_0       & \mathbf{0}_m   \\ \mathbf{0}_m &\mathbf{0}_m\end{pmatrix}  \left( \mathbf{P}_2 ^{t} \right)^T ,  \end{align*}
where $ \mathbf{P}_2 \triangleq  \rho_2^{-1/2}\begin{pmatrix}
(1+\beta) (\mathbf{I}_m-\alpha  \mathbf{H})  &~ -\beta (\mathbf{I}_m-\alpha  \mathbf{H})
\\  \mathbf{I}_m   &~\mathbf{0}_m \end{pmatrix}$ and $\rho( \mathbf{P}_2)<1.$
 \end{theorem}
{\bf Proof.} By using  \eqref{Ass2},
  we can rewrite \eqref{Alg21}    as
   \begin{align}\label{Alg21-y}
y_{k+1}-x^*&  = (\mathbf{I}_m-\alpha  \mathbf{H})(x_k-x^*) -\alpha  ( \mathbf{D}(x_k)(x_k-x^*) +w_{k,N_k}).
\end{align}
Note by  \eqref{Alg22}  that $ x_{k }-x^* = (1+\beta) (y_{k }-x^*)-\beta (y_{k-1}-x^*).$
This together with \eqref{Alg21-y} produces
  \begin{align*}
y_{k+1}-x^*& = (\mathbf{I}_m-\alpha  \mathbf{H})( (1+\beta) (y_{k }-x^*)-\beta (y_{k-1}-x^*)) -\alpha  ( \mathbf{D}(x_k)(x_k-x^*) +w_{k,N_k})
\\&=  (1+\beta) (\mathbf{I}_m-\alpha  \mathbf{H})  (y_{k }-x^*)-\beta (\mathbf{I}_m-\alpha  \mathbf{H}) (y_{k-1}-x^*) -\alpha  ( \mathbf{D}(x_k)(x_k-x^*) +w_{k,N_k})    .
\end{align*}
Define $z_{k+1} \triangleq     \begin{pmatrix}
y_{k+1}-x^* \\
y_{k}-x^* \end{pmatrix} $ and $\mathbf{H}_2 \triangleq  \begin{pmatrix}
(1+\beta) (\mathbf{I}_m-\alpha  \mathbf{H})  & ~-\beta (\mathbf{I}_m-\alpha  \mathbf{H})
\\  \mathbf{I}_m   &  ~\mathbf{0}_m \end{pmatrix} $. Then based on  the above equation, we obtain the following recursion.
\begin{align}  \label{re-z}
z_{k+1}&  =\mathbf{H}_2 z_k   -  \alpha \begin{pmatrix}   \mathbf{D}(x_k)(x_k-x^*)  +w_{k,N_k}    \\ 0\end{pmatrix}  .
\end{align}

The eigenvalue decomposition of $\mathbf{H}$ is given by $\mathbf{H}=\mathbf{U}\bm{\Lambda} \mathbf{U}^T$,  where $\bm{\Lambda} \triangleq {\rm diag}\{\lambda_1, \lambda_2,\cdots, \lambda_m\}$ and $\mathbf{U}$ is orthogonal. This allows us to rewrite  {$\mathbf{H}_2$ as}
\begin{align*} \mathbf{H}_2 = \begin{pmatrix}\mathbf{U} &  {\bf 0}_m
\\   {\bf 0}_m   & \mathbf{U}
  \end{pmatrix}    \begin{pmatrix}
(1+\beta) (\mathbf{I}_m-\alpha   \bm{\Lambda})  & ~-\beta (\mathbf{I}_m-\alpha  \bm{\Lambda})
\\  \mathbf{I}_m   &  ~\blue{{\bf 0}_m } \end{pmatrix}  \begin{pmatrix}\mathbf{U} & \blue{{\bf 0}_m }
\\  \blue{{\bf 0}_m }& \mathbf{U}
  \end{pmatrix}^T .  \end{align*}
     Then the eigenvalues of  the matrix  $ \mathbf{H}_2  $ are given by  the roots of the equation  $det\big(vI_{2m}- \mathbf{H}_2 \big)=0$, i.e.,
\[ det \begin{pmatrix}
v \mathbf{I}_m-(1+ \beta  )(\mathbf{I}_m-\alpha   \bm{\Lambda}) & ~\beta (\mathbf{I}_m-\alpha  \bm{\Lambda})
\\-  \mathbf{I}_m &  v \mathbf{I}_m   \end{pmatrix}=0. \]
By the property $det \begin{pmatrix}
A& B
\\C &  D  \end{pmatrix}=det(AD-BC)$ when  the blocks $A,B,C,D$ are square matrices of the same size and $CD=DC$ \cite{silvester2000determinants},
we have\[ det \begin{pmatrix}
v\big(v \mathbf{I}_m-(1+ \beta  )(\mathbf{I}_m-\alpha   \bm{\Lambda}) \big) + \beta  (\mathbf{I}_m-\alpha  \bm{\Lambda})   \end{pmatrix}=0. \]
Since the matrix in the above determinant  is a diagonal matrix,  it can be described by the following  characteristic  equations for $i = 1, \cdots, m$.
\begin{align}\label{find-root}v \big(v -(1+\beta )(1-\alpha   \lambda_i)\big) +\beta (1-\alpha   \lambda_i)=v^2 -(1+\beta )(1-\alpha   \lambda_i)v +\beta (1-\alpha   \lambda_i)=0. \end{align}

By   $\gamma^2=\alpha \eta$ and $\beta = {1 -\gamma \over 1+ \gamma},$  we will  show that  the discriminant of the above  quadratic equation is nonpositive for each $i=1,\cdots,m$.
\begin{align*}
\Delta_i & =(1+\beta)^2 (1-\alpha    \lambda_i)^2 -4  \beta (1-\alpha   \lambda_i )=4(1-\alpha   \lambda_i ) \left(  {(1-\alpha   \lambda_i )  \over (1+ \gamma)^2}-{1 -\gamma \over 1+ \gamma} \right )
\\& ={ 4(1-\alpha   \lambda_i ) \over (1+ \gamma)^2}  \left(   \gamma^2-\alpha   \lambda_i \right )=
{ 4\alpha (1-\alpha   \lambda_i ) \over (1+ \gamma)^2}  \left(  \eta-  \lambda_i \right ) \leq 0,
\end{align*}
where the last inequality holds by the fact that for any $   i=1, \dots,m, $  $1-\alpha   \lambda_i  \geq  1-  \lambda_i /L\geq 0$ and  $   \eta -\lambda_i \leq 0. $  Consequently, \eqref{find-root} has two complex roots ${ (1+\beta) (1-\alpha    \lambda_i) \pm \mathbf{i} \sqrt{-\Delta_i } \over 2}$,  where $\mathbf{i}$ denotes the imaginary part.  Thus,  the magnitude  of the roots is  $
\sqrt{ \left(\tfrac{ (1+\beta) (1-\alpha    \lambda_i)}{2} \right)^2 -\tfrac{\Delta_i}{4}}= \sqrt{\beta (1-\alpha   \lambda_i ) }  .$  Since $\gamma^2=\alpha \eta$, $\beta ={1 -\gamma \over 1+ \gamma},$  $\rho_2 \in (1-\gamma,1)$, and $\lambda_i  \in [\eta, L]$,  we obtain that
\begin{align} \label{sp-h}
\rho(\mathbf{H}_2)& =\max_i  \sqrt{\beta (1-\alpha   \lambda_i ) } \leq \sqrt{\beta (1-\alpha  \eta) }
=\sqrt{\tfrac{(1 -\gamma)(1-\gamma^2)}{1+ \gamma}  }  =1-\gamma <\rho_2 <\rho_2^{1/2}.
\end{align}
Define $\varepsilon_k \triangleq \rho_2^{-k/2} z_k.$  By multiplying  both sides of \eqref{re-z} by $\rho_2^{-(k+1)/2}$,    we obtain the following recursion {with $\varsigma_{k+1} \triangleq -\alpha  \rho_2^{-(k+1)/2}  \mathbf{D}(x_k)(x_k-x^*) :$}
\begin{equation*}
\begin{split}
\varepsilon_{k+1}&  =\rho_2^{-1/2}\mathbf{H}_2 \varepsilon_k   -  \alpha  \rho_2^{-(k+1)/2} \begin{pmatrix}   \mathbf{D}(x_k)(x_k-x^*) +w_{k,N_k}    \\ 0\end{pmatrix}
  = \mathbf{P}_2 \varepsilon_k   -  \alpha  \rho_2^{-(k+1)/2} \begin{pmatrix}    w_{k,N_k}    \\ 0\end{pmatrix}+\begin{pmatrix}   \varsigma_{k+1}    \\ 0\end{pmatrix} .
\end{split}
\end{equation*}

We proceed to show that $ \varsigma_{k}  \xlongrightarrow [k \rightarrow \infty]{P} 0 $.
 By Proposition \ref{prp2}, we see that $\mathbb{E}[\rho_2^{-k}\|x_{k}-x^*\|^2] \leq c $  for some constant $c>0 $   and
 $\| x_k -x^*\| =o_P(1)$ by the   Markov's inequality. Then by  invoking the bound that  $ \| \mathbf{D}(x) \| \leq R_D \| x- x^*\|  $,  we achieve   $\| \mathbf{D}(x_k)\|= o_P(1).$
Similarly to the procedures for proving \eqref{bd-absek},  we can  also  show  that
 $\| \rho_2^{-k/2}(x_k-x^*)\|=O_P(1)$.
Therefore, we conclude that
  \begin{align*}
  & \| \varsigma_{k+1}\| \leq \alpha  \rho_2^{-1/2} \|\mathbf{D}(x_k)\|  \| \rho_2^{-k/2}(x_k-x^*)\|= o_P(1) O_P(1) \overset{ \eqref{bigO} }{=}  o_P(1)
  \Rightarrow  \varsigma_{k+1}  \xlongrightarrow [k \rightarrow \infty]{P} 0  .
  \end{align*}

 Because  $\rho(\mathbf{P}_2)  = \rho_2^{-1/2} \rho( \mathbf{H}_2)<1$   by \eqref{sp-h},   by
  {setting $\mathbf{P}=\mathbf{P}_2,~ \mathbf{G}= \begin{pmatrix} \mathbf{I}_m    \\ \mathbf{0}_m\end{pmatrix} $ and} using Lemma \ref{lem-CLT2}, we obtain that $ \alpha^{-1}  \varepsilon_k \xlongrightarrow [k \rightarrow \infty]{d}   N(0,\bm{\Sigma}_2) $
with $\bm{\Sigma}_2  \triangleq \lim_{k\to \infty}\sum_{t=0}^{k}  \mathbf{P}_2 ^{t}  \mathbf{G} \mathbf{S}_0    \mathbf{G}^T \left( \mathbf{P}_2^t \right)^T   .$
Then by $\mathbf{G} \mathbf{S}_0    \mathbf{G}^T= \begin{pmatrix} \mathbf{S}_0      & \mathbf{0}_m   \\ \mathbf{0}_m &\mathbf{0}_m\end{pmatrix}$ and $\varepsilon_k=\rho_2^{-k/2} z_k=\rho_2^{-k/2}  \begin{pmatrix}
y_{k }-x^* \\ y_{k-1}-x^* \end{pmatrix}$, we prove the result.
\hfill  $\blacksquare$

\begin{remark}
 {By setting   $\alpha={1\over L}$  and   $\rho_2
=1-{1\over a\sqrt{\kappa}} $ for some $a>1$ in Algorithm \ref{Alg_2}, we have that
$1-\gamma=1-{1\over \sqrt{\kappa}}<\rho_2,$ and    Theorem \ref{thm-CLT1-Alg2} holds.
Thus, we obtain  $L  \left(1-{1\over a\sqrt{\kappa}}\right)^{-k/2} \begin{pmatrix}
y_{k }-x^* \\
y_{k-1}-x^* \end{pmatrix} \xlongrightarrow [k \rightarrow \infty]{d}   N(0,\bm{\Sigma}_2) $.
Hence the following holds with $V \sim  N(0,\mathbf{I}_{2m}).$}
  \begin{align*}
  & \begin{pmatrix}
y_{k }  \\
y_{k-1} \end{pmatrix} \overset {\tiny D}{\approx } \begin{pmatrix}
 x^* \\
 x^* \end{pmatrix} +{1\over L }\left( 1-{1\over  a\sqrt{\kappa}}\right)^{k/2}\bm{\Sigma}_2^{1/2} V \mbox{~for large~}k
\end{align*}
\end{remark}

 In a similar fashion,   we establish  the central limit theorem for the stochastic heavy ball method (Algorithm \ref{Alg_3})  \red{on quadratic objective functions} with geometrically increasing batch-sizes  based on Proposition~\ref{prp-alg3-Q}  and   Lemma \ref{lem-CLT2}.

 \begin{theorem}[CLT  for  Alg. \ref{Alg_3} with Geometrically Increasing $N_k$ on Quadratic Functions]\label{thm-CLT1-Alg3}
Suppose that  Assumptions \ref{ass-fun}  and \ref{ass-CLT}(ii) hold, and   $f$ is a quadratic function with   $  \nabla^2 f(x) \equiv \mathbf{H} $ for any $x\in \mathbb{R}^m.$
Consider  Algorithm \ref{Alg_3}, where  $   \alpha \triangleq { 4 \over (\sqrt{\eta}+\sqrt{L})^2}$,  $\beta \triangleq   \left(1-{2 \over \sqrt{\kappa} +1}\right)^2,$  and   $N_{k }\triangleq \lceil \rho_3^{- (k+1)} \rceil $  with   $\rho_3 \in \left( \beta,1\right)$ for any $k\geq 0.$ Then
 \begin{align*}
  &\alpha^{-1} \rho_3^{-k/2} \begin{pmatrix}
x_{k }-x^* \\
x_{k-1}-x^* \end{pmatrix} \xlongrightarrow [k \rightarrow \infty]{d}   N(0,\bm{\Sigma}_3) {~\rm with~}\bm{\Sigma}_3 \triangleq  \sum_{t=0}^{\infty}  \mathbf{P}_3 ^{t} \begin{pmatrix} \mathbf{S}_0       & \mathbf{0}_m   \\ \mathbf{0}_m &\mathbf{0}_m\end{pmatrix}  \left( \mathbf{P}_3 ^{t} \right)^T,\end{align*}
where $\mathbf{P}_3 \triangleq  \rho_3^{-1/2}\begin{pmatrix} (1+\beta) \mathbf{I}_m -\alpha    \mathbf{H} &~-\beta  \mathbf{I}_m\\
\mathbf{I}_m &~\mathbf{0}_m
\end{pmatrix}$ and $\rho( \mathbf{P}_3)<1$.
 \end{theorem}
{\bf Proof.} As has been shown in the proof of Lemma \ref{Lem-HB-Q} (Equation \eqref{hb_recur1}),  we  can rewrite the recursion \eqref{Alg3} in a  matrix form as follows:
\begin{align*}
\begin{pmatrix}
 x_{k+1}-x^*
\\ x_k-x^* \end{pmatrix} =\underbrace{\begin{pmatrix} (1+\beta) \mathbf{I}_m -\alpha   \mathbf{H} &-\beta  \mathbf{I}_m\\
\mathbf{I}_m &\mathbf{0}_m
\end{pmatrix} }_{\triangleq \mathbf{H}_3}\begin{pmatrix}
 x_{k }-x^*
\\ x_{k-1}-x^* \end{pmatrix}  -\alpha\begin{pmatrix}     w_{k,N_k} \\ 0\end{pmatrix} .
\end{align*}
Define $\varepsilon_k \triangleq \rho_3^{-k/2} \begin{pmatrix}
 x_{k }-x^*
\\ x_{k-1}-x^* \end{pmatrix} .$  By multiplying   the above equation with $\rho_3^{-(k+1)/2}$,  one obtains
\begin{equation*}
\begin{split}
\varepsilon_{k+1}&  =  \mathbf{P}_3   \varepsilon_k   -  \alpha  \rho_3^{-(k+1)/2} \begin{pmatrix} \mathbf{I}_m  \\ \mathbf{0}_m\end{pmatrix}   w_{k,N_k}  ,
\end{split}
\end{equation*}
 where $ \mathbf{P}_3 = \rho_3^{-1/2} \mathbf{H}_3 $.   Thus, $\rho(\mathbf{P}_3)= \rho_3^{-1/2}  \rho(\mathbf{H}_3) <1$
  by recalling  $\rho(\mathbf{H}_3)= \sqrt{\beta}$  from  \eqref{eigen-alg3}. This together with  Assumption \ref{ass-CLT}(ii) proves    the result  by  using Lemma \ref{lem-CLT2}.
  \hfill $\blacksquare$

\begin{remark}\label{remark8}
 Suppose we set $\rho_3 = \left(1-{1 \over \sqrt{\kappa} +1}\right)^2$ in Algorithm \ref{Alg_3},  then $\rho_3>\beta $ and  Theorem \ref{thm-CLT1-Alg3} holds;  {i.e., $\tfrac{   (\sqrt{\eta}+\sqrt{L})^2}{4}  \left(1-\tfrac{1}{\sqrt{\kappa} +1}\right)^{-k} \begin{pmatrix}
x_{k }-x^* \\
x_{k-1}-x^* \end{pmatrix} \xlongrightarrow [k \rightarrow \infty]{d}   N(0,\bm{\Sigma}_3)  $.}
If we denote $V \sim  N(0,\mathbf{I}_{2m})$, then
\begin{align*}
& \begin{pmatrix}
x_{k }  \\
x_{k-1} \end{pmatrix} \overset {\tiny D}{\approx } \begin{pmatrix}
 x^* \\
 x^* \end{pmatrix} +\tfrac{ 4}{(\sqrt{\eta}+\sqrt{L})^2} \left( 1-\tfrac{1}{\sqrt{\kappa} +1}\right)^{k }\bm{\Sigma}_3^{1/2} V \mbox{~for large~}k.
\end{align*}
\end{remark}

Next,  based on Proposition \ref{prp-alg3} \red{and   Lemma \ref{lem-CLT2}}, we  also show the CLT  of  Algorithm  \ref{Alg_3}  on non-quadratic \red{objective}  functions.
 \begin{theorem}[CLT  for  Alg. \ref{Alg_3} with Geometrically Increasing $N_k$ on Non-Quadratic Functions]\label{thm-CLT1-Alg3-NQ}
Let  Assumption \ref{ass-fun} and \ref{ass-CLT}  hold. Consider Algorithm  \ref{Alg_3}, where     $   \beta = |1-\sqrt{\alpha \eta}|^2$ and  $\alpha\in (0, \tfrac{2( 1-\beta) }{ L+\eta})$.   Let  $q=\max\{q_1,q_2\}$ with
$q_1$ and $q_2$ defined in \eqref{def-q1q2}.   Set  $N_{k } \triangleq \lceil \rho_4^{-(k+1)} \rceil $  with   $\rho_4 \in \left(q,1\right)$.  Then
 \begin{align*}
  &\alpha^{-1} \rho_4^{-k/2} \begin{pmatrix}
x_{k }-x^* \\
x_{k-1}-x^* \end{pmatrix} \xlongrightarrow [k \rightarrow \infty]{d}   N(0,\bm{\hat{\Sigma}}_3) {~\rm with~}\bm{\hat{\Sigma}}_3 \triangleq  \sum_{t=0}^{\infty}  \mathbf{P}_4 ^{t} \begin{pmatrix} \mathbf{S}_0       & \mathbf{0}_m   \\ \mathbf{0}_m &\mathbf{0}_m\end{pmatrix}  \left( \mathbf{P}_4 ^{t} \right)^T,\end{align*}
where $\mathbf{P}_4 \triangleq  \rho_4^{-1/2}\begin{pmatrix} (1+\beta) \mathbf{I}_m -\alpha    \mathbf{H} &~-\beta  \mathbf{I}_m\\
\mathbf{I}_m &~\mathbf{0}_m
\end{pmatrix}$ and $\rho( \mathbf{P}_4)<1$.
 \end{theorem}
  {\bf Proof.}  From  \eqref{Alg3} and Assumption \ref{ass-CLT}(a) it follows  that
\begin{align*}
x_{k+1}-x^*& =x_k-x^* +\beta(x_k-x_{k-1}) -\alpha    \mathbf{H}   ( x_k -x^*) -\alpha \mathbf{D}(x_k)(x_k-x^*)  -\alpha  w_{k,N_k}.
\end{align*}Then we may rewrite the recursion in a  matrix form as follows:
\begin{align*}
\begin{pmatrix}
 x_{k+1}-x^*
\\ x_k-x^* \end{pmatrix} =\underbrace{\begin{pmatrix} (1+\beta) \mathbf{I}_m -\alpha   \mathbf{H} &-\beta  \mathbf{I}_m\\
\mathbf{I}_m &\mathbf{0}_m
\end{pmatrix} }_{\triangleq \mathbf{H}_4}\begin{pmatrix}
 x_{k }-x^*
\\ x_{k-1}-x^* \end{pmatrix}  -\alpha\begin{pmatrix}    \mathbf{D}(x_k)(x_k-x^*) +  w_{k,N_k} \\ 0\end{pmatrix} .
\end{align*}
Define $\varepsilon_k \triangleq \rho_4^{-k/2} \begin{pmatrix}
 x_{k }-x^*
\\ x_{k-1}-x^* \end{pmatrix} .$  By multiplying   the above equation with $\rho_4^{-(k+1)/2}$,  one obtains
\begin{equation*}
\begin{split}
\varepsilon_{k+1}&  =  \mathbf{P}_4   \varepsilon_k   -  \alpha  \rho_4^{-(k+1)/2} \begin{pmatrix} \mathbf{I}_m  \\ \mathbf{0}_m\end{pmatrix}   w_{k,N_k}  +\begin{pmatrix}   \varsigma_{k+1}    \\ 0\end{pmatrix} {\rm~with~} \varsigma_{k+1} \triangleq -\alpha  \rho_4^{-(k+1)/2}    \mathbf{D}(x_k)(x_k-x^*) ,
\end{split}
\end{equation*}
 where $ \mathbf{P}_4 = \rho_4^{-1/2} \mathbf{H}_4$.

    Next, we show that $\rho(\mathbf{H}_4)\leq  \sqrt{q} $ similar to how we showed tha     $\rho(\mathbf{H}_2)\leq 1-\gamma$ in Theorem \ref{thm-CLT1-Alg2}.
Recall the eigenvalue decomposition   $\mathbf{H}=\mathbf{U}\bm{\Lambda} \mathbf{U}^T$,  where $\mathbf{U}$ is orthogonal and $\bm{\Lambda} \triangleq {\rm diag}\{\lambda_1, \lambda_2,\cdots, \lambda_m\}$ with   satisfy $0<\eta \leq \lambda_m\leq \lambda_{m-1} \leq \dots\leq \lambda_2\leq \lambda_1\leq L.$ Then we can rewrite $$\mathbf{H}_4= \begin{pmatrix}  \mathbf{U} & {\bf 0}_m \\  {\bf 0}_m  & \mathbf{U}
\end{pmatrix}\begin{pmatrix} (1+\beta) \mathbf{I}_m -\alpha \bm{\Lambda}  &~ -\beta  \mathbf{I}_m\\
\mathbf{I}_m &~ \mathbf{0}_m
\end{pmatrix}\begin{pmatrix}  \mathbf{U} & {\bf 0}_m \\   {\bf 0}_m  & \mathbf{U}\end{pmatrix}^T.$$
 Similarly to the derivation of  \eqref{find-root}, the  eigenvalues of $\mathbf{H}_4$ can be described by the following  characteristic  equations.
 \[  p_i(v)=(v-  (1+\beta    -\alpha  \lambda_i)) v +\beta=v^2 -(1+\beta    -\alpha  \lambda_i) v+\beta=0,\quad i = 1, \cdots, m .\]

Since $0<\eta \leq \lambda_m\leq \lambda_{m-1} \leq \dots\leq \lambda_2\leq \lambda_1\leq L,$
we have $|1-\sqrt{\alpha  \lambda_i}| \leq \max\{|1-\sqrt{\alpha \eta}| ,  |1-\sqrt{\alpha L}| \}.$
Note from  $\alpha\in (0, \tfrac{2( 1-\beta) }{ L+\eta})$ that
$\alpha-\left(\tfrac{2}{\sqrt{L}+\sqrt{\eta}}\right)^{2}<\tfrac{2}{ L + \eta} -\tfrac{4}{ L + \eta+2\sqrt{\eta L}} \leq 0.$
Thus,   $|1-\sqrt{\alpha \eta}| \geq  |1-\sqrt{\alpha L}| $ and  $|1-\sqrt{\alpha  \lambda_i}| \leq  |1-\sqrt{\alpha \eta}| =\beta.$
    Then similar to \eqref{dictri-Q}, the discriminant of the  equation $p_i(v)  =0 $  is
\begin{align*}
\Delta_i & =(1+\beta    -\alpha  \lambda_i)^2 -4  \beta = \beta^2-2\beta(1    +\alpha  \lambda_i)+ (1    -\alpha  \lambda_i) ^2
= (\beta -  (1    +\alpha  \lambda_i))^2-4\alpha  \lambda_i \leq 0.
\end{align*}
    Then the  spectral radius  of $\mathbf{H}_4$ is
\begin{align}\label{case1} \rho(\mathbf{H}_4) = \sqrt{\beta} \geq  |1-\sqrt{\alpha \eta}| =1-\sqrt{\alpha \eta} {\rm~since~}\alpha \eta<1.
\end{align}

Since $\kappa>1 $ and $\beta \in (0,1)$, we have  $(1-\beta)(\kappa+1)+2\kappa \beta=\kappa+1+\beta(\kappa-1)>\kappa+1>\tfrac{4\kappa}{\kappa+1}$.
Therefore,
\begin{align*}
& (1-\beta)^2(\kappa+1)+2\kappa \beta (1-\beta) >\tfrac{4\kappa (1-\beta)}{\kappa+1}=\tfrac{4  (1-\beta)L}{L+\eta}>2\alpha L .
\\& \Rightarrow 1-\beta>\tfrac{2\alpha L}{(1-\beta)(\kappa+1)+2\kappa \beta }=\tfrac{2\alpha L \eta}{(1-\beta)(L+\eta)+2L \beta }
\\& \Rightarrow \beta<1-\tfrac{2\alpha L \eta}{(1-\beta)(L+\eta)+2L \beta } \overset{\eqref{def-q1q2}}{\leq} q_2.
\end{align*}
This together with   \eqref{case1} implies that
$  \rho(\mathbf{H}_4)<\sqrt{q_2}\leq \sqrt{q}$ since $q=\max\{q_1,q_2\}$.

    From Proposition \ref{prp-alg3} it follows that  $ \mathbb{E}[\rho_4^{-k}\|x_{k}-x^*\|^2]  \leq C $ for some constant $C>0$  and hence
$ \| x_k - x^*\| =o_P(1)$.   Then by  invoking the bound that  $ \| \mathbf{D}(x) \| \leq R_D \| x- x^*\|  $, we obtain that  $\| \mathbf{D}(x_k)\|=  o_P(1) .$  Similarly to the procedures for proving \eqref{bd-absek},  we can   show  that
 $\| \rho_4^{-k/2}(x_k-x^*)\|=O_P(1)$.  Therefore,
  \begin{align*}
  & \| \varsigma_{k+1}\| \leq \alpha  \rho_3^{-1/2} \|\mathbf{D}(x_k)\|  \| \rho_3^{-k/2}(x_k-x^*)\|= o_P(1) O_P(1) \overset{ \eqref{bigO} }{=}  o_P(1)
  \Rightarrow  \varsigma_{k+1}  \xlongrightarrow [k \rightarrow \infty]{P} 0  .
  \end{align*}    Since   $\rho(\mathbf{P}_4)= \rho_4^{-1/2}  \rho(\mathbf{H}_4)= \rho_4^{-1/2}\sqrt{q}<1$  by $  \rho(\mathbf{H}_4)<  \sqrt{q}$ and $\rho_4>q$,   the result follows by invoking Lemma \ref{lem-CLT2}.
  \hfill $\blacksquare$

    Since    $   \beta = |1-\sqrt{\alpha \eta}|^2$ and  $\alpha\in (0, \tfrac{2( 1-\beta) }{ L+\eta})$  imply that $\beta$ is defined in an implicit sense,  we add the following remark to clarify their existence and show how does the spectral  radius depends on the condition number $\kappa.$
\begin{remark}\label{rem9} By substituting   $   \beta = |1-\sqrt{\alpha \eta}|^2$ into the upper bound of $\alpha$, it  requires that
\[\alpha < \tfrac{2( 1- |1-\sqrt{\alpha \eta}|^2) }{ L+\eta} \Leftrightarrow \alpha( L+\eta) < 4\sqrt{\alpha \eta}-2\alpha \eta
 \Leftrightarrow \alpha( L+3\eta) < 4\sqrt{\alpha \eta}  \Leftrightarrow \alpha <{16\eta \over( L+3\eta)^2 }.\]
   Then by \eqref{case1}, we achieve  $ \rho(\mathbf{H}_4)\geq 1-\sqrt{\alpha \eta}>1-\sqrt{{16\eta^2 \over( L+3\eta)^2 }}
    =1-\sqrt{{16  \over( \kappa+3 )^2 }} = 1- \tfrac{4}{\kappa+3} \ge  1-\tfrac{4}{\kappa} $. 
    In comparison with the geometric rate $ 1- 1/\sqrt{\kappa+1}$ of Algorithm \ref{Alg_3} for quadratic functions,  the geometric parameter  of Algorithm \ref{Alg_3}  non-quadratic functions  will not exceed $1-\tfrac{4}{\kappa}$.  Therefore, it does not lead to an acceleration compared with variable sample-size stochastic gradient algorithm (Algorithm \ref{Alg_1}),
which is consistent with the statement in \cite[Lemma 2.5]{goujaud2023provable}.
\end{remark}

   \section{Central Limit Theorems on Polynomial  Batch-size.}\label{sec:CLT-poly}

 There are many settings where a geometric increase in $N_k$
is impractical. For instance,   the
generation of a sampled gradient is computationally
expensive; an example of this commonly arises in simulation
optimization problems in the context of large manufacturing
or queueing simulations. To this end, we consider the use of
polynomial increases in sample-size, an avenue that allows
for more gentle growth,  and proceed to   investigate  the
rate  of  convergence and the associated central limit
statements for Algorithms~\ref{Alg_1}--\ref{Alg_3}. Polynomial rates have been studied in~\cite{lei18asynchronous} as well as~\cite{pasupathy18sampling} but remain unaware of rate and complexity statements as well as the ensuing CLTs in accelerated and heavy-ball regimes.
 \subsection{ Rate of Convergence.}
 We first recall  some  preliminary results     that find utility in  the
 rate analysis of the proposed algorithms with  the polynomially increasing batch-sizes.
\begin{lemma} \label{lem-recur}  (i)   \cite[Eqn. (17)]{lei18game}  For any $q\in (0,1)$ and $v>0$, there  holds   \[\sum_{t=1}^{k  } q^{k  -t} t^{-v}    \leq  q^{k  }
\tfrac{e^{2v}q^{-1}-1 }{ 1-q}+\tfrac{2  k^{-v}}{q \ln(1/q) }.\]
\noindent  (ii) \cite[Lemma 4]{lei18game}  For any $q\in (0,1)$ and $v>0$,  $q^x \leq c_{q,v}  {x^{-v}}$ for all $x >0$ where $c_{q,v} \triangleq   e^{- v  }{\left(\tfrac{v}{ \ln(1/q)}\right)^{ v  }}.$
\end{lemma}
Based on Lemmas \ref{Lemma_1}-\ref{Lem-HB}  in Section \ref{Sec:Alg} and Lemma \ref{lem-recur}, we can establish polynomial  rates of convergence of the iterates generated by the three proposed methods. {Omitted proofs  are included   in the supplementary material for purposes of completeness.}
 \begin{proposition}[Rate statement for Algorithm  \ref{Alg_1} under polynomially increasing $N_k$]  \label{prp3}    Suppose Assumption \ref{ass-fun} holds  and that   $N_{k }   \triangleq \lceil (k+1)^v\rceil $ for some $v>0$. Consider  Algorithm  \ref{Alg_1} with $  \alpha \in (0, {2\over \eta+L}]$.   Define  $q\triangleq 1-{2 \alpha \eta L \over \eta+L} $ and $c_{q,v} \triangleq   e^{- v  }{\left(\tfrac{v}{ \ln(1/q)}\right)^{ v  }}.$   Then
\begin{align}\label{rate-pol}
 \mathbb{E}[\| x_{k}-x^*\|^2  ]  \leq  \left(\underbrace{  c_{q,v} \mathbb{E}[ \| x_0 - x^*\|^2] + \tfrac{\alpha^2\nu^2 c_{q,v}(e^{2v}q^{-1}-1)}{1-q}+\tfrac{2\alpha^2\nu^2}{q \ln(1/q) }}_{C_v}\right) k^{-v},\quad \forall k\geq 1 . \end{align}
 Then  the   number of iterations  and sampled gradients  required    to  obtain an $\epsilon $-optimal solution in the mean-squared sense (i.e. $\mathbb{E}[\| x  -x^*\|^2  ] \leq \epsilon$) are   $  \mathcal{O}  \big(v \left ({1/\epsilon} \right)^{1/v} \big)$ and  $  \mathcal{O}  \big( e^v v^v \left ({1/\epsilon} \right)^{1+1/v} \big)$, respectively. \end{proposition}

 \begin{proposition}[Rate statement	 for  Algorithm  \ref{Alg_2} under polynomially increasing $N_k$]  \label{prp4}Let Algorithm  \ref{Alg_2}
be applied to   \eqref{Problem1}, where $N_{k }   \triangleq \lceil (k+1)^v\rceil $ with some   $v>0$. Suppose Assumption \ref{ass-fun} holds and   $   \alpha \in ( 0, {1\over L}]$.
Define  $\gamma\triangleq \sqrt{\alpha \eta}$ and  $\beta \triangleq {1 -\gamma \over 1+ \gamma}.$  Then there exists a constant $C (v) >0$ such that \[
   \mathbb{E}[f(y_{k }) ]-f^* \leq   C (v)  k^{-v},\quad \forall k\geq 1 .\]
\end{proposition}

\begin{remark} \label{rem-poly}   Since $f(x)-f (x^*) \geq {\eta \over 2} \|x-x^*\|^2,$   we have
  $\mathbb{E}[\|y_k-x^*\|^2] \leq \tfrac{2 C_v }{\eta}k^{-v}.$ Then from \eqref{Alg22} it follows that    $\mathbb{E}[\|x_{k}-x^*\|^2] \leq 2(1+\beta)^2\mathbb{E}[\|y_k-x^*\|^2] +2\beta^2 \mathbb{E}[\|y_{k-1}-x^*\|^2] \leq c k^{-v}  $ for some   $c>0$.
Thus,  the mean-squared error of Algorithm  \ref{Alg_2} also  displays  the polynomial rate of convergence similar to that shown in Proposition \ref{prp3} for Algorithm  \ref{Alg_1}. Because   the mean-squared convergence implies converges in probability,    the sequences   $\{x_k\}$ and $\{y_k\}$ generated by Algorithm   \ref{Alg_2} satisfy    $ x_k  \xlongrightarrow [k \rightarrow \infty]{P}  x^* $ and $ y_k  \xlongrightarrow [k \rightarrow \infty]{P}  x^* $ when the conditions of Proposition \ref{prp4} hold.
\end{remark}

\begin{proposition}[Rate of Alg.  \ref{Alg_3}  on quadratic Functions with polynomially increasing $N_k$]  \label{prp2-alg3-Q}
Suppose that Assumption \ref{ass-fun} holds and that   {$f$ is a quadratic function with   $  \nabla^2 f(x) \equiv \mathbf{H} $ for any $x\in \mathbb{R}^m.$}     Consider Algorithm  \ref{Alg_3}, where $   \alpha \triangleq \tfrac{ 4}{(\sqrt{\eta}+\sqrt{L})^2}$,  $\beta \triangleq    ( \tfrac{\sqrt{\kappa} -1}{\sqrt{\kappa} +1} )^2 $,  and  $N_{k }   \triangleq \lceil (k+1)^v\rceil $ with  $v>0$. Then    there exists a constant $C(v)>0$ such that
$ \mathbb{E} \left [ \left \| \begin{pmatrix}
  x_{k+1}-x^*
\\ x_k-x^* \end{pmatrix} \right\|^2  \right]\leq  C(v) (k+1)^{-v},~ \forall k\geq 0. $
\end{proposition}

\begin{proposition}[Rate of Alg.  \ref{Alg_3}  on Non-Quadratic Functions with polynomially increasing $N_k$]  \label{prp2-alg3}
Let  Assumption \ref{ass-fun} hold. Consider Algorithm  \ref{Alg_3}, where     $   \beta \in (0,1)$, $\alpha\in (0, \tfrac{2( 1-\beta) }{ L+\eta})$, and  $N_{k }   \triangleq \lceil (k+1)^v\rceil $ with  $v>0$. Then    there exists a constant $C(v)>0$ such that
\[  \mathbb{E}[\|x_k-x^*\|^2]\leq C(v) (k+1)^{-v},~ \forall k\geq 0 .\]
\end{proposition}

\subsection{Central Limit Theorems under Polynomially Increasing $N_k$.}

In this part, we first establish  the  asymptotic normality of  a     time-varying linear  recursion and   provide   the proof  in Appendix \ref{app-lem-CLT3}.  This result will be applied in  {proving} Theorems \ref{thm-CLT3}-\ref{thm-CLT2-Alg3}.

\begin{lemma}\label{lem-CLT3}  Suppose  that  the square matrix $  \mathbf{A}$ satisfies $\rho(  \mathbf{A}) <1$,    $N_k=\lceil ( k+1)^v \rceil$ with $v>0$, and $\{ w_{k,N_k} \} $ satisfies Assumption \ref{ass-CLT}(ii).  Let the sequence  $\{e_k\}$ be generated by
\begin{align}\label{poly-recursion4}
e_{k+1}& =  \mathbf{A}_k    e_k -\alpha(k+1)^{v/2}  \mathbf{G}w_{k ,N_k }+ \zeta_{k+1}, ~~\mathbb{E} [\|e_0\|^2]<\infty,
\end{align}
where  $\mathbf{A}_0=\mathbf{A}$, $ \mathbf{A}_k \triangleq   {\left(k+1 \over k\right)}^{v/2} \mathbf{A}$  for any $k\geq 1 $, and $\zeta_{k}    \xlongrightarrow [k \rightarrow \infty]{P} 0 $.
Then
 \begin{align*}
  &\alpha^{-1}   e_k\xlongrightarrow [k \rightarrow \infty]{d}   N(0,\bm{\Sigma})
~ {\rm with} ~  \bm{\Sigma}  \triangleq \lim_{k\to \infty}\sum_{t=1}^{k} {\left(k  \over  t\right)}^{v}  \mathbf{A} ^{k-t}  \mathbf{G} \mathbf{S}_0    \mathbf{G}^T \left( \mathbf{A}^T \right)^{k-t}  .\end{align*}
\end{lemma}

Based on  Proposition \ref{prp3} and  Lemma \ref{lem-CLT3},  by using
Assumption \ref{ass-CLT}, we are now  ready to derive the associated  central
limit theorem  for Algorithm \ref{Alg_1} with polynomially increasing
batch-sizes.

 \begin{theorem}[CLT  of Algorithm \ref{Alg_1} with  Polynomially increasing    Batch-sizes]\label{thm-CLT3} Suppose that Assumptions \ref{ass-fun}  and \ref{ass-CLT} hold.  Consider Algorithm \ref{Alg_1}, where   $ \alpha \in (0, {2\over \eta+L}]$ and  $N_{k }   \triangleq \lceil (k+1)^v\rceil $ with some $v>0$. Define $\mathbf{A} \triangleq  \mathbf{I}_m-\alpha  \mathbf{H}$. Then $\rho( \mathbf{A})<1$ and
 \[\alpha^{-1} k^{v/2}( x_{k}-x^*)  \xlongrightarrow [k \rightarrow \infty]{d}   N(0,  \bm{\Sigma}_4) ~
{\rm with} ~ \bm{\Sigma}_4\triangleq \lim_{k\to \infty}\sum_{t=1}^k  {\left(k  \over  t\right)}^{v} \mathbf{A}^{k-t}  \mathbf{S}_0 \mathbf{A}^{k-t}  . \]
 \end{theorem}
 {\bf Proof. }    We begin by noting that  $x_k\xlongrightarrow [k \rightarrow \infty]{P} x^*$  and
 \eqref{recursion3}   holds.
Define $e_0\triangleq x_0-x^* $, $e_{k} \triangleq k^{v/2}(x_{k}-x^*)$ for any $k\geq 1$, and $\zeta_{k+1} \triangleq -\alpha(k+1)^{v/2} \ \mathbf{D}(x_k )(x_k-x^*) $ for any $k\geq 0$. By multiplying  both sides of  \eqref{recursion3} with $(k+1)^{v/2}$,   and using   $ \mathbf{A}_k =  {\left(k+1 \over k\right)}^{v/2}(\mathbf{I}_m-\alpha  \mathbf{H}) $, we achieve
 \begin{align} \label{def-ek}
e_{k+1}& =  {\left( \tfrac{k+1 }{ k}\right)}^{v/2} (\mathbf{I}_m-\alpha  \mathbf{H}) k^{v/2}(x_{k}-x^*) -\alpha(k+1)^{v/2}( \mathbf{D}(x_k )(x_k-x^*) +w_{k,N_k}) \notag
\\& =   \mathbf{A}_k  e_k - \alpha(k+1)^{v/2}w_{k,N_k}+\zeta_{k+1} , \quad \forall k\geq 1  .
\end{align}
 By setting  $k=0$ in   \eqref{recursion3},     and using $ \mathbf{A}_0=\mathbf{A}=\mathbf{I}_m-\alpha  \mathbf{H}$,  we see that
 \begin{align*} e_1=x_1-x^*=\mathbf{A}_0(x_0-x^*)-\alpha  w_{1,N_1}-\alpha \mathbf{D}(x_0 )(x_0-x^*) .
 \end{align*}
Then by $\zeta_1=-\alpha  \mathbf{D}(x_0 )(x_0-x^*) , $  we see that the
\eqref{def-ek} also  holds for $k=0.$ Hence the recursion \eqref{def-ek}  holds for any $k\geq 0.$
From  \eqref{bd-P} it follows  that  the symmetric matrix $\mathbf{A}$ satisfies
 $$\rho(\mathbf{A})=\|\mathbf{A}\|_2 =\| \mathbf{I}_m-\alpha  \mathbf{H}  \|_2= \rho^{1/2} \| \mathbf{P}_1\|_2 \leq q^{1/2} <1.$$


 We conclude from  Proposition \ref{prp3}    that $\mathbb{E}[k^v\| x_{k}-x^*\|^2   ] \leq c $   for some constant $c>0,$ and $ \| x_k - x^*\| =o_P(1)$.
Then by  invoking the bound that  $ \| \mathbf{D}(x) \| \leq R_D \| x- x^*\|  $,  we  obtain $\| \mathbf{D}(x_k)\|= o_P(1).$ Similarly to the procedures for proving  \eqref{bd-absek},  we can   show  that $ \|  k^{v/2}(x_{k}-x^*)\|=O_P(1)$.   Therefore,
\begin{align}\label{limit-zeta2}
\| \zeta_{k+1} \| \leq \alpha {\left (1+\tfrac{1}{ k}\right)}^{v/2} \|\mathbf{D}(x_k)\| \|  k^{v/2}(x_{k}-x^*)\| =  o_P(1) O_P(1) \overset{ \eqref{bigO} }{=}  o_P(1)
    \Rightarrow \zeta_{k+1}  \xlongrightarrow [k \rightarrow \infty]{P} 0 .
\end{align}

 Therefore,  by using Lemma \ref{lem-CLT3}  with
$\mathbf{G}=\mathbf{I}_m, $ we conclude that
$\alpha^{-1}e_k \xrightarrow [k \rightarrow \infty]{d}   N(0,\bm{\Sigma}_4) $.
 By {the definition  $  e_k=k^{v/2}( x_{k}-x^*)$,}  we obtain the result.  \hfill $\blacksquare$

Similarly to Theorem \ref{thm-CLT3},  based on Proposition \ref{prp4}
and   Lemma \ref{lem-CLT3},   we    now  establish  the central
limit theorem for Algorithm \ref{Alg_2}  under the assumption of
 polynomially increasing  batch-sizes.

 \begin{theorem}[CLT of Algorithm \ref{Alg_2} under Polynomially increasing  Batch-sizes]\label{thm-CLT2-Alg2} Suppose Assumptions \ref{ass-fun}  and \ref{ass-CLT} hold.  Consider Algorithm \ref{Alg_2}, where  $   \alpha \in (0, {1\over  L}]$ and $N_{k }=\lceil (k+1)^v\rceil $  with   some $v>0$.  Define $\gamma \triangleq \sqrt{ \alpha \eta } $, $\beta \triangleq {1 -\gamma \over 1+ \gamma},$  and $\mathbf{H}_2 \triangleq  \begin{pmatrix}
(1+\beta) (\mathbf{I}_m-\alpha  \mathbf{H})  & ~-\beta (\mathbf{I}_m-\alpha  \mathbf{H})
\\  \mathbf{I}_m   &~\mathbf{0}_m \end{pmatrix}.$ Then $\rho( \mathbf{H}_2)<1$ and
 \begin{align*}
  &\alpha^{-1} k^{v/2} \begin{pmatrix}
y_{k }-x^* \\
y_{k-1}-x^* \end{pmatrix} \xlongrightarrow [k \rightarrow \infty]{d}   N(0, \bm{\Sigma}_5)  {\rm~with~}\bm{\Sigma}_5 \triangleq \lim_{k\to \infty} \sum_{t=1}^{k}   {\left(k  \over  t\right)}^{v}   \mathbf{H}_2 ^{k-t} \begin{pmatrix} \mathbf{S}_0       & \mathbf{0}_m   \\ \mathbf{0}_m &\mathbf{0}_m\end{pmatrix}  \left( \mathbf{H}_2 ^T\right)^{k-t}  .\end{align*}
 \end{theorem}
{\bf Proof.}     Define $z_{k+1} \triangleq     \begin{pmatrix}
y_{k+1}-x^* \\
y_{k}-x^* \end{pmatrix} $, $\varepsilon_0 \triangleq  z_0$, and $\varepsilon_k \triangleq k^{v/2} z_k$ for any $  k\geq 1.$
Therefore, by defining $ \varsigma_{k+1} \triangleq     -  \alpha   (k+1)^{v/2}   \mathbf{D}(x_k)(x_k-x^*) $,
 and  multiplying  both sides of  \eqref{re-z} by  $(k+1)^{v/2}$,     we obtain that
\begin{equation*}
\begin{split}
\varepsilon_{k+1}&  = {\left(k+1 \over k\right)}^{v/2}\mathbf{H}_2 \varepsilon_k   -  \alpha  (k+1)^{v/2} \begin{pmatrix}  \mathbf{I}_m  \\ \mathbf{0}_m\end{pmatrix}   w_{k ,N_k }  +  \begin{pmatrix}  \varsigma_{k+1}  \\ 0 \end{pmatrix}  .
\end{split}
\end{equation*}
By defining    $\mathbf{G}\triangleq  \begin{pmatrix} \mathbf{I}_m    \\ \mathbf{0}_m\end{pmatrix}$, $\mathbf{A}_0  \triangleq \mathbf{H}_2 ,   $  and $ \mathbf{A}_k \triangleq   {\left(k+1 \over k\right)}^{v/2}  \mathbf{H}_2 $, there holds
\begin{equation}  \label{thm4-vare}
\begin{split}
\varepsilon_{k+1}&  = \mathbf{A}_k \varepsilon_k   -  \alpha  (k+1)^{v/2}  \mathbf{G}  w_{k ,N_k } +\begin{pmatrix}  \varsigma_{k+1}  \\ 0 \end{pmatrix},\quad \forall k\geq 1.
\end{split}
\end{equation}
From  \eqref{re-z} it is  seen  that    \eqref{thm4-vare}    holds for $k=0 $ as well.
Thus,  the recursion \eqref{thm4-vare} holds for any $k\geq 0.$

 Note from  Remark \ref{rem-poly} that  $\mathbb{E}[k^v\| x_{k}-x^*\|^2   ] \leq c $   for some constant $c>0,$ and $ \|x_k - x^*\|=o_P(1)$. Then by  invoking the bound that  $ \| \mathbf{D}(x) \| \leq R_D \| x- x^*\|  $,  we derive $\| \mathbf{D}(x_k)\|= o_P(1).$
 Similarly to the procedures for proving \eqref{bd-absek},  we can   show  that
$ \| k^{v/2}(x_{k}-x^*)\|=O_P(1)$. Therefore,
\[\| \varsigma_{k+1}   \| \leq \alpha {\left (1+\tfrac{1}{ k}\right)}^{v/2} \|\mathbf{D}(x_k)\| \|  k^{v/2}(x_{k}-x^*)\| =  o_P(1) O_P(1) \overset{ \eqref{bigO} }{=}  o_P(1).\]
    Hence $\varsigma_{k}   \xlongrightarrow [k \rightarrow \infty]{P} 0.$
Then by     using Lemma \ref{lem-CLT3} and  \eqref{sp-h}, we obtain that
\[\alpha^{-1}  \varepsilon_k \xlongrightarrow [k \rightarrow \infty]{d}   N(0, \bm{\Sigma}_5),
 {\rm~where~}   \bm{\Sigma}_5    \triangleq \lim_{k\to \infty}\sum_{t=0}^{k}
  {\left( \tfrac{k }{  t}\right)}^{v}  \mathbf{H}_2 ^{k-t}  \mathbf{G} \mathbf{S}_0    \mathbf{G}^T \left( \mathbf{H}_2^T \right)^{k-t}  .\]
Then by the fact that $\mathbf{G} \mathbf{S}_0    \mathbf{G}^T= \begin{pmatrix} \mathbf{S}_0      & \mathbf{0}_m   \\ \mathbf{0}_m &\mathbf{0}_m\end{pmatrix}$ and $\varepsilon_k=k^{v/2}  \begin{pmatrix}
y_{k }-x^* \\ y_{k-1}-x^* \end{pmatrix}$,  we prove the result.
\hfill  $\blacksquare$

 \begin{remark}\label{rem8}
  Suppose we set $\alpha={2\over  L+\eta}$ in    Algorithm \ref{Alg_1}.  Then $q=1-{2 \alpha \eta L \over \eta+L}=\left({\kappa-1\over \kappa+1} \right)^2$.
By \eqref{bd-P}, it is seen that  the matrix $\mathbf{A}$, defined in Theorem \ref{thm-CLT3}, satisfies  $\|  \mathbf{A}  \|_2  \leq      q ^{1/2}={\kappa-1\over \kappa+1}=1-{2\over \kappa+1}.$
  Suppose $\alpha={1\over  L}$ in Algorithm \ref{Alg_2}.  Then   $\gamma=1/\sqrt{\kappa} $ and by \eqref{sp-h},  we know that $\mathbf{H}_2$, defined in Theorem \ref{thm-CLT1-Alg2}, satisfies $\| \mathbf{H}_2 \|_2  =1-\gamma =1-{1\over \sqrt{\kappa}}.$ Thus, we conclude from Theorems \ref{thm-CLT3} and \ref{thm-CLT2-Alg2}
that both  the unaccelerated gradient method  (Algorithm \ref{Alg_1})
and its    accelerated  counterpart (Algorithm \ref{Alg_2}) have convergence rates  with the same order   $\|x_k-x^*\|=\mathcal{O}(k^{-v/2})$;
however the   accelerated scheme  has a smaller constant than its unaccelerated counterpart  since  $\| \mathbf{H}_2 \|_2<\| \mathbf{A} \|_2$ due to the fact that ${2\over \kappa+1} \leq {1\over \sqrt{\kappa}}$.   \end{remark}

Similar to Theorem \ref{thm-CLT3},  based on Proposition  \ref{prp2-alg3-Q} for quadratic \red{objective} functions  (resp. Proposition \ref{prp2-alg3} for non-quadratic \red{objective}  functions)
and   Lemma \ref{lem-CLT3},   we    now  state  the central limit theorem for  Algorithm \ref{Alg_3} with polynomially  increasing batch-size (proof omitted).

 \begin{theorem}[CLT of Alg.  \ref{Alg_3} on Quadratic Functions with  polynomially increasing $N_k$]\label{thm-CLT2-Alg3}
 {Let Assumptions \ref{ass-fun}  and \ref{ass-CLT}(ii) hold. Suppose that  $f$ is a quadratic function with   $  \nabla^2 f(x) \equiv \mathbf{H} $ for any $x\in \mathbb{R}^m.$
Consider  Algorithm \ref{Alg_3}, where  $   \alpha \triangleq \tfrac{ 4}{ (\sqrt{\eta}+\sqrt{L})^2}$,  $\beta \triangleq   \left(\tfrac{\sqrt{\kappa} -1 }{ \sqrt{\kappa} +1}\right)^2,$  and   $N_k=\lceil (k+1)^v\rceil,~ v>0 $. Set $\mathbf{H}_3 \triangleq   \begin{pmatrix} (1+\beta) \mathbf{I}_m -\alpha    \mathbf{H} &~-\beta  \mathbf{I}_m\\ \mathbf{I}_m &~\mathbf{0}_m \end{pmatrix}.$ Then
 $\rho( \mathbf{H}_3)<1$ and  \begin{align*}
    &\alpha^{-1} k^{v/2} \begin{pmatrix}
x_{k }-x^* \\
x_{k-1}-x^* \end{pmatrix} \xlongrightarrow [k \rightarrow \infty]{d}   N(0, \bm{\Sigma}_6) ~{ \rm with~} \bm{\Sigma}_6 \triangleq \lim_{k\to \infty}\sum_{t=1}^{k}   {\left(k  \over  t\right)}^{v}   \mathbf{H}_3 ^{k-t} \begin{pmatrix} \mathbf{S}_0       & \mathbf{0}_m   \\ \mathbf{0}_m &\mathbf{0}_m\end{pmatrix}  \left( \mathbf{H}_3 ^T\right)^{k-t}.\end{align*}}
 \end{theorem}

 \begin{theorem}[CLT of Alg.  \ref{Alg_3} on Non-Quadratic Functions with  polynomially increasing $N_k$] \label{thm8}
 Let  Assumptions \ref{ass-fun}   and \ref{ass-CLT} hold. Consider Algorithm  \ref{Alg_3}, where     $   \beta = |1-\sqrt{\alpha \eta}|^2$, $\alpha\in (0, \tfrac{2( 1-\beta) }{ L+\eta})$, and  $N_k=\lceil (k+1)^v\rceil,~ v>0 $.  Set $\mathbf{H}_4 \triangleq  \begin{pmatrix} (1+\beta) \mathbf{I}_m -\alpha    \mathbf{H} &~-\beta  \mathbf{I}_m\\
\mathbf{I}_m &~\mathbf{0}_m
\end{pmatrix}$.  Then $\rho( \mathbf{H}_4)<1$ and
 \begin{align*}
    &\alpha^{-1} k^{v/2} \begin{pmatrix}
x_{k }-x^* \\
x_{k-1}-x^* \end{pmatrix} \xlongrightarrow [k \rightarrow \infty]{d}   N(0, \bm{\hat{\Sigma}}_6) ~{ \rm with~} \bm{\hat{\Sigma}}_6 \triangleq \lim_{k\to \infty}\sum_{t=1}^{k}   {\left(k  \over  t\right)}^{v}   \mathbf{H}_4 ^{k-t} \begin{pmatrix} \mathbf{S}_0       & \mathbf{0}_m   \\ \mathbf{0}_m &\mathbf{0}_m\end{pmatrix}  \left( \mathbf{H}_4 ^T\right)^{k-t}.\end{align*}
 \end{theorem}

\section{Confidence Regions of the Optimal Solution.} \label{sec:CI}

 A crucial motivation for developing CLTs lies in developing
confidence statements. In this section,  we  proceed to construct  the
confidence regions  for  the optimal solution $x^*$. Note that the limiting covariance matrix is dependent on
the Hessian at the solution, which is unavailable. Furthermore, we do
     not have an consistent  estimator for  the covariance matrix.  Yet, in
the absence of such an estimate, we proceed to develop rigorous
confidence statements, adopting an   approach developed in~\cite{hsieh2002recent}.

 Since Algorithms \ref{Alg_1}, \ref{Alg_2}, and  \ref{Alg_3} lead to  similar central limit results,
we show how to construct  confidence regions  merely for the sequence
$\{x_k\}$ generated by Algorithm \ref{Alg_1}.  The  simulation
framework  in \cite{hsieh2002recent}  lies in generating $n$
independent  replications of Algorithm \ref{Alg_1}, leading to $n$ copies of the  random iterate  $ x_k $,  denoted by
 $x_{1k},\cdots, x_{nk}.$
Then the sample mean  and  the covariance estimator    are respectively defined as
\begin{align}\label{est-covariance}
\bar{\mathbf{x}}_k={1\over n}\sum_{i=1}^n x_{ik} \mbox{ ~and } \quad \mathbf{S}_k={1\over n-1}\sum_{i=1}^n (x_{ik}-\bar{\mathbf{x}}_k) (x_{ik}-\bar{\mathbf{x}}_k)^T.
\end{align}
  Based on  \cite[Theorem 4 and Corollary 1]{bodnar2016singular} and  \cite[Theorem 2]{hsieh2002recent},
 we   may achieve    the following result. The proof  is  given in Appendix \ref{app-cf-clt} for completeness.

\begin{proposition}\label{cf-clt}
Consider Algorithm \ref{Alg_1}. Suppose that the conditions of Theorem  \ref{thm-CLT1}  hold and  $n\geq m+1$. Then \\
(i) $\sqrt{n} \alpha^{-1} \rho_1^{-k/2}( \bar{\mathbf{x}}_k-x^*)  \xlongrightarrow [k \rightarrow \infty]{d}   N(0,\bm{\Sigma}_1).$ \\
(ii)  $n ( \bar{\mathbf{x}}_k-x^*)^T\mathbf{S}_k^{-1}   ( \bar{\mathbf{x}}_k-x^*)  \xlongrightarrow [k \rightarrow \infty]{d} {m(n-1) \over n-m}F(m,n-m), $ where $F(m,n-m)$ denotes the  $F$-distribution  with $(m,n-m)$ degrees of freedom.
\end{proposition}

Proposition  \ref{cf-clt}  can be used to construct the   confidence region  of the optimal solution.   Define
\begin{align}\label{def-conf-region}
X_{mk}(z) \triangleq \left \{ x \, \mid \,  n ( \bar{\mathbf{x}}_k-x )^T\mathbf{S}_k^{-1}   ( \bar{\mathbf{x}}_k-x ) \leq  \tfrac{m(n-1) }{ n-m}z \right \},
\end{align}
where $z$ is selected such that $\mathbb{P}(F(m,n-m) \leq z) \geq 1-\delta$ with  some $\delta\in(0,1)$.   Then we have the  following corollary.

\begin{corollary} (\cite[Proposition 3]{hsieh2002recent}) \label{cor2} Consider Algorithm \ref{Alg_1}. Suppose that  $n\geq m+1$ and the conditions of Theorem  \ref{thm-CLT1}  hold.  Then the  confidence region    $X_{mk}(z) $   defined in  \eqref{def-conf-region}  is asymptotically correct, i.e.,
$ \lim\limits_{k \rightarrow \infty}\mathbb{P}\big(x^*\in X_{mk}(z) \big) = 1-\delta.$
\end{corollary}

The above result asserts that  the  estimated confidence region $X_{mk}(z) $     asymptotically covers the  optimal solution  $x^*$  with probability  $100(1-\delta)\%$. The approach is easily implementable  because it   merely  requires $n$ independent   replications of  Algorithm \ref{Alg_1},
 while without  requiring  a consistent estimator  of the covariance  matrix of the stationary normal distribution.
The   confidence regions of  Algorithm \ref{Alg_2} and Algorithm \ref{Alg_3}  can be constructed in a similar way.
 When it is expensive to run multiple independent
     trials, an alternative  avenue for reducing the complexity requirements can be found in~\cite{zhu2019constructing},  where the authors employ  a batch-means method for  constructing the
     confidence region.  This framework is reliant on a cancellation approach to ``cancel'' out the covariance matrix that is hard to estimate.
  This extension will be considered in  future work.

 In the following, we show that
the sequence  $\sqrt{\tfrac{n(n-m) }{ n-1}}S_k^{-1/2}(\bar{x}_k-x^*)$  converges  to  a suitably defined  multivariate $t$ distribution, the proof  of which  is given in Appendix \ref{app-cf-clt2}.

\begin{proposition}\label{cf-clt2}
Consider Algorithm \ref{Alg_1}. Suppose that the conditions of Theorem  \ref{thm-CLT1}  hold and  $n> m+2$. Then
\begin{align}\label{t-dis}
 \sqrt{\tfrac{n(n-m)}{ n-1}}    \mathbf{S}_k^{-1/2}( \bar{\mathbf{x}}_k-x^*)   \xlongrightarrow [k \rightarrow \infty]{d}   T_{n-m}(0,\mathbf{I}_m,m),
 \end{align}
  where a $m$-variate random vector $X\sim T_v(\mu,\bm{\Lambda},m)$ (i.e., $X$ is a $t$ distribution with mean $\mu$, covariance matrix $v(v-2)^{-1}\bm{\Lambda}, v>2$) if it has a probability density function $ h $ given by
\begin{align}\label{def-prob}
 h(u) =\tfrac{\Gamma((v+m)/2)}{(\pi v)^{v/2} \Gamma(v/2) |\bm{\Lambda}|^{1/2}}
 \left\{ 1+\tfrac{(u-\mu)^T\bm{\Lambda}^{-1} (u-\mu)}{ v}\right\}^{-\tfrac{v+m}{2}}, v>2.
 \end{align}
\end{proposition}
The above result can also be adopted to construct the confidence regions. Define
\begin{align*}
\tilde{X}_{mk}(\mathbf{U}) \triangleq \left \{ x \, \mid \,  \sqrt{\tfrac{n(n-m)}{ n-1}}  \mathbf{S}_k ^{-1/2}( \bar{\mathbf{x}}_k-x) \in \mathbf{U} \right \},
\end{align*}
where the region $ \mathbf{U}$ is selected such that $\mathbb{P}\big(  T_{n-m}(0,\mathbf{I}_m,m) \in  \mathbf{U}\big) \geq 1-\delta$ with  some $\delta\in(0,1)$.
Similarly to Corollary 2, the  confidence region    $\tilde{X}_{mk}(\mathbf{U}) $   defined above  is asymptotically correct, i.e.,
$\lim\limits_{k \rightarrow \infty}\mathbb{P}\big(x^*\in \tilde{X}_{mk}(\mathbf{U})) \big) = 1-\delta.$ However, it might be harder to construct the confidence regions from the multivariate $t$ distribution than that from the $F$-distribution  since   the region $ \mathbf{U}$   might not be easily obtained.

\section{Numerical Simulations.} \label{sec:simu}

   In this section, we carry out  simulations for    the    parameter estimation problem.   We aim to   estimate  the    unknown  $m $-dimensional parameter $x^*$ based on the  gathered scalar  measurements $ \{d_k\}_{k\geq 1}$ given by
$  d_{k}=u_{k}^Tx^* +\nu_{k},$  where   $u_{k} \in  \mathbb{R}^m$ denotes the regression vector and $\nu_{ k}  \in \mathbb{R}$ denotes   the local observation noise.
Assume that $\{u_k\}$ and $\{\nu_k\}$ are mutually independent   i.i.d.  Gaussian sequences with  distributions $N(\mathbf{0}, R_u)$ and    $N(0,\sigma_{\nu}^2)$, respectively. Suppose   the   covariance matrix $  R_u$    is  positive definite.  Then we might  model the parameter  estimation problem  as   the following stochastic optimization problem:
 \begin{equation}\label{filter1}
  \min_{x\in \mathbb{R}^m } ~  f(x)\triangleq  \mathbb{E}   \big[\| d_{ k}-u_{k}^Tx\|^2\big].
  \end{equation}
Thus,   $f(x)=(x-x^*)^T R_u(x-x^*)+\sigma_{\nu}^2 $ and $\nabla f(x)= R_u(x-x^*).$
 Because the Hessian matrix $R_u$ of the objective function  $f$ is positive  definite, $x^*$ is the unique optimal solution to  \eqref{filter1}.
Suppose that  we can   observe the  regressor  $u_{k}$ and  the measurement $d_{k}$,  then the noisy observation  of the   gradient  $\nabla f(x) $ can  be  constructed as $ u_ku_k^T x-d_k u_k$. Set the dimension of $x^*$ as $m=5.$ We run Algorithm \ref{Alg_1} with  $\alpha={2\over L+\eta}$, Algorithm \ref{Alg_2} with $\alpha={1 \over L}$  and  $\beta = {\sqrt{\kappa} -1 \over\sqrt{\kappa}+1},$ and  Algorithm \ref{Alg_3}    with   $\alpha = { 4 \over (\sqrt{\eta}+\sqrt{L})^2}$ and $\beta =  \left(  { \sqrt{\kappa} -1 \over \sqrt{\kappa} +1}\right)^2,$    where   $x_{0}=y_0=0$ and the batch-size  $N_k=\lceil \rho^{-k}\rceil$  with $\rho={\kappa^2 \over (\kappa+1)^2}$. In the remainder of this section, VSS-SGD, VSS-ACC, VSS-HB represent the abbreviations of Variance-Reduced SGD (Algorithm   \ref{Alg_1}) ,  Variance-Reduced Accelerated SGD (Algorithm   \ref{Alg_2}) ,  Variance-Reduced Heavy-Ball SGD (Algorithm   \ref{Alg_3}),  respectively.

 \begin{figure}[htbp]
 \centering
 \begin{minipage}{0.33\linewidth}
       \centering
  \includegraphics[width=2in]{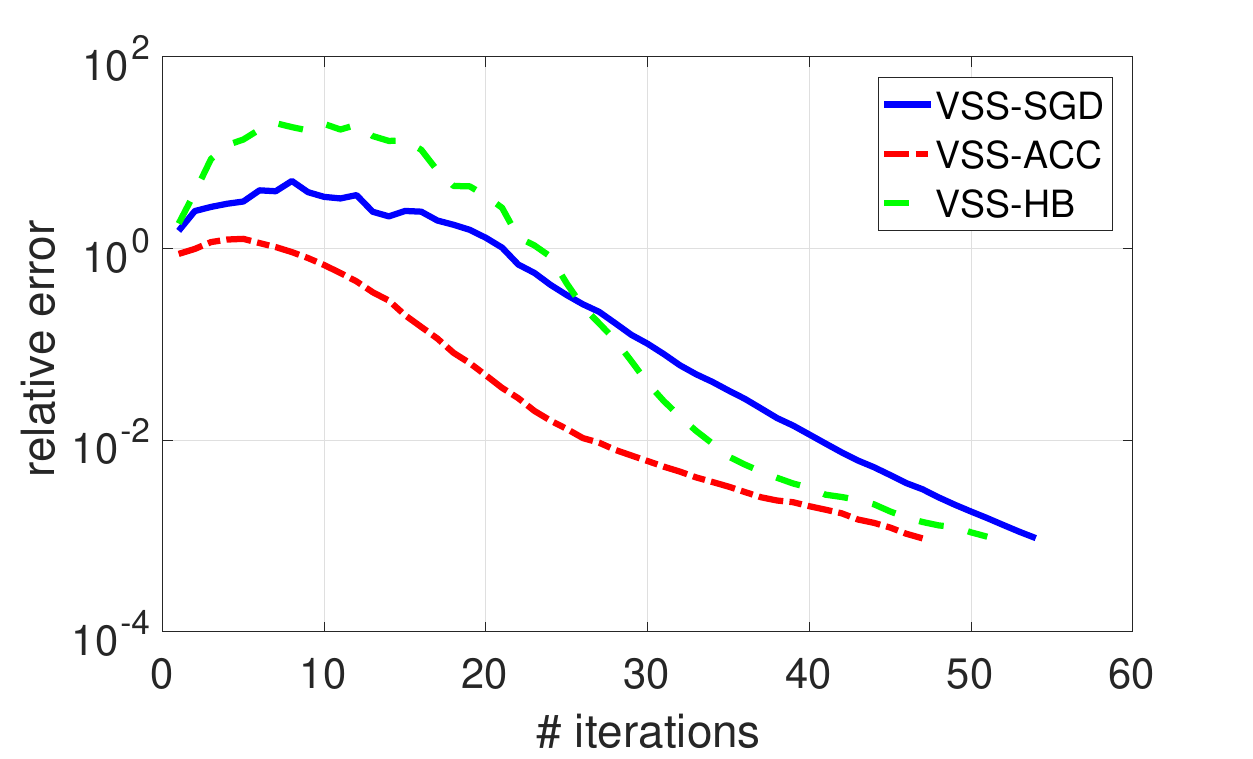}
  \caption{Linear Rate  }
   \label{GOne}
    \end{minipage}%
 \begin{minipage}{0.03\linewidth}
\end{minipage}%
  \begin{minipage}{0.33\linewidth}
       \centering
 \includegraphics[width=2in]{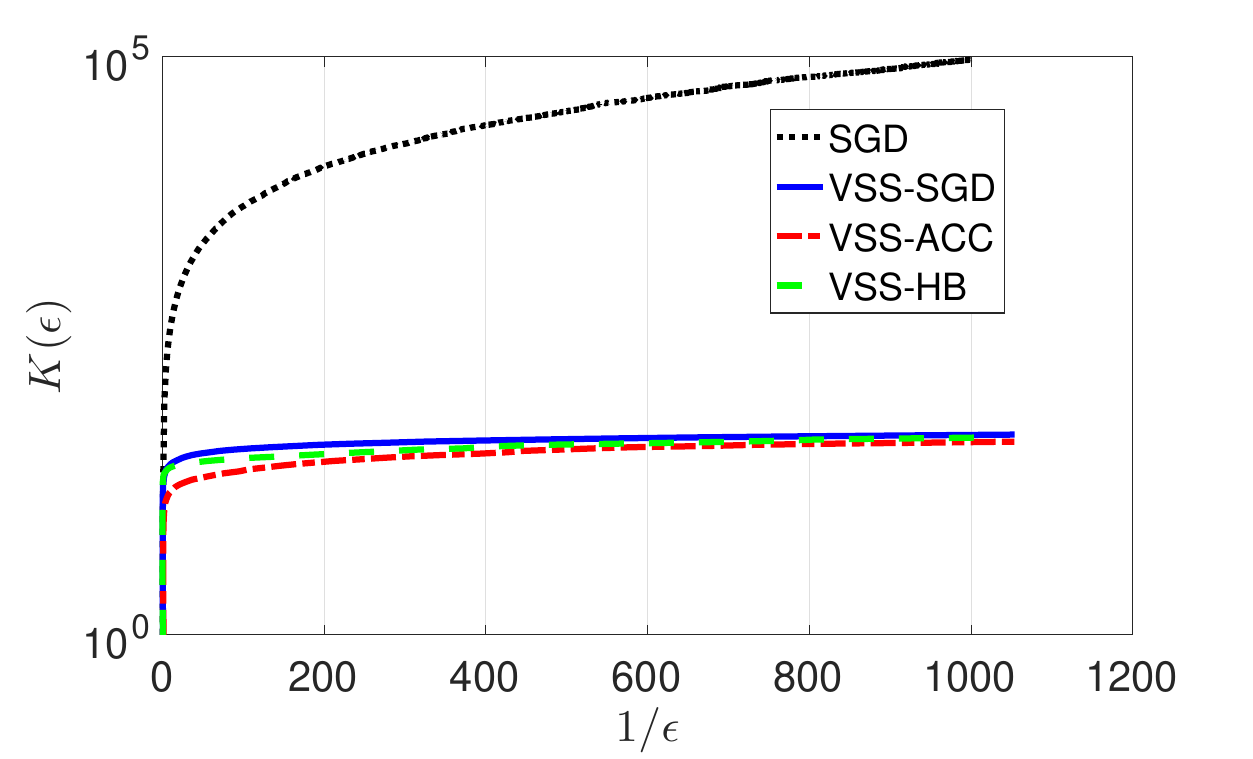}
 \caption{Iteration  Complexity}     \label{Giter}
    \end{minipage}
  \begin{minipage}{0.33\linewidth}
       \centering
 \includegraphics[width=2in]{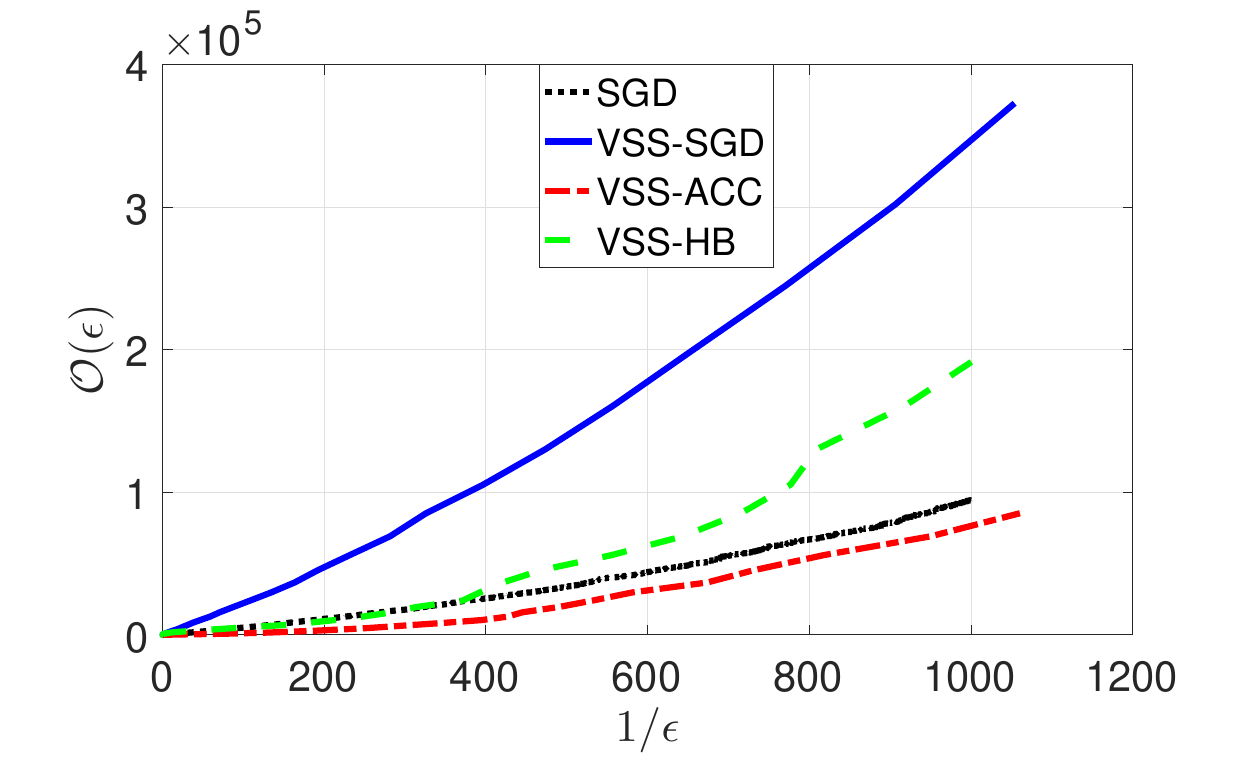}
 \caption{Oracle Complexity}     \label{GTwo}
    \end{minipage}%
\end{figure}

 {\bf Convergence rate, iteration and oracle complexity. } We run
Algorithms  \ref{Alg_1},    \ref{Alg_2},   \ref{Alg_3},  and
the standard SGD algorithm $ x_{k+1}=x_k - \alpha_k  \nabla f(x_k,\xi_{j,k})$
with  $\alpha_k=R_u^{-1}/k$   setting to be the optimal tuning steplength, and
terminate the schemes when   $ { \mathbb{E}[\| x_k-x^*\|_2] \over
\|x^*\|_2}\leq 10^{-3}. $   We then examine their  empirical rate  of
convergence,   iteration and oracle complexity. Here the empirical mean is
calculated by averaging across  100 trajectories. The convergence rate of   the
relative error $ { \mathbb{E}[\| x_k-x^*\|_2] \over \|x^*\|_2} $   is   shown
in Figure \ref{GOne},  which demonstrates   that the iterates   generated by
Algorithms  \ref{Alg_1}-\ref{Alg_3}  converge  in mean to the  optimal solution
at a linear rate. We see that the accelerated  scheme  (Algorithm
\ref{Alg_2}) has the fastest empirical rate,  while the heavy ball method
(Alg.~\ref{Alg_3}) tends to stabilize later in the process.  The empirical relationship between the accuracy
$\epsilon$ and $K(\epsilon)$ is shown in Figure \ref{Giter}, where
$K(\epsilon)$ denotes  the   number of iterations  required   to make  $  {
\mathbb{E}[\| x_k-x^*\|_2] \over \|x^*\|_2} <\epsilon.$ It is seen that  the
standard SGD algorithm requires far more iterations than the proposed  variance reduced schemes  for obtaining  an approximate solution with the same accuracy.  The empirical relationship between $\epsilon$ and $O(\epsilon)$ is shown
in Figure \ref{GTwo}, where  $O(\epsilon)$ denotes  the   number of    sampled
gradients required   to make  $  { \mathbb{E}[\| x_k-x^*\|_2] \over \|x^*\|_2}
<\epsilon.$ We observe that for obtaining an estimate with the same accuracy,
the accelerated  scheme  (Alg.~\ref{Alg_3}) requires the smallest number of sampled gradients,
while  the variable sample-size SGD method (Alg.~\ref{Alg_1})   requires
more  sampled gradients than  standard SGD.

\begin{figure}[htbp]
       \centering
  \includegraphics[width=6in]{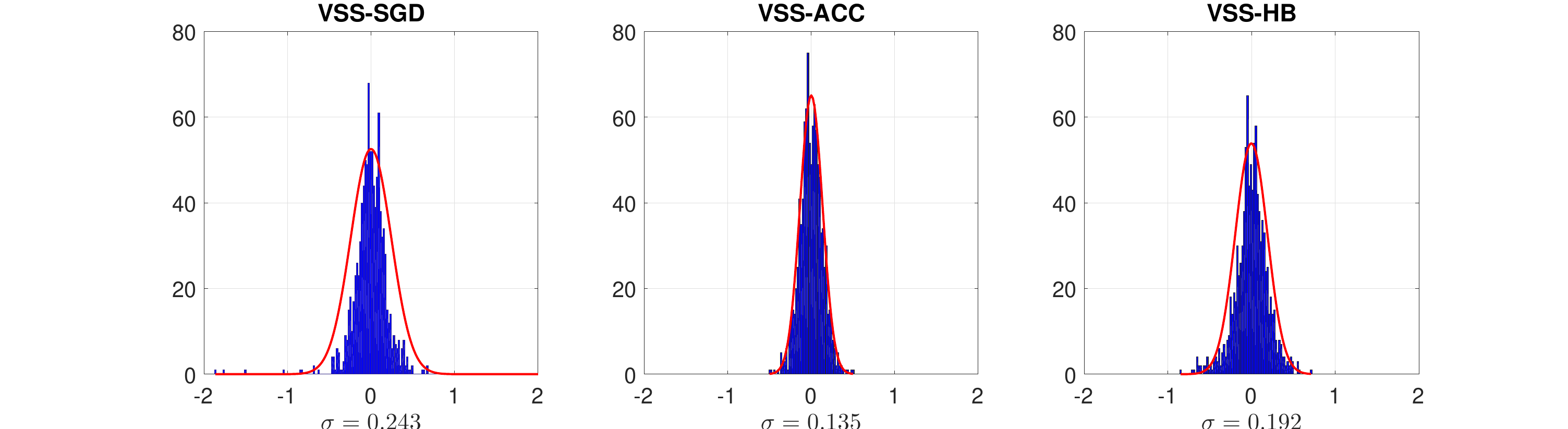}
 \caption{ Histograms of   $\rho^{-k/2}(x_k^5-x^*) $ at $k=50$ along fitted normal distributions }      \label{GThree}
\vspace{-0.2in}
\end{figure}

\begin{figure}[htbp]
\centering
 \includegraphics[width=6in]{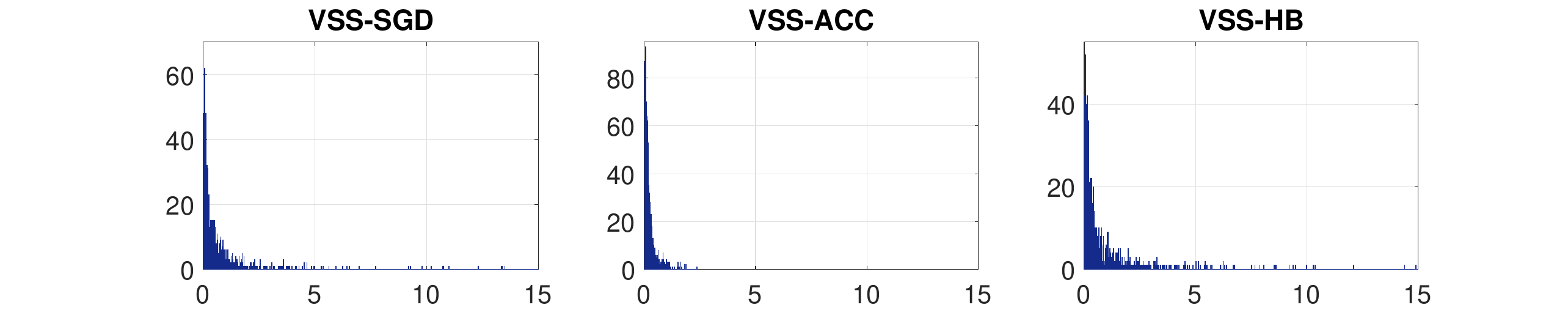}
 \caption{ Histograms of $\rho^{-k}(f(x_k)-f(x^*))$ at $k=50$ }     \label{Ggap}
\end{figure}

  {\bf Limiting  distributions.}   We run $1000$ independent sample paths of      Algorithms  \ref{Alg_1},   \ref{Alg_2}, and   \ref{Alg_3}  and terminate  at $k=50$.  The  empirically obtained  largest eigenvalue of the   covariance matrix $\mathbf{S}_k$ (estimated by \eqref{est-covariance})  at $k=50$
are  $1.648\times 10^{-5} ,  1.649\times 10^{-6} , 1.825\times 10^{-6} $, respectively.     This might  imply that the   accelerated   scheme (Algorithm  \ref{Alg_2}) has the best  performance measured by the covariance of the stationary distribution.  Since the unknown
parameter $x^*$ is multi-dimensional, we merely show  the limiting  distribution of one component of the   rescaled error    $\rho^{-k/2}(x_{ k}-x^*).$ The histograms  of    $\rho^{-k/2}(x_{ k}^5-x^*)$  at    $k=50 $ are shown in   Figure  \ref{GThree}   along with the fitted normal distribution (the red curve), where $x_k^5 $ denotes the fifth component of $x_k.$    It is also seen that the rescaled error is (approximately) normally distributed and among  these, the  accelerated   scheme  has the smallest variance. In addition,    the histograms of the rescaled  suboptimality gap $\rho^{-k}(f(x_k)-f(x^*))$ at $k=50$  are shown in Figure \ref{Ggap}.  We can further conclude that the  empirical sub-optimality gap of the   accelerated   scheme  is the best among  the proposed variable sample-size schemes. We further run  the SGD method  $ x_{k+1}=x_k - R_u^{-1} \nabla f(x_k,\xi_{j,k})/k$.  The histogram  along with the fitted normal distribution     of   $\sqrt{N}(x_N^5-x^*)$     is   displayed in Figure  \ref{GFour},  while the  empirical sub-optimality gap $N(f(x_N)-f(x^*))$  is shown in Figure  \ref{Gsub},  where $N\triangleq \sum_{k=1}^{50} N_k $ denotes the number of sampled gradients utilized by    Algorithms  \ref{Alg_1}-\ref{Alg_3}.

 \begin{figure}[htbp]
  \begin{minipage}{0.56\linewidth}
       \centering
   \includegraphics[width=2.8in]{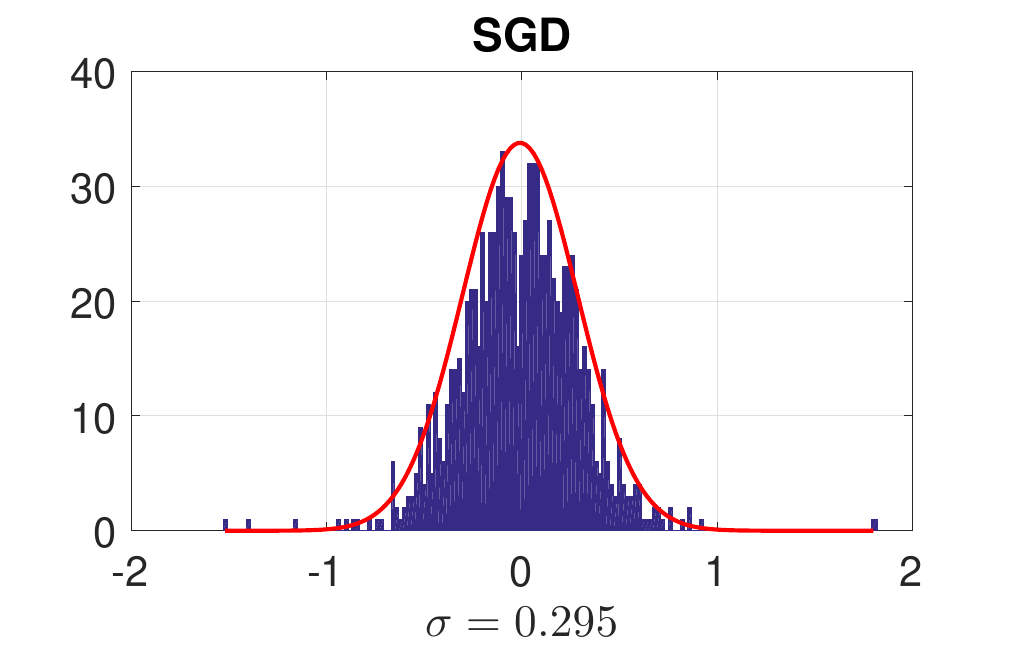}
 \caption{ Hist.    along fitted normal for $\sqrt{N}(x_N^5-x^*)$ }     \label{GFour}
    \end{minipage}%
 \begin{minipage}{0.02\linewidth}
\end{minipage}%
\begin{minipage}{0.4\linewidth}
       \centering
   \includegraphics[height=1.8in,width=2.5in]{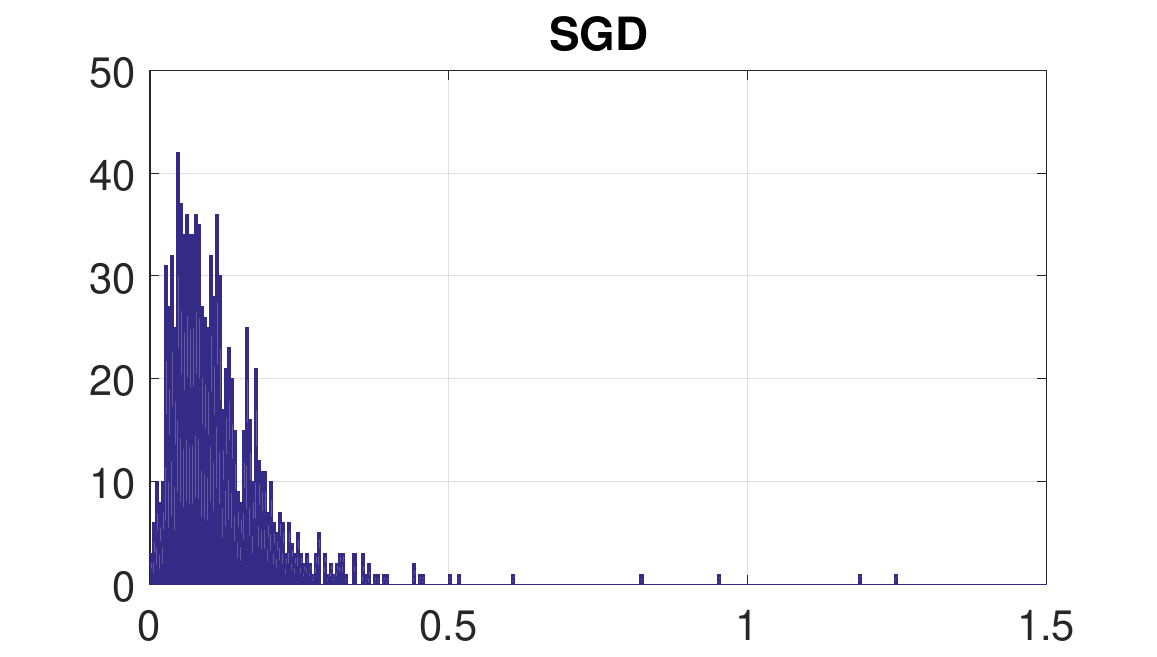}
 \caption{ Hist.     of  $N(f(x_N)-f(x^*))$ }     \label{Gsub}
    \end{minipage}%
\vspace{-0.2in}
\end{figure}

  {\bf Coverage probability of the constructed  confidence region.} In  a single  replication,  we \vus{generate $n$ independent sample paths for}     Algorithms  \ref{Alg_1},    \ref{Alg_2}  and   \ref{Alg_3}   and terminate {each sample path} when the total  number of sampled gradients used reaches $N_{\max}.$
Then we can  construct the $95\%$ confidence region   by  \eqref{def-conf-region} and check whether the true parameter lies in the  constructed  confidence region.  We estimate the  coverage probability (i.e., the proportion of  replications  that the confidence region  contains the true value)  by  conducting    1000 replications.  The estimated coverage probability for  Algorithms \ref{Alg_1}-\ref{Alg_3}  with different selection of  sample paths $n$ and  simulation budget $N_{\max} $ is shown in Table \ref{Tab1}.
It shows the impact of the number of independent runs on the quality of the confidence intervals.  It can be  seen that the coverage probabilities are getting closer to the nominal level of  $95\%$ when the number of  sample paths $n$ grows larger. In particular, we see that $n \geq 10$ seems to produce relatively accurate confidence bands but for $n \leq  8$, the intervals did drop in quality (measured in terms of coverage probability). Thus, we conclude from empirical simulations  that a modest  number of independent runs can generate relatively high quality confidence intervals.
 While it is seen from Table \ref{Tab1} that  the simulation budget $N_{\max}$ does not  significantly impact  the coverage probability; however, a  larger  $N_{\max}$ does lead to narrowing of   the confidence region.


  \begin{table} [htbp]
 \centering
 \small
 \begin{tabular}{|c|c|c|c|c|c|c|c|c|}
 \hline   \multicolumn{2}{|c|}{ } & VSS-SGD &VSS-ACC& VSS-HB  \\  \hline
{\multirow{3}{*}{$N_{\max}=10^3$}}  &  $n=6$ &  $    0.671\pm  0.221$ &$ 0.652	\pm 0.2259$  &$  0.665	\pm
    0.223$ \\   \cline{2-5}
 &  $n=8 $     & $  0.904	\pm  0.0869 $  & $  0.872\pm0.1117$ &$   0.893	\pm   0.0995$    \\ \cline{2-5}
 &  $n=10 $     & $ 0.946	\pm   0.0511$  & $  0.923\pm 0.0591$ &$ 0.951	\pm   0.0466$     \\ \cline{2-5}
 &  $ n=15$    &  $    0.973	\pm   0.0263$ &$  0.966\pm     0.0329$  &$ 0.961\pm   0.0375$    \\ \hline
{\multirow{3}{*}{$N_{\max}=10^4$}}  &  $n=6$ &  $   0.664\pm   0.2233$ &$ 0.62	\pm    0.2358$  &$   0.646 \pm     0.2289 $ \\   \cline{2-5}
 &  $n=8 $     & $  0.872	\pm  0.117 $  & $  0.877\pm  0.108$ &$   0.896	\pm   0.0933$    \\ \cline{2-5}
 &  $n=10 $     & $ 0.929\pm    0.066$  & $      0.944 \pm 0.0592$ &$ 0.949	\pm   0.0484$     \\ \cline{2-5}
 &  $ n=15$    &  $       0.965	\pm     0.0338$ &$      0.946 \pm       0.0511$  &$ 0.961\pm   0.0375$    \\ \hline
{\multirow{3}{*}{$N_{\max}=10^5$}}  &  $n=6$ &  $ 0.651	\pm 0.2274$ &$ 0.656	\pm 0.2259$  &$ 0.664	\pm  0.2233$ \\   \cline{2-5}
 &  $n=8 $     & $ 0.87	\pm 0.1132$  & $ 0.871\pm0.1125$ &$   0.896	\pm   0.0933$    \\ \cline{2-5}
 &  $n=10 $     & $ 0.912	\pm 0.083$  & $ 0.937\pm 0.0591$ &$ 0.919	\pm 0.0745$    \\ \cline{2-5}
 &  $ n=15$    &  $ 0.952	\pm 0.0457$ &$ 0.964	\pm 0.0347$  &$ 0.956	\pm 0.0421$    \\ \hline
 \end{tabular}
  \caption{The estimated coverage probability  of  Algorithms \ref{Alg_1}-\ref{Alg_3}. The ideal coverage probability is 0.95. }  \label{Tab1}
\end{table}

  {\bf Polynomially increasing batch-sizes. } We run Algorithms  \ref{Alg_1}-\ref{Alg_3} with the batch-size  increasing  at a polynomial rate $N_k=\lceil k^{v}\rceil ,v=2 $.  The convergence rate of   the relative error $ { \mathbb{E}[\| x_k-x^*\|_2] \over \|x^*\|_2} $   is   shown in Figure \ref{GPoly1}, while the    histograms  of    $k^{ v/2}(x_{ k}^5-x^*)$  at    $k=100 $  along with the fitted normal distribution are shown in   Figure  \ref{GPoly2}.
It can be seen that the  accelerated  scheme  has the best performance because it displays  fastest convergence rate and the rescaled error has smallest variance.

 \begin{figure}[htbp]
     \begin{minipage}{0.45\linewidth}
       \centering
 \includegraphics[width=2.6in]{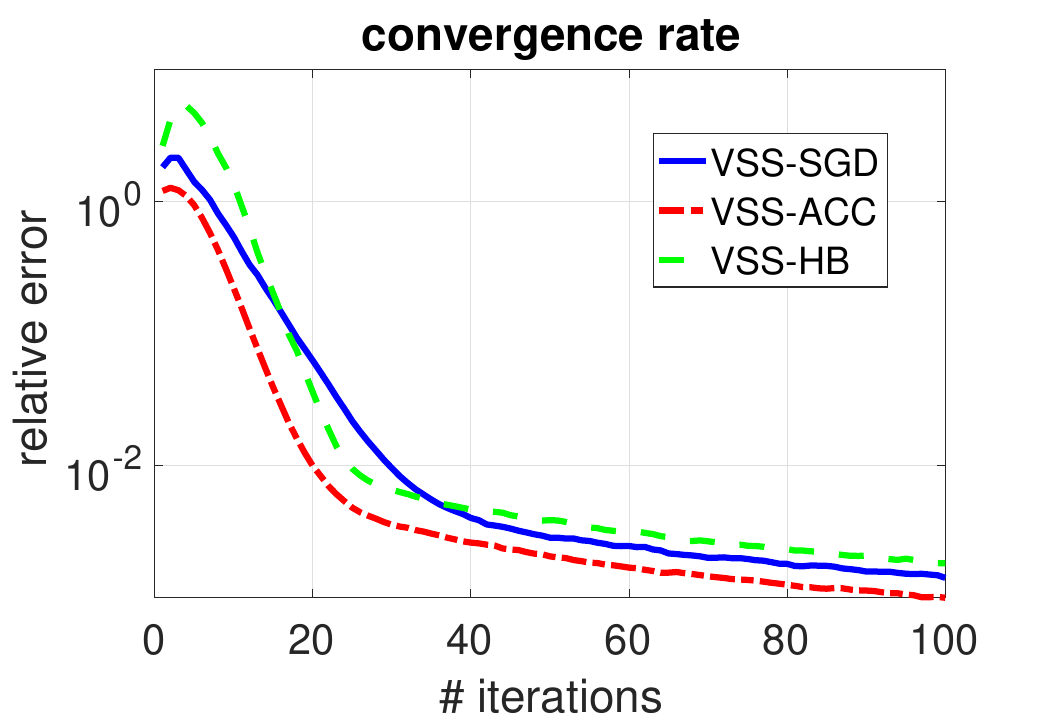}
 \caption{Poly. Batch-size}     \label{GPoly1}
    \end{minipage}%
 \begin{minipage}{0.06\linewidth}
\end{minipage}%
  \begin{minipage}{0.48\linewidth}
       \centering
   \includegraphics[width=2.6in]{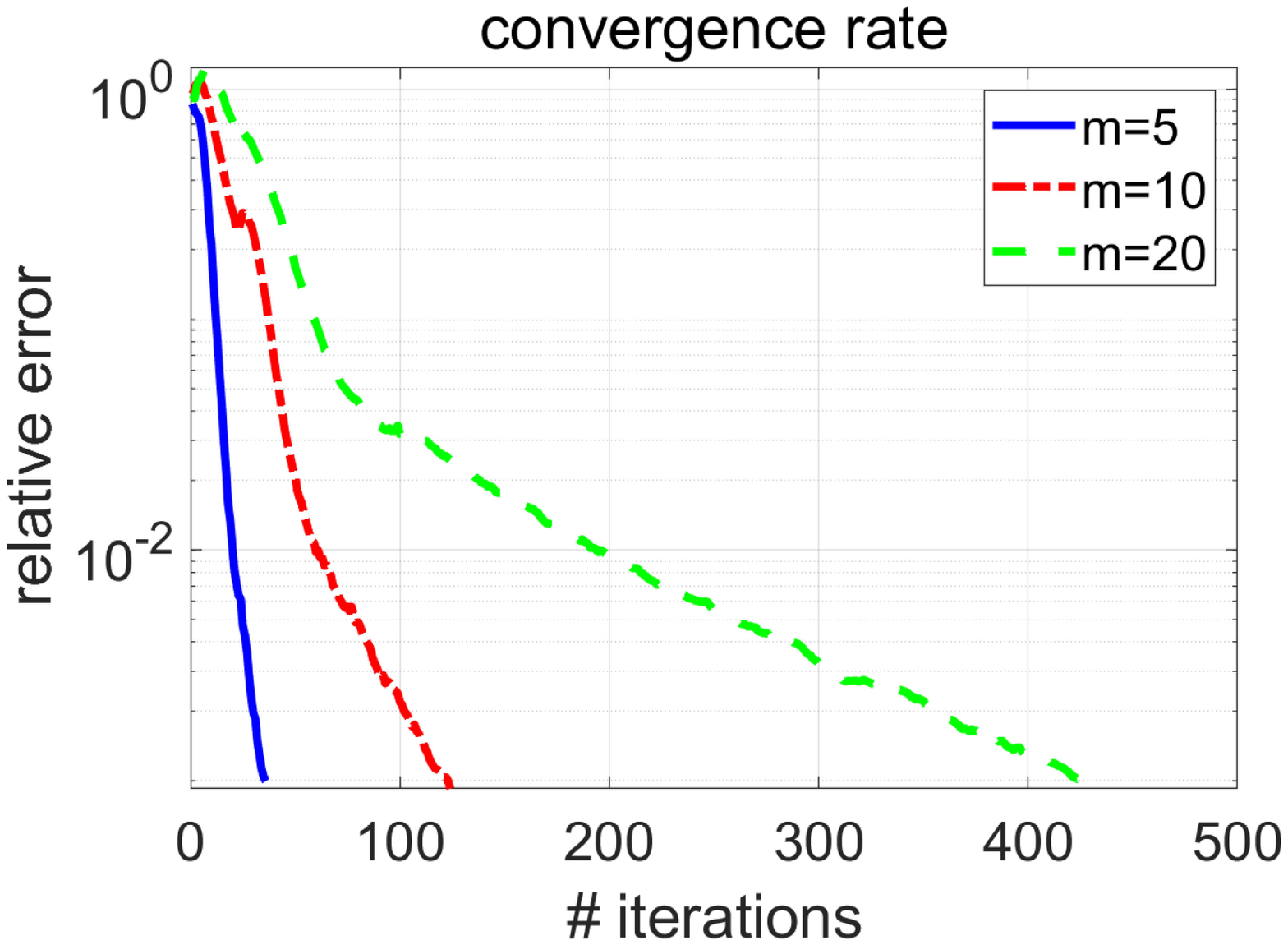}
 \caption{ Convergence rate of VSS-ACC (Alg. \ref{Alg_2}) for different problem dimension $m$}     \label{GDimen}
    \end{minipage}%
\end{figure}

\begin{figure}[htbp]
       \centering
   \includegraphics[width=6in]{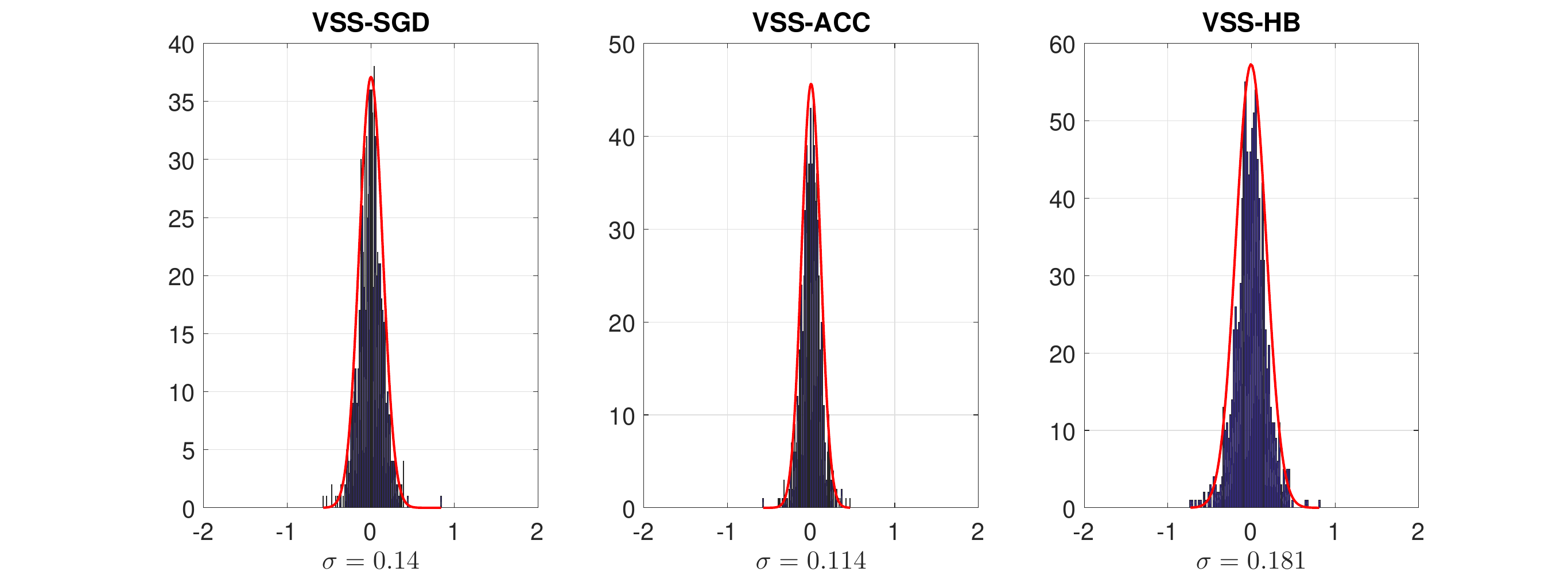}
 \caption{ Histograms        of  $k^{ v/2}(x_{ k}^5-x^*)$  at $k=100$ }     \label{GPoly2}
\end{figure}

 {\bf The effect of the problem dimension.}
We take Algorithm \ref{Alg_2} as an example.
   By adding  {further} numerical simulations,  we examine the impact of $m$  on  algorithm performance.  {In particular, we}  implement    Algorithm \ref{Alg_2} for $m=5,10,20 $, where
  $\alpha={1 \over L},~\beta = {\sqrt{\kappa} -1 \over\sqrt{\kappa}+1},$
    $N_k=\max \big \{\lceil m^{3/4}\rceil,  \lceil \rho^{-k}\rceil \big \}$  with $\rho={\kappa^2\over (\kappa+1)^2 }$, and the  initial values  $x_{0}=y_0=0$.
    The  {empirical counterpart of the } convergence rate of   the relative error $ { \mathbb{E}[\| x_k-x^*\|_2] \over \|x^*\|_2} $   is   shown in Figure \ref{GDimen}, {demonstrating}
     that the larger $m$ leads to a worse rate of convergence.
     We additionally test the coverage probability of the constructed confidence region for  Algorithm \ref{Alg_2} {for different choices of }  $m.$ In a single  replication,  we run     Algorithm    \ref{Alg_2}     for   $n=30$  independent sample paths and terminate them when the total  number of sampled gradients used reaches $5000$.  Then we can  construct the $95\%$ confidence region   by  \eqref{def-conf-region} and check whether the true parameter lies in the  constructed  confidence region.  Finally, we estimate the  coverage probability   by   1000 replications.  The estimated coverage probabilities for
     $m=5,10,20$   are   respectively $       0.956	\pm     0.0421$,  $      0.953 \pm  0.0448     $,  and $ 0.945\pm   0.052$. 

%
%
%

\section{Conclusions}

In this work, we considered the
strongly convex stochastic optimization problem and
proposed three classes of  variance reduced stochastic
gradient algorithms (unaccelerated, accelerated, and
heavy ball methods), where the  unavailable exact gradient
is approximated by an increasing batch of sampled gradients.
We then establish rate and complexity guarantees.
Further, we establish amongst the first formal central limit
theorems for  all the three schemes when the batch-size
increased at either a geometric or a polynomial  rate.  The
covariance matrix specifies how problem   structure
(including the strong convexity parameter,  the Lipschitz
constant, and the Hessian matrix) and distribution  of
gradient   noise influences the algorithm  performance.  In
addition,  we provide an avenue to construct  the
confidence region  of the  optimal solution based on the
central limit theorems. The paper concludes with an application of the  proposed
schemes to the stochastic  parameter estimation problem to
validate our theoretical  findings.  Yet much remains to be understood about the how one may develop analogs of such statements in constrained and nonsmooth regimes.
One avenue for addressing such challenges is via a smoothing framework~\cite{beck17fom} whereby a stochastic gradient update is taken with respect to a smoothed objective and the smoothing parameter is progressively reduced (cf.~\cite{jalilzadeh2018optimal}).
It is  also of interest to explore whether the central limit result on the function value  can be obtained  under   weaker conditions. Finally,  we intend to examine if   confidence regions can be designed based on the  batch-means method~\cite{zhu2019constructing} when  generating independent trials is expensive.

\def\cprime{$'$}

\appendix 

\section{Proof of Lemma \ref{Lem-HB-Q}} \label{App-Lem-HB-Q}

    {First, we introduce an  important relationship between  the  matrix norm  and the spectral  radius of the matrix,
    originally proven  by Gelfand and also shown in
\cite[Lemma 1 in  Section 2.1]{polyak1987introduction}.

\begin{lemma} \label{lemma-Gel} (Gelfand's formula). It holds that
$\rho(\mathbf{P} )=\lim_{k\to \infty}\|\mathbf{P}^k \| ^{1/k}$, i.e., the spectral radius of $\mathbf{P}$  gives the
asymptotic growth rate of  $\| \mathbf{P}^k\|$. Then for any $\iota >0, $ there exists a constant   $c=c(\iota)$ such that
 $\|\mathbf{P}^k \| \leq c(\iota)(\rho(\mathbf{P} )+\iota)^k,\quad \forall k\geq 1.$
\end{lemma}

}

{\bf  Proof of Lemma \ref{Lem-HB-Q}.}
Note that $\nabla f(x)-\nabla f(x^*)=  \mathbf{H}  (x-x^*)$ since $f(x) $ is quadratic and $\nabla^2 f(x)\equiv \mathbf{H}.$ Then  by using  \eqref{Alg3} and  $\nabla f(x^*)=0$, we have
\begin{align*}
x_{k+1}-x^*& =x_k-x^* +\beta(x_k-x_{k-1}) -\alpha     (\nabla f(x_k)-\nabla f(x^*)) -\alpha  w_{k,N_k}
\\& =x_k-x^* +\beta(x_k-x_{k-1}) -\alpha  \mathbf{H}  (x_k-x^*) -\alpha w_{k,N_k}
\\& =\big((1+\beta) \mathbf{I}_m -\alpha   \mathbf{H} \big) (x_k-x^* )  -\beta (x_{k-1}-x^*)  -\alpha w_{k,N_k}.
\end{align*}
Denote by $z_{k}=\begin{pmatrix}
 x_{k}-x^*
\\ x_{k-1}-x^* \end{pmatrix}$. We  then write this  recursion in  matrix form.
\begin{align}\label{hb_recur1}
z_{k+1} &=  \underbrace{\begin{pmatrix} (1+\beta) \mathbf{I}_m -\alpha    \mathbf{H}  &~ -\beta  \mathbf{I}_m\\
\mathbf{I}_m &~ \mathbf{0}_m
\end{pmatrix} }_{ \triangleq \mathbf{T}  } z_k+\begin{pmatrix}  -\alpha   w_{k,N_k} \\ 0\end{pmatrix}  .
\end{align}

By the eigenvalue decomposition, $\mathbf{H}  =\mathbf{U}\bm{\Lambda} \mathbf{U}^T$,  where $\mathbf{U}$ is orthogonal and
 $\bm{\Lambda} \triangleq {\rm diag}\{\lambda_1, \lambda_2,\cdots, \lambda_m\}$ with   $  \lambda_i \in [\eta,  L],i=1,\cdots, m$ being the eigenvalues of  $\mathbf{H} $.
 {Then we can rewrite $\mathbf{T}  $ as $$ \begin{pmatrix}  \mathbf{U} &  {\bf 0}_m \\  {\bf 0}_m  & \mathbf{U}
\end{pmatrix}\begin{pmatrix} (1+\beta) \mathbf{I}_m -\alpha \bm{\Lambda}  &~ -\beta  \mathbf{I}_m\\
\mathbf{I}_m &~ \mathbf{0}_m
\end{pmatrix}\begin{pmatrix}  \mathbf{U} & {\bf 0}_m  \\  {\bf 0}_m  & \mathbf{U}\end{pmatrix}^T.$$
Then the eigenvalues of  the matrix  $\mathbf{T} $ are the roots of the equation  $det\big(vI_{2m}-\mathbf{T}\big)=0$, i.e.,
$$det \begin{pmatrix}
v \mathbf{I}_m-(1+ \beta  )\mathbf{I}_m+\alpha   \bm{\Lambda} & ~\beta \mathbf{I}_m
\\-  \mathbf{I}_m &  v \mathbf{I}_m   \end{pmatrix}=0. $$
By the property $det \begin{pmatrix}
A& B
\\C &  D  \end{pmatrix}=det(AD-BC)$ when  the blocks $A,B,C,D$ are square matrices of the same size and $CD=DC$ \cite{silvester2000determinants},
we have $ det \begin{pmatrix}
v\big(v \mathbf{I}_m-(1+ \beta  )\mathbf{I}_m+\alpha   \bm{\Lambda}  \big) + \beta  \mathbf{I}_m  \end{pmatrix}=0. $
Since the matrix in the determinant  is diagonal,  it can be described by the following  characteristic  equations.
 \[  p_i(v)=(v-  (1+\beta    -\alpha  \lambda_i)) v +\beta=v^2 -(1+\beta    -\alpha  \lambda_i) v+\beta=0 ,\quad i = 1, \cdots, m.\]
}
 From   $ \eta  \leq   \lambda_i  \leq L$ {and $\alpha<4/L$} it follows  that  for any $i=1,\cdots,m:$
 $1-\sqrt{\alpha L} \leq  1-\sqrt{\alpha  \lambda_i } \leq  1-\sqrt{\alpha \eta} $, hence $ | 1-\sqrt{\alpha  \lambda_i }|^2 \leq \max \{|1-\sqrt{\alpha \eta}|^2, |1-\sqrt{\alpha L}|^2\}=\beta<1.$
{Thus,}  the discriminant of the  equation $p_i(v)$,  denoted by $\Delta_i$, is nonpositive.
\begin{align}\label{dictri-Q}
\Delta_i & =(1+\beta    -\alpha  \lambda_i)^2 -4  \beta =(1-\beta)^2-2(1+\beta) \alpha  \lambda_i +(\alpha  \lambda_i)^2\notag
\\&  \leq  (1-| 1-\sqrt{\alpha  \lambda_i }|^2)^2 -2(1+| 1-\sqrt{\alpha  \lambda_i }|^2) \alpha \lambda_i +(\alpha  \lambda_i)^2\notag
\\&  =  (1+ | 1-\sqrt{\alpha  \lambda_i }|^2)^2-4| 1-\sqrt{\alpha  \lambda_i }|^2-2(1+| 1-\sqrt{\alpha  \lambda_i }|^2) \alpha \lambda_i +(\alpha  \lambda_i)^2\notag
 \\&  =  (1-\alpha\lambda_i+ | 1-\sqrt{\alpha  \lambda_i }|^2)^2-4| 1-\sqrt{\alpha  \lambda_i }|^2   = 0.
\end{align}
Hence {$p_i(v)=0$} has two complex roots ${  1+\beta -\alpha    \lambda_i  \pm \mathbf{i} \sqrt{-\Delta_i } \over 2}$,  where $\mathbf{i}=\sqrt{-1}$.  Thus,  the magnitude  of the roots is  $
\sqrt{ \left(\tfrac{  1+\beta -\alpha    \lambda_i}{2} \right)^2 -\tfrac{\Delta_i}{4}}= \sqrt{\beta  }.$
{Then    
 \begin{align}  \label{eigen-alg3}
\rho( \mathbf{T})   = \sqrt{\beta}, \quad \forall k\geq 1.\end{align}
This together with  Lemma \ref{lemma-Gel} implies  that for  any $ \iota\in (0,1-\sqrt{\beta}), $ there exists a constant   $c=c(\iota)$ such that
 \begin{align}  \label{eigen-alg3v2} \|\mathbf{T}^k \| \leq c(\iota)(\sqrt{\beta}+\iota) ^k,\quad \forall k\geq 1.
 \end{align}
 We can recursively write \eqref{hb_recur1} as
 \begin{align}\label{hb_recur2}
z_{k+1} & =  \mathbf{T}^{k+1} z_0-  \sum_{t=0}^{k} \mathbf{T}^{k-t} \alpha   w_{t,N_t} .\notag
\\  \Rightarrow  \| z_{k+1}\|^2 &=  z_0^T (\mathbf{T}^{k+1})^T \mathbf{T}^{k+1} z_0+\alpha^2\left\| \sum_{t=0}^{k} \mathbf{T}^{k-t}    w_{t,N_t} \right\|^2-2 \alpha   (\mathbf{T}^{k+1} z_0)^T \sum_{t=0}^{k} \mathbf{T}^{k-t} w_{t,N_t}.
\end{align}
Note from \eqref{ass-bd-noise} that
\begin{align*}
&\mathbb{E}[w_{i,N_i}^T   w_{j,N_j}]=\mathbb{E}\left[\mathbb{E}[w_{i,N_i}^T   w_{j,N_j} \mid \mathcal{F}_j]\right]=\mathbb{E}\left[w_{i,N_i}^T \mathbb{E}[  w_{j,N_j} |\mathcal{F}_j]\right]=0  {\rm~ for~ any~}  i<j,
 \\&  \mathbb{E}[z_0^T (\mathbf{T}^{k+1})^T \sum_{t=0}^{k} \mathbf{T}^{k-t} w_{t,N_t}]
 = \sum_{t=0}^{k}\mathbb{E}[z_0^T (\mathbf{T}^{k+1})^T  \mathbf{T}^{k-t} w_{t,N_t}] =\sum_{t=0}^{k}\mathbb{E}\left[\mathbb{E}[z_0^T (\mathbf{T}^{k+1})^T  \mathbf{T}^{k-t} w_{t,N_t}\mid \mathcal{F}_t]\right]=0.
 \end{align*}
Then by using  \eqref{ass-bd-noise}, we derive
\begin{align*}
\mathbb{E}\left[\left\|\sum_{t=0}^{k}  \mathbf{T}^{k-t}  w_{t,N_t}\right\|^2\right]& =\sum_{t=0}^{k} \mathbb{E}\left[\|  \mathbf{T}^{k-t}  w_{t,N_t}\|^2\right]
+2\sum_{i<j}  \mathbb{E}\left[  (\mathbf{T}^{k-i}  w_{i,N_i})^T(\mathbf{T}^{k-j}  w_{j,N_j})\right]
\\&=\sum_{t=0}^{k} \mathbb{E}\left[\|  \mathbf{T}^{k-t}  w_{t,N_t}\|^2\right] \leq \sum_{t=0}^{k} \|  \mathbf{T}^{k-t}  \|^2  {\nu^2\over N_t}
\overset{\eqref{eigen-alg3v2}}{\leq }\sum_{t=0}^{k}  c(\iota)^2(\sqrt{\beta}+\iota) ^{2(k-t)}{\nu^2\over N_t}.
\end{align*}
Then by using \eqref{hb_recur2}  and  \eqref{eigen-alg3v2}, we obtain that
\begin{align*}
\mathbb{E}[\| z_{k+1}\|^2] & \leq   \|  \mathbf{T}^{k+1}\|^2 \mathbb{E}[ \| z_0    \|^2] +\alpha^2\mathbb{E}\left[ \left\| \sum_{t=0}^{k}  \mathbf{T}^{k-t} w_{t,N_t} \right\|^2\right]
\\ & \leq  c(\iota)^2(\sqrt{\beta}+\iota)^{2(k+1)}  \mathbb{E}\left[\left\| z_{0}\right\|^2\right]+\alpha^2  \nu^2 c(\iota)^2 \sum_{t=0}^{k} (\sqrt{\beta}+\iota) ^{2(k-t)}{1\over N_t}.
\end{align*}
By recalling  $ \mathbb{E}[\| z_{0}\|^2]= \left [ \left\| \begin{pmatrix}
 x_{0}-x^*
\\ x_{-1}-x^* \end{pmatrix} \right \|^2 \right] =2 \mathbb{E}[\|  x_{0}-x^*\|^2]$ from $x_{-1}=x_0,$ the result follows.}
  \hfill $\blacksquare$

\section{Proof of Lemma \ref{lem-CLT2}}\label{app-lem-CLT2}
Define the sequence $\{u_k\}$ as follows:
\begin{align}\label{recursion-u}
u_{k+1}& =    \mathbf{P}   u_k -\alpha \rho^{-(k+1)/2} \mathbf{G}  w_{k,N_k}, \quad u_0=0.
\end{align}
By combining \eqref{recursion-u}   with  \eqref{recursion-e}, we obtain the following recursion:
 \begin{align}\label{recursion-u'}
e_{k+1}-u_{k+1}& =    \mathbf{P} (e_k-u_k) +\zeta_{k+1}    =   \mathbf{P}^{k+1}(e_0-u_0) +\sum_{t=0}^ k    \mathbf{P}^{k-t} \zeta_{t+1}  .
\end{align}

{Lemma \ref{lemma-Gel}   together with $\rho(\mathbf{P} )   < 1 $ implies that
for any $\varrho\in (\rho(\mathbf{P} ),1)$, there exists a constant  $c(\varrho)$ such that
\begin{align}\label{bd-Pk}
\|\mathbf{P}^k \| \leq c(\varrho)\varrho^k,\quad \forall k\geq 1.
\end{align}Then by using \eqref{bd-Pk}
 and  $\zeta_k \xlongrightarrow [k \rightarrow \infty]{P}  0$,  we obtain  from  \eqref{recursion-u'}  that
\begin{equation*} 
\begin{split}
\| e_{k+1}-u_{k+1}\| & \leq   \|  \mathbf{P}^{k+1} \| \| e_0-u_0\| +\sum_{t=0}^ k   \|  \mathbf{P}^{k-t} \| \|\zeta_{t+1} \|
  \leq c(\varrho)\varrho^{k+1}   \| e_0-u_0\|+c(\varrho) \sum_{t=0}^ k \varrho^{k-t }  o(1) \xlongrightarrow [k \rightarrow \infty]{P}   0.
\end{split}
\end{equation*}
Thus,  $  e_k  $  defined by \eqref{recursion-e}  and  $u_{k}$ produced by \eqref{recursion-u}  converge to the same distribution asymptotically, if indeed the limit exists. 
 This follows from the fact that   $X_k  \xlongrightarrow [ ]{d}  X,~\|X_k -Y_k\| \xlongrightarrow [ ]{P} 0 \Rightarrow Y_k  \xlongrightarrow [ ]{d}  X$ (see \cite[Theorem 2.7(iv)]{van2000asymptotic}).}

From \eqref{recursion-u}  and  $u_0=0$ it follows that
$$u_{k+1}  =    \mathbf{P}  ^{k+1} u_0 -\alpha \sum_{t=0}^ k    \mathbf{P}  ^{k-t}  \mathbf{G} \rho^{-(t+1)/2} w_{t ,N_{t }}
 =     -\alpha \sum_{t=0}^ k    \mathbf{P}  ^{k-t} \mathbf{G}  \rho^{-(t+1)/2} w_{t,N_t}.$$
This implies that
 \begin{align}\label{recursion-u2}
 ~& \alpha^{-1} u_k   =-\sum_{t=0}^ {k-1}    \mathbf{P}  ^{k-1-t}  \mathbf{G} \rho^{-(t+1)/2} w_{t,N_t} =-\sum_{t=1}^ {k}    \mathbf{P}  ^{k-t} \mathbf{G}  \rho^{-t/2} w_{t-1,N_{t-1}} .  \end{align}
  Define $\xi_{kt}\triangleq   - \mathbf{P}  ^{k-t}  \mathbf{G} \rho^{-t/2} w_{t-1,N_{t-1}} $ for any $t: 1\leq t\leq k.$
 We intend to apply Lemma \ref{lem-CLT}. Therefore, we have to check conditions \eqref{cond1}-\eqref{cond3}.  By using  Assumption \ref{ass-fun}(iii), $N_{t }\geq \rho^{-(t+1)}  $, \eqref{ass-bd-noise},  and  \eqref{bd-Pk},  we obtain that  $ \mathbb{E}[\bxi_{kt}| \bxi_{k1},\cdots, \bxi_{k,t-1}]=0$ and
{ \begin{align*}  \mathbb{E}[\|\bxi_{kt}\|^2|  \bxi_{k1},\cdots, \bxi_{k,t-1}] & \leq  \|  \mathbf{P}  ^{k-t}\|^2 \| \mathbf{G}\|^2  \rho^{-t}  \mathbb{E}[\|w_{t-1,N_{t-1}} \|^2| \mathcal{F}_{t-1}]
\\&  \leq c(\varrho)^2
\varrho^{2(k-t)}\| \mathbf{G}\|^2  \rho^{-t}   {\nu^2 \over N_{t-1}}\leq c(\varrho)^2 \nu^2 \| \mathbf{G}\|^2 \varrho^{2(k-t)}  ~a.s..
\\ \Longrightarrow & \sum_{t=1}^k \mathbb{E}[\| \bxi_{kt}\|^2] \leq c(\varrho)^2 \nu^2 \| \mathbf{G}\|^2  \sum_{t=1}^k
\varrho^{2(k-t)} ={c(\varrho)^2\nu^2\varrho^2 \| \mathbf{G}\|^2  \over 1-\varrho^2}<\infty.
\end{align*}} Thus, \eqref{cond1} holds.
   By  {recalling the definitions of  $  \mathbf{S}_{kt}$ and $  \mathbf{R}_{kt}$    in \eqref{def-RS},   using $ \sum_{t=1}^k    \| \mathbf{P}  ^{k-t}\|^2 \leq \sum_{t=1}^k  c(\varrho) ^2 \varrho^{2(k-t)}={ c(\varrho)^2\over 1-\varrho^2} $ from \eqref{bd-Pk} and   \eqref{Ass3}, and   $\rho^{-t}/N_{t-1}  \xlongrightarrow [t \rightarrow \infty] {}1$,}  we have
   \begin{align*}
  &\sum_{t=1}^k  \mathbb{E} [\|  \mathbf{R}_{kt} -\mathbf{S}_{kt}\| ]
\leq \sum_{t=1}^k    \| \mathbf{P}  ^{k-t}  \mathbf{G}\|^2 \mathbb{E}  \left[   \Big\|  N_{t-1}  \mathbb{E}\big[w_{t-1,N_{t-1}}w_{t-1,N_{t-1}}^T| \mathcal{F}_{t-1}\big] - N_{t-1}\mathbb{E}\big[w_{t-1,N_{t-1}}w_{t-1,N_{t-1}}^T\big ]    \Big \|  \right]
  \\& +\sum_{t=1}^k    \| \mathbf{P}  ^{k-t}  \mathbf{G}\|^2 \mathbb{E}  \left[ \Big\|  \left(N_{t-1}  \mathbb{E}\big[w_{t-1,N_{t-1}}w_{t-1,N_{t-1}}^T| \mathcal{F}_{t-1}\big] - N_{t-1}\mathbb{E}\big[w_{t-1,N_{t-1}}w_{t-1,N_{t-1}}^T\big ] \right) \left (\rho^{-t}/N_{t-1}-1\right)   \Big \| \right]
\\&   \leq   \| \mathbf{G}\|^2  \sum_{t=1}^k {c(\varrho) ^2}\varrho^{2(k-t)}  o(1)  \xlongrightarrow [k \rightarrow \infty]{} 0.
   \end{align*}
   This verifies the second equality  in \eqref{cond2}. We now verify the first  equality  in \eqref{cond2}.
\begin{align} \label{recursion-cov}
 \mathbf{S}_k&=\sum_{t=1}^k  \mathbf{S}_{kt}=\sum_{t=1}^k  \mathbf{P}  ^{k-t}  \mathbf{G}  \mathbb{E}\left [ \rho^{-t}w_{t-1,N_{t-1}}w_{t-1,N_{t-1}}^T \right]  (  \mathbf{P}  ^{k-t} \mathbf{G})^T    \notag
\\ & =\sum_{t=1}^k  \mathbf{P}  ^{k-t} \mathbf{G}  \mathbf{S}_0    (  \mathbf{P}  ^{k-t} \mathbf{G})^T+ \sum_{t=1}^k  \mathbf{P}  ^{k-t}  \mathbf{G}\mathbf{S}_0    (  \mathbf{P}  ^{k-t}  \mathbf{G})^T  \left (\rho^{-t}/N_{t-1}-1\right) \\&
\quad  +\sum_{t=1}^k  \mathbf{P}  ^{k-t} \mathbf{G}\left(  N_{t-1}\mathbb{E}\big[w_{t-1,N_{t-1}}w_{t-1,N_{t-1}}^T \big] -\mathbf{S}_0\right)   (  \mathbf{P}  ^{k-t} \mathbf{G})^T  \rho^{-t}  /N_{t-1} , \notag
\end{align}
    where the second and last terms tend  to zero by using   {$ \sum_{t=1}^k    \| \mathbf{P}  ^{k-t}\|^2 \leq { c(\varrho)^2\over 1-\varrho^2} $},   $\rho^{-t}/N_{t-1}  \xlongrightarrow [t \rightarrow \infty] {}1$, and   \eqref{Ass3}.
The   first term on the right-hand side of \eqref{recursion-cov} can be written  as
 \begin{align}\label{sum-s}
 \mathbf{\tilde{S}}_k\triangleq  \sum_{t=1}^k  \mathbf{P}  ^{k-t}  \mathbf{G}\mathbf{S}_0    (  \mathbf{P}  ^{k-t} \mathbf{G})^T = \sum_{t=0}^{k-1}  \mathbf{P}  ^{t} \mathbf{G} \mathbf{S}_0  \mathbf{G}^T ( \mathbf{P}  ^{t}  )^T.
 \end{align}
{By  \eqref{bd-Pk},  $  \sum_{t=0}^{\infty}  \mathbf{P}  ^{t}  \mathbf{G} \mathbf{S}_0   \mathbf{G}^T   \left( \mathbf{P}^t\right)^T $ is well defined and denoted by $\bm{\Sigma}$.
Then  $  \bm{\Sigma}- \mathbf{\tilde{S}}_k = \sum_{t=k}^{\infty}  \mathbf{P}  ^{t}  \mathbf{G} \mathbf{S}_0   \mathbf{G}^T   \left( \mathbf{P}^t\right)^T $ and hence  we derive from \eqref{bd-Pk} that
\begin{align*}
\| \bm{\Sigma}- \mathbf{\tilde{S}}_k\| &\leq  \sum_{t=k}^{\infty}  \| \mathbf{P}  ^{t} \mathbf{G} \mathbf{S}_0  \mathbf{G}^T ( \mathbf{P}  ^{t}  )^T\| \leq \| \mathbf{G} \mathbf{S}_0  \mathbf{G}^T\|   \sum_{t=k}^{\infty}\| \mathbf{P}^t\| ^{2}
\\&\leq \| \mathbf{G} \mathbf{S}_0  \mathbf{G}^T\|  c(\varrho)^2 \sum_{t=k}^{\infty}\varrho^{2t}
=\| \mathbf{G} \mathbf{S}_0  \mathbf{G}^T\|  \tfrac{c(\varrho)^2}{1-\varrho^2}\varrho^{2k} \to 0 {\rm~when~ } k\to \infty.
 \end{align*}
Therefore,  the limit of  \eqref{sum-s}  exists and equals  $\bm{\Sigma} .$ This
  together with \eqref{recursion-cov} and that  the second and last terms on the right hand side of  \eqref{recursion-cov} tend  to zero, we achieve
  $$ \lim_{k\to \infty}  \mathbf{S}_k \triangleq\bm{\Sigma}  = \lim_{k\to \infty}\sum_{t=0}^{\infty}  \mathbf{P}  ^{t}  \mathbf{G} \mathbf{S}_0   \mathbf{G}^T   \left( \mathbf{P}^t\right)^T  .$$ }

 Finally, we have to verify the Lindeberg condition  \eqref{cond3}.
 By  using   $\xi_{kt}= -   \mathbf{P}  ^{k-t}  \mathbf{G} \rho^{-t/2}w_{t-1,N_{t-1}} $, $N_{t-1}=\lceil \rho^{-t} \rceil $, {and  \eqref{bd-Pk}, we obtain that
\begin{align}\label{bd-xi}
 \|\xi_{kt}\|   \leq \|  \mathbf{P}^{k-t}  \| \| \mathbf{G}\|   \rho^{-t/2} \|w_{t-1,N_{t-1}}  \| \leq  c(\varrho) \varrho ^{k-t} \| \mathbf{G}\|  \sqrt{N_{t-1}} \|w_{t-1,N_{t-1}}  \| .
 \end{align}
Hence for any $\epsilon >0,$
\begin{align}\label{bd-xi2} \Big \{ \| \xi_{kt}\| >\epsilon   \Big\}  \subset     \Big \{  \sqrt{N_{t-1}} \|w_{t-1,N_{t-1}}  \|  > \epsilon \| \mathbf{G}\| ^{-1 } c(\varrho)^{-1} \varrho^{-(k-t)}   \Big \}.
\end{align}
Because  for any $t\geq 1$, $\varrho^{-(k-t)} \xlongrightarrow [k \rightarrow \infty] {} \infty$.
Then using \eqref{Ass32}, we obtain that
$$ \sup_{t\geq 1}\mathbb{E} \left [   N_{t-1} \|w_{t-1,N_{t-1}}  \|^2  I_{  [ \sqrt{N_{t-1}} \|w_{t-1,N_{t-1}}  \|  >  \epsilon \| \mathbf{G}\|^{-1 }c(\varrho)^{-1}  \varrho^{-(k-t)}  ]}\right]  \xlongrightarrow [k \rightarrow \infty] {}0.$$
 Consequently, for any $\epsilon>0,$ by using \eqref{bd-xi} and \eqref{bd-xi2}, the following holds:
\begin{equation*}
\begin{split}   & \sum_{t=1}^k \mathbb{E}[\| \bxi_{kt} \|^2I_{[\|\bxi_{kt}\|\geq \epsilon]}]
 \leq   \sum_{t=1}^k\| \mathbf{G}\|^2 c(\varrho)^2 \varrho^{2(k-t)}  \mathbb{E} \left [   N_{t-1} \|w_{t-1,N_{t-1}}  \|^2  I_{  \big [ \sqrt{N_{t-1}} \|w_{t-1,N_{t-1}}  \|  > \epsilon\| \mathbf{G}\|^{-1 } c(\varrho)^{-1}  \varrho^{-(k-t)} \big]}\right]
\\& = \sum_{t=1}^k\| \mathbf{G}\|^2 c(\varrho)^2 \varrho^{2(k-t)}o(1) \xlongrightarrow [k \rightarrow \infty] {} 0 .\end{split}
  \end{equation*}}
  Thus,  the conditions \eqref{cond1}-\eqref{cond3} hold. Then  by using Lemma  \ref{lem-CLT}, the fact that $  e_k  $ and  $u_{k}$ have the same limit distribution, and $\alpha^{-1} u_k   =- \sum_{t=1}^ {k}  \xi_{kt} $,
    we  proves Lemma \ref{lem-CLT2}.
\hfill $\blacksquare$

\section{Proof of Lemma \ref{lem-CLT3}.}\label{app-lem-CLT3}

 We observe that   $\rho(\mathbf{A} )<1$ together with  Lemma \ref{lemma-Gel}  implies that
for any $\varrho\in (\rho(\mathbf{A} ),1)$, there exists a constant $c(\varrho)$ such that
\begin{align}\label{bd-Ak}
\|\mathbf{A}^k \| \leq c(\varrho)\varrho^k,\quad \forall k\geq 1.
\end{align}
   Define $\bm{\Phi}_{k,j} \triangleq \mathbf{A}_k \cdots \mathbf{A}_j$ with $\bm{\Phi}_{j,j+1} \triangleq \mathbf{I}_m$.
   Recall from the definitions $\mathbf{A}_0=\mathbf{A}$ and $ \mathbf{A}_k \triangleq   {\left(k+1 \over k\right)}^{v/2} \mathbf{A} $,
   we obtain that  $ \bm{\Phi}_{k,0}  = \left(k+1 \right)^{v/2} \mathbf{A}^{k } $  and  $ \bm{\Phi}_{k,t+1}  = {\left(k+1 \over t+1\right)}^{v/2} \mathbf{A}^{k -t} $. This   combined   with \eqref{bd-Ak} produces
  \begin{align}
  \|\bm{\Phi}_{k,0}\|& \leq  c(\varrho)  \varrho^k \left(k+1 \right)^{v/2} ,\quad \quad \forall k\geq 0, \label{bd-Phik0}\\
\mbox{ and } \|\bm{\Phi}_{k,t+1}\| &\leq  c(\varrho) {\left(k+1 \over t+1\right)}^{v/2} \varrho^{k-t},\quad \forall k\geq t\geq 0.\label{bd-Phik}
\end{align}
For any given $a>0,$  define  $\tilde{q}\triangleq \varrho ^{a}$ and $\tilde{v}\triangleq av/2.$  Then by using \eqref{bd-Phik} and Lemma \ref{lem-recur},
we derive
 \begin{align}\label{limit-Phi}
& \sum_{t=0}^{k} \|\bm{\Phi}_{k,t+1}\|^a
\leq   c(\varrho)^a \sum_{t=0}^{k}  {\left(k+1 \over t+1\right)}^{av/2} \varrho ^{a(k+1-(t+1)) }
\leq c(\varrho)^a  \left(k+1 \right)^{\tilde{v}} \sum_{t=1}^{k+1}  t^{-\tilde{v}}  \tilde{q} ^{k+1-t} \notag
 \\&  \leq  c(\varrho)^a  \left(k+1 \right)^{\tilde{v}}  \left( \tilde{q}^{k+1 } \tfrac{e^{2\tilde{v}}\tilde{q}^{-1}-1}{ 1-\tilde{q}}+\tfrac{2   (k+1)^{-\tilde{v}}}{ \tilde{q} \ln(1/\tilde{q}) } \right) \notag
  \\& \leq   c(\varrho)^a   \left(  \left(k+1 \right)^{\tilde{v}} \tilde{q}^{k+1 } \tfrac {e^{2\tilde{v}}\tilde{q}^{-1}-1}{ 1-\tilde{q}}+\tfrac{2  }{ \tilde{q} \ln(1/\tilde{q}) } \right)
  \leq    c(\varrho)^a   \left( c_{\tilde{q},\tilde{v}} \tfrac{e^{2\tilde{v}}\tilde{q}^{-1}-1}{ 1-\tilde{q}}+\tfrac{2  }{ \tilde{q} \ln(1/\tilde{q}) } \right)    ,\quad \forall k\geq 1. \end{align}
Define an  auxiliary   sequence $\{u_k\}$ by
 \begin{align}\label{recursion-poly-u}
u_{k+1}& =    \mathbf{A}_k  u_k - \alpha (k+1)^{v/2} \mathbf{G} w_{k,N_k}, \quad u_0=0.
\end{align}
This combined with  \eqref{poly-recursion4}    produces     the following recursion:
 \begin{align*}
e_{k+1}-u_{k+1}& =   \mathbf{A}_k(e_k-u_k) +\zeta_{k+1}
 \\ 	&   =\bm{\Phi}_{k,0} (e_0-u_0) +\sum_{t=0}^ k   \bm{\Phi}_{k,t+1} \mathbf{G}  \zeta_{t+1}=\bm{\Phi}_{k,0}e_0 +\sum_{t=0}^ k   \bm{\Phi}_{k,t+1}   \zeta_{t+1}.
\end{align*}
Then by using  \eqref{bd-Phik0}, \eqref{limit-Phi}, $\mathbb{E} [\|e_0\|^2]<\infty,$ and $\zeta_{k+1} \xlongrightarrow [k \rightarrow \infty]{P} 0 $, we conclude that
 \begin{align*}
\| e_{k+1}-u_{k+1}\|& =   \| \bm{\Phi}_{k,0} \| \| e_0 \| +\sum_{t=0}^ k  \|  \bm{\Phi}_{k,t+1}\|  \| \zeta_{t+1} \|\xlongrightarrow [k \rightarrow \infty]{P} 0  .
\end{align*}
This implies that  $e_k$ defined by   \eqref{poly-recursion4}   and    $u_{k}$ defined as in \eqref{recursion-poly-u} have the same limit distribution if exists.
Thus, it remains to find the stationary distribution of $u_{k}$.

 From \eqref{recursion-poly-u} it follows  that
$$ u_{k+1}=\bm{\Phi}_{k,0}u_0- \alpha\sum_{t=0}^ k   \bm{\Phi}_{k,t+1}\mathbf{G}  (t+1)^{v/2}w_{t,N_t} =- \alpha\sum_{t=0}^ k   \bm{\Phi}_{k,t+1}\mathbf{G}  (t+1)^{v/2}w_{t,N_t}.$$
 This implies that
\begin{align} \label{def-poly-u}   \alpha^{-1} u_k =- \alpha\sum_{t=0}^ {k-1}   \bm{\Phi}_{k-1,t+1}\mathbf{G}   (t+1)^{v/2} w_{t,N_t}  =-\sum_{t=1}^ {k}  \bm{\Phi}_{k-1,t} \mathbf{G} ~ t^{v/2} w_{t-1,N_{t-1}}   . \end{align}

 We intend to apply Lemma \ref{lem-CLT} by  defining $\xi_{kt}\triangleq   -\bm{\Phi}_{k-1,t}\mathbf{G}   t^{v/2} w_{t-1,N_{t-1}}  $  and check conditions  \eqref{cond1},  \eqref{cond2}, and  \eqref{cond3}.
   Using $N_{t-1}   \triangleq \lceil t^v\rceil \geq t^v$,   \eqref{limit-Phi}, and Assumption \ref{ass-fun}(iii),  we can verify  \eqref{cond1}.
 Also, using \eqref{Ass3} and \eqref{limit-Phi}, the definitions of $  \mathbf{S}_{kt}$ and $  \mathbf{R}_{kt}$ in \eqref{def-RS}, the second equality  of  \eqref{cond2} holds. We now verify the first  equality  in \eqref{cond2}.
\begin{align} \label{poly-recursion-cov}
 \mathbf{S}_k&=\sum_{t=1}^k  \mathbf{S}_{kt}=\sum_{t=1}^k \bm{\Phi}_{k-1,t} \mathbf{G} \mathbb{E}\left [t^{v}  w_{t-1,N_{t-1}}  w_{t-1,N_{t-1}} ^T \right]  \mathbf{G}^T   \bm{\Phi}_{k-1,t}^T  \notag
\\ & =\sum_{t=1}^k \bm{\Phi}_{k-1,t} \mathbf{G}  \mathbf{S}_0\mathbf{G} ^T \bm{\Phi}_{k-1,t}^T + \sum_{t=1}^k \bm{\Phi}_{k-1,t} \mathbf{G} \mathbf{S}_0  \mathbf{G} ^T \bm{\Phi}_{k-1,t}^T  \left (t^v/N_t-1\right) \\&
\quad  +\sum_{t=1}^k  \bm{\Phi}_{k-1,t}\mathbf{G}  \left(  N_{t-1}\mathbb{E}\left [w_{t-1,N_{t-1}}  w_{t-1,N_{t-1}} ^T \right] -\mathbf{S}_0\right) \mathbf{G}^T \bm{\Phi}_{k-1,t}^T t^v /N_{t-1} , \notag
\end{align}
    where the second and last terms tend  to zero by    \eqref{Ass3}, $t^v/N_{t-1}  \xrightarrow [t \rightarrow \infty] {}1$, and  $ \sum_{t=1}^k    \|\bm{\Phi}_{k-1,t}\|^2 <\infty$ from \eqref{limit-Phi}.
   While  the    first term on the right-hand side of \eqref{poly-recursion-cov},  by using  \eqref{limit-Phi} one obtains
   \begin{align*}
&  \sum_{t=1}^k \| \bm{\Phi}_{k-1,t} \mathbf{G} \mathbf{S}_0 \mathbf{G}^T \bm{\Phi}_{k-1,t}^T  \|\leq \|  \mathbf{G} \mathbf{S}_0 \mathbf{G}^T\|  \sum_{t=1}^k \| \bm{\Phi}_{k-1,t}\|^2
<\infty..
 \end{align*}
Because  $ \sum_{t=1}^k \| \bm{\Phi}_{k-1,t} \mathbf{G} \mathbf{S}_0 \mathbf{G}^T \bm{\Phi}_{k-1,t}^T  \|$ is   monotonically   increasing and bounded,   its limit   exists.  Thus, the limit of  $  \sum_{t=1}^k \bm{\Phi}_{k-1,t}  \mathbf{G} \mathbf{S}_0 \mathbf{G}^T \bm{\Phi}_{k-1,t}^T= \sum_{t=1}^k  {\left(k  \over  t\right)}^{v} \mathbf{A}^{k-t}   \mathbf{G} \mathbf{S}_0 \mathbf{G}^T (\mathbf{A}^{k-t} )^T $  exists and is denoted by $\bm{\Sigma}  .$
Then $ \mathbf{S}_k$ as defined as in  \eqref{poly-recursion-cov}
satisfies  that $$   \lim\limits_{k\to \infty}\mathbf{S}_k = \lim\limits_{k\to \infty}\sum_{t=1}^k  {\left(k  \over  t\right)}^{v} \mathbf{A}^{k-t}  \mathbf{G} \mathbf{S}_0 \mathbf{G}^T (\mathbf{A}^T )^{k-t} \triangleq \bm{\Sigma} .$$ Thus, the first equality  in \eqref{cond2} holds.

 Finally,  the Lindeberg condition  \eqref{cond3} can similarly   validated as  that of Lemma \ref{lem-CLT2}.     Then all conditions of Lemma  \ref{lem-CLT} hold. Thus, by  Lemma \ref{lem-CLT} and   \eqref{def-poly-u},  we conclude that  $\alpha^{-1} u_k   \xlongrightarrow [k \rightarrow \infty]{d}   N(0,\bm{\Sigma})$.  Because     $e_k$ defined by   \eqref{poly-recursion4}     and  $u_{k}$  defined as in \eqref{recursion-poly-u}  have the same limit distribution,  Lemma \ref{lem-CLT3}  is then proved.  \hfill $\blacksquare$

 \section{Proof of Proposition \ref{cf-clt}} \label{app-cf-clt}
  Based on Theorem \ref{thm-CLT1} it is seen that
\begin{align}\label{clt} e_{ik} \triangleq \alpha^{-1} \rho_1^{-k/2}( x_{ik}-x^*)  \xlongrightarrow [k \rightarrow \infty]{d}   Y_i \sim N(0,\bm{\Sigma}_1),\quad  i=1,\cdots, n,\end{align}
where $Y_1,\cdots, Y_n$   are  $n$   i.i.d.  random vectors with distribution $N(0,\bm{\Sigma}_1)$. Define
\begin{align}\label{def-es}
\bar{ \mathbf{e}}  \triangleq{1\over n} \sum_{i=1}^n  Y_i {\rm~and~}
\mathcal{\bm{S}}   \triangleq  {1\over n-1}  \sum_{i=1}^n\big (Y_i-
\bar{ \mathbf{e}}  \big)\big (Y_i- \bar{ \mathbf{e}}  \big)^T  .
\end{align} Recall from \cite[Theorem 4]{bodnar2016singular}  that  $ \bar{\mathbf{e}} $ and  $\mathcal{\bm{S}} $ are independently distributed,  and
\begin{align}\label{est-es}  \bar{ \mathbf{e}}  \sim N(0,\bm{\Sigma}_1/n) {\rm~~and ~~}(n-1) \mathcal{\bm{S}}   \sim \mathcal{\bm{W}}_m(\bm{\Sigma}_1,n-1),
\end{align}
where $ \mathcal{\bm{W}}_m(\bm{\Sigma}_1,n-1)$ denotes the  $m$-dimensional Wishart distribution with $n-1$  degrees of freedom and the matrix parameter   $\bm{\Sigma}_1.$  Because $\bm{\Sigma}_1$ is invertible and $n\geq m+1$, the random  matrix  $\mathcal{\bm{S}} $ is
almost surely invertible.

 (i)   Denote by    $\mathbf{e}=col\{e_1,\cdots,e_n\} \triangleq (e_1^T,\cdots, e_n^T)^T \in \mathbb{R}^{mn}$ with $e_i \in \mathbb{R}^m,~i=1,\cdots, n$.  Note that    $g_1(\mathbf{e} )\triangleq  {1\over n}\sum_{i=1}^n e_i  $   is a continuous function.
 Since $e_{ik},~i=1,\cdots,n$ are mutually independent,    from \eqref{clt}  it follows that
\begin{align}\label{lim-concated-e}
col\{e_{1k},\cdots,e_{nk}\} \xlongrightarrow [k \rightarrow \infty]{d}   col\{ Y_1,\cdots, Y_n\}  .
\end{align} Then  by the continuous mapping theorem \cite[Theorem 1.14]{dasgupta2008asymptotic} and \eqref{est-es},  we derive
\begin{align}\label{clt2} &\sqrt{n} \alpha^{-1} \rho_1^{-k/2}( \bar{\mathbf{x}}_k-x^*)   =
\sqrt{n} {1\over n} \sum_{i=1}^n e_{ik}  \xlongrightarrow [k \rightarrow \infty]{d} \sqrt{n} {1\over n} \sum_{i=1}^n  Y_i =\sqrt{n} \bar{ \mathbf{e}}    \sim N(0,\bm{\Sigma}_1).
\end{align}

(ii)
 Since   $rank(\bm{\Sigma}_1)=m$ and $\sqrt{n} \bar{ \mathbf{e}}    \sim N(0,\bm{\Sigma}_1) $  by \eqref{clt2}, we obtain
 \begin{align}\label{chi-square2}
 U_1\triangleq  (\sqrt{n} \bar{ \mathbf{e}})^T\bm{\Sigma}_1^{-1} (\sqrt{n} \bar{ \mathbf{e}}) \sim  \chi^2(m)  , \end{align}  where $ \chi ^2( m)$ denotes the   chi-squared distribution with $m$ degrees of freedom. Because $rank(\bm{\Sigma}_1)=m,$   $ \bar{\mathbf{e}} \in \mathbb{R}^m$ and  $\mathcal{\bm{S}} $ are independently distributed with $ \bar{\mathbf{e}}^T\bm{\Sigma}_1$ being non-zero with probability one,  from    \cite[Corollary 1]{bodnar2016singular} it follows that
\begin{align} \label{chi-square}
U_2\triangleq  { \bar{ \mathbf{e}}^T\bm{\Sigma}_1^{-1}   \bar{ \mathbf{e}}   \over \bar{ \mathbf{e}} ^T \left( (n-1) \mathcal{\bm{S}}  \right)^{-1}  \bar{ \mathbf{e}}}\sim  \chi ^2(n-m)
\end{align}
 is  independent of $ \bar{ \mathbf{e}}  $.   Hence $U_1$ and $U_2$ are independent.

From  \eqref{clt} and \eqref{est-covariance}  it follows that
\begin{equation}\label{est-covariance2}
\begin{split}
 \bar{\mathbf{x}}_k-x^*&={1\over n}\sum_{i=1}^n  ( x_{ik}-x^*)=\alpha\rho_1^{k/2} \sum_{i=1}^n e_{ik},
 \\   \mathbf{S}_k&={1\over n-1}\sum_{i=1}^n (x_{ik}-\bar{\mathbf{x}}_k)
(x_{ik}-\bar{\mathbf{x}}_k)^T={ \alpha^2\rho_1^{k}\over n-1}  \sum_{i=1}^n \big ( e_{ik}-{1\over n}\sum_{i=1}^n e_{ik} \big)    \big ( e_{ik}-{1\over n}\sum_{i=1}^n e_{ik} \big).
\end{split}
\end{equation}
 Note that  $g_2(\mathbf{e})   \triangleq    {1\over n-1}  \sum_{i=1}^n\left ( e_i-{1\over n}\sum_{i=1}^n e_i \right)\left ( e_i-{1\over n}\sum_{i=1}^n e_i \right)^T  $ is a continuous function.  Because the matrix inverse functional is continuous in
a neighborhood of any non-singular matrix, and  $$ \sum_{i=1}^n\left ( Y_i-{1\over n}\sum_{i=1}^n Y_i \right)\left ( Y_i-{1\over n}\sum_{i=1}^n Y_i \right)^T$$  is almost surely invertible from  \eqref{est-es}, we conclude that  $  \big(g_2(\mathbf{e}) \big)^{-1}  $  is almost surely continuous in a neighborhood of   $col\{ Y_1,\cdots, Y_n\}.$   Hence,  $g(\mathbf{e}) =g_1(\mathbf{e})^T \big(g_2(\mathbf{e}) \big)^{-1} g_1(\mathbf{e})$  is almost surely continuous in a neighborhood of   $col\{Y_1,\cdots, Y_n\}.$   Therefore,    by the continuous mapping theorem \cite[Theorem 1.14]{dasgupta2008asymptotic},  and \eqref{lim-concated-e},    we have
\begin{equation}\label{cmt1}
\begin{split} & n ( \bar{\mathbf{x}}_k-x^*)^T\mathbf{S}_k^{-1}   ( \bar{\mathbf{x}}_k-x^*)  \\&
\overset{\eqref{est-covariance2}}{=}n \left({\sum_{i=1}^n e_{ik}\over n}\right)^T
\Big({1\over n-1} \sum_{i=1}^n \big ( e_{ik}-{1\over n}\sum_{i=1}^n e_{ik} \big)    \big ( e_{ik}-{1\over n}\sum_{i=1}^n e_{ik} \big) ^T\Big)^{-1} {\sum_{i=1}^n e_{ik}\over n}
 \\& \xlongrightarrow [k \rightarrow \infty]{d}
   n \tfrac{ \sum_{i=1}^nY_i}{ n}
    \left(\tfrac{1}{n-1}\sum_{i=1}^n\left (Y_i-\tfrac{\sum_{i=1}^nY_i  }{ n} \right)\left (Y_i-\tfrac{\sum_{i=1}^nY_i  }{ n}  \right)^T \right)^{-1} \tfrac{ \sum_{i=1}^nY_i}{ n}
  \\& \overset{\eqref{def-es}}{=}n(n-1) \bar{ \mathbf{e}}^T \left( (n-1) \mathcal{\bm{S}}  \right)^{-1}  \bar{ \mathbf{e}}
=   (n-1){ n\bar{ \mathbf{e}}^T\bm{\Sigma}_1^{-1}   \bar{ \mathbf{e}}   \over   \tfrac{ \bar{ \mathbf{e}}^T\bm{\Sigma}_1^{-1}   \bar{ \mathbf{e}} }{ \bar{ \mathbf{e}} ^T \big( (n-1) \mathcal{\bm{S}}  \big)^{-1}  \bar{ \mathbf{e}}}}
\\& = (n-1) {U_1 \over U_2} \sim {m(n-1) \over n-m} F(m,n-m),
\end{split}
\end{equation}
where the last one holds because $U_1=n\bar{ \mathbf{e}}^T\bm{\Sigma}_1^{-1}   \bar{ \mathbf{e}} \sim  \chi^2(m)$  by \eqref{chi-square2},  $ U_2={ \bar{ \mathbf{e}}^T\bm{\Sigma}_1^{-1}   \bar{ \mathbf{e}}   \over \bar{ \mathbf{e}} ^T \left( (n-1) \mathcal{\bm{S}}  \right)^{-1}  \bar{ \mathbf{e}}} \sim  \chi ^2(n-m) $   by \eqref{chi-square},   $U_1 $ and $U_2$ are independent. Thus,  $F(d_1, d_2) $ arises as the ratio of two appropriately scaled chi-squared variates \cite{degroot2012probability}.
 \footnote{ If the random variables $U_1\sim \chi^2(d_1)$ and $U_2\sim \chi^2(d_2)$  are independent, then ${U_1 /d_1 \over U_2 /d_2} \sim F(d_1, d_2).$}
\hfill $\blacksquare$

\section{Proof of Proposition \ref{cf-clt2}} \label{app-cf-clt2}
%

Since $\bm{\Sigma}_1$ is symmetric and invertible,   by the property of Wishart distribution, we have that
\begin{align*}
 \bm{\Sigma}_1^{-1}(n-1) \mathcal{\bm{S}}\bm{\Sigma}_1^{-1}  \overset{ \eqref{est-es}}{\sim} \bm{\Sigma}_1^{-1} \mathcal{\bm{W}}_m(\bm{\Sigma}_1,m+(n-m)-1)\bm{\Sigma}_1^{-1}
\sim  \mathcal{\bm{W}}_m(\bm{\Sigma}_1^{-1},m+(n-m)-1).
\end{align*}
Note by \eqref{clt2} that $ \sqrt{n(n-m)} {\bm{\Sigma}_1}^{-1/2} \bar{ \mathbf{e}} \sim N(0,(n-m)\mathbf{I}_m).$ Recall from \cite[Theorem 4]{bodnar2016singular}  that  $ \bar{\mathbf{e}} $ and  $\mathcal{\bm{S}} $ are independent.
Then by \cite[Representation B]{lin1972some},
we conclude that
\begin{align*}
\left((n-1)\bm{\Sigma}_1^{-1} \mathcal{\bm{S}}\bm{\Sigma}_1^{-1} \right)^{-1/2} \sqrt{n(n-m)} {\bm{\Sigma}_1}^{-1/2} \bar{ \mathbf{e}} \sim T_{n-m}(0,\bm{\Sigma}_1  ,m). \end{align*}
Thus, by the probability   density  function of  $ T_v(\mu,\bm{\Lambda},m)$  defined in \eqref{def-prob}, we see that
 \begin{align*}
{\bm{\Sigma}_1}^{-1/2}  \left((n-1)\bm{\Sigma}_1^{-1} \mathcal{\bm{S}}\bm{\Sigma}_1^{-1} \right)^{-1/2} \sqrt{n(n-m)} {\bm{\Sigma}_1}^{-1/2} \bar{ \mathbf{e}} \sim T_{n-m}(0, \mathbf{I}_m ,m). \end{align*}
Hence
\begin{align}\label{Tdis}
\left((n-1) \mathcal{\bm{S}}  \right)^{-1/2} \sqrt{n(n-m)}  \bar{ \mathbf{e}} \sim T_{n-m}(0,\mathbf{I}_m  ,m). \end{align}

Similarly to the derivation of \eqref{cmt1}, we can show that $ \big(g_2(\mathbf{e}) \big)^{-1/2} g_1(\mathbf{e})$  is almost surely continuous in a neighborhood of   $col\{Y_1,\cdots, Y_n\}.$ Then  by  the continuous mapping theorem
\cite[Theorem 1.14]{dasgupta2008asymptotic}   and \eqref{lim-concated-e},    we
achieve  \eqref{t-dis} by the following.
\begin{align*}
& ( (n-1) \mathbf{S}_k)^{-1/2}    \sqrt{n(n-m) }  ( \bar{\mathbf{x}}_k-x^*)
\\& \overset{\eqref{est-covariance2}}{=}
\Big((n-1){1\over n-1} \sum_{i=1}^n \big ( e_{ik}-{1\over n}\sum_{i=1}^n e_{ik} \big)    \big ( e_{ik}-{1\over n}\sum_{i=1}^n e_{ik} \big) ^T\Big)^{-1/2}\sqrt{n(n-m) } {\sum_{i=1}^n e_{ik}\over n}
 \\& \xlongrightarrow [k \rightarrow \infty]{d}
    \left((n-1) \tfrac{1}{n-1}\sum_{i=1}^n\left (Y_i-\tfrac{\sum_{i=1}^nY_i  }{ n} \right)\left (Y_i-\tfrac{\sum_{i=1}^nY_i  }{ n}  \right)^T \right)^{-1/2} \sqrt{n(n-m) } \tfrac{ \sum_{i=1}^nY_i}{ n}
 \\& \overset{\eqref{def-es}}{=} \left((n-1)\mathcal{\bm{S}}\right) ^{-1/2} \sqrt{n(n-m) } \bar{ \mathbf{e}}  \overset{\eqref{Tdis}}{\sim} T_{n-m}(0,\mathbf{I}_m,m).
 \qquad\qquad\qquad \qquad\qquad\qquad \blacksquare
\end{align*}

\newpage

 \noindent {\bf SUPPLEMENTARY MATERIAL.}  \\

\renewcommand\thesection{\Alph{section}}

\numberwithin{equation}{section}

The  proof of Lemma \ref{lem-acc} is motivated by \cite[Section 2.2]{nesterov2013introductory} and \cite[Section 3.6.2]{bubeck2014theory}.
\section{Proof of Lemma \ref{lem-acc} } \label{app-lem-acc}
  We   define $\phi_k(x)$ and $p_k$ as follows.
\begin{align} \label{def-ph1}& \phi_0(x)=f(x_0)+{\eta \over 2} \|x-x_0\|^2,
\\  \label{def-ph2}&\phi_{k+1}(x) =(1-\gamma) \phi_k(x) +\gamma \left(f(x_{k })   +(x-x_k)^T h(x_k)+{\eta \over 2} \|x-x_k\|^2\right), \\&
p_{k+1} =  (1-\gamma ) p_k + \left( \alpha +{(1-\gamma)\gamma \over 2\eta}\right) \|w_{k,N_k}\|^2 +  \alpha  w_{k,N_k} ^T \nabla f(x_k)  , ~ p_0=0, \label{def-p}
\end{align}
where $h(x_k) \triangleq {x_k-y_{k+1} \over \alpha} =\nabla f(x_k)  + w_{k,N_k}.$

We first show  by induction that  for any $ k\geq 0$,  $ \nabla^2 \phi_{k }(x)= \eta \mathbf{I}_m. $
  By \eqref{def-ph1} it is seen that $ \nabla^2  \phi_0(x)= \eta \mathbf{I}_m$. Suppose $ \nabla^2 \phi_{k }(x)= \eta \mathbf{I}_m $, then  by \eqref{def-ph2}, we obtain that
$$  \nabla^2 \phi_{k+1}(x)=(1-\gamma)   \nabla^2 \phi_{k }(x) +\gamma \eta \mathbf{I}_m= \eta \mathbf{I}_m.$$
Thus, $ \nabla^2 \phi_{k }(x)= \eta \mathbf{I}_m$ for any $ k\geq 0$.  Because  $\phi_k$  is a quadratic function,   we can rewrite $\phi_k(x)$  as
\begin{align} \label{can-form} \phi_{k }(x)=\phi_k^* + {\eta \over 2} \|x-v_k\|^2  {~\rm with~} v_k \triangleq \argmin_x \phi_{k }(x), \quad \forall k\geq 0.\end{align}

We proceed to give a recursive form for $v_{k+1}$ and $\phi_{k+1}^*$. Noting from  \eqref{can-form}  that
$\nabla  \phi_k(x)=\eta (x-v_k).$
Then by   using the first-order optimality condition  $\nabla  \phi_{k+1}(x) =0$ of the unconstrained  convex optimization $\min_x \phi_{k+1}(x)$,   and
the   definition of  $\phi_{k+1}(x)$ in  \eqref{def-ph2}, we obtain that
$$ \nabla  \phi_{k+1}(x)  =(1-\gamma) \eta (x-v_k)+\gamma h(x_k) + \gamma  \eta (x-x_k)=0 .$$
 By rearranging the previous equation, we have that
 \begin{align}
 &  v_{k+1}=(1-\gamma)   v_k+ \gamma  x_k -\gamma h(x_k) /\eta . \label{def-v}
\end{align}
By  using  \eqref{def-ph2} and  \eqref{can-form},     evaluating  $\phi_{k+1}(x)$ at $x=x_k$  we obtain that
\begin{equation}\label{exp-phi}\begin{split}
 \phi_{k+1}^*  &  = \phi_{k +1}(x_k) -{\eta \over 2} \|x_{k }-v_{k+1}\|^2
 = (1-\gamma ) \phi_k(x_k) +\gamma  f(x_{k })  -{\eta \over 2} \|x_{k }-v_{k+1}\|^2
\\&= (1-\gamma)    \phi^*_k+{\eta  (1-\gamma)        \over 2} \|x_k-v_k\|^2 +\gamma  f(x_{k })  -{\eta \over 2} \| v_{k+1}-x_{k }\|^2 .
\end{split} \end{equation}
Note by \eqref{def-v} that
\begin{align*}
   \|v_{k+1}-x_k\| ^2 &=  \| (1-\gamma)   (v_k-x_k) -\gamma h(x_k) /\eta\|^2
\\&  =(1-\gamma) ^2 \|   v_k-x_k\|^2 + { \gamma^2 \over \eta^2}  \| h(x_k) \|^2- {2 \gamma  (1-\gamma) \over \eta}  (v_k-x_k)^Th(x_k) .
\end{align*}
This together  with \eqref{exp-phi} leads to
\begin{align}\label{rec-phi}
  \phi_{k+1}^* &  =  (1-\gamma)    \phi^*_k+\gamma  f(x_{k })+{\eta \gamma (1-\gamma)        \over 2} \|x_k-v_k\|^2
        - { \gamma^2 \over 2\eta }  \|h(x_k) \|^2  + \gamma  (1-\gamma)   (v_k-x_k)^Th(x_k) .
\end{align}

We  then  show by induction that the following holds  for any $k \geq 0$.
\begin{align}\label{re-vx}
v_k-x_k={1\over \gamma} (x_k-y_k) .
\end{align}
From  \eqref{def-ph1}   it is seen that  the minimizer of $\phi_0$   is  $v_0=x_0$. Then by the initial condition $x_0=y_0$,  we conclude that \eqref{re-vx} holds for $k=0.$  We inductively assume that \eqref{re-vx} holds for $k$, and proceed to prove that \eqref{re-vx} holds for $k+1$. By  substituting   $v_k=x_k+(x_k-y_k)/\gamma $ into \eqref{def-v},  one obtains
$v_{k+1}    =  (1-\gamma)  (x_k+(x_k-y_k)/\gamma  )+ \gamma  x_k -\gamma h(x_k)  /\eta. $ Hence
\begin{align*}
v_{k+1}-x_{k+1} &    =   {1 \over \gamma}  \left( x_k -\gamma ^2  h(x_k)  /\eta\right)  -\left({1 \over \gamma} -1 \right) y_k -x_{k+1}   .
\end{align*}
This together with $h(x_k) = {x_k-y_{k+1} \over \alpha}$ and  $\gamma =\sqrt{ \alpha \eta }$ produces
\begin{align*}
v_{k+1}-x_{k+1} &    =   {1 \over \gamma}  \left( x_k - \alpha h(x_k)   \right)  -\left({1 \over \gamma} -1 \right) y_k -x_{k+1}
\\& {{\tiny \eqref{Alg22}} \atop = }{1 \over \gamma}    y_{k+1}-\left({1 \over \gamma} -1 \right) { (1+\beta) y_{k+1}-x_{k+1 } \over \beta}  -x_{k+1} ={1\over \gamma} (x_{k+1}-y_{k+1})  ,
\end{align*}
where the last equality holds by $\beta = {1 -\gamma \over 1+ \gamma} $.  Thus, we have shows that \eqref{re-vx} holds for any $k\geq 0$.

Then by substituting\eqref{re-vx} into \eqref{rec-phi}, and using  $  h(x_k) =\nabla f(x_k)+ w_{k,N_k},$   we obtain that
\begin{align}\label{rec-phi2}
  \phi_{k+1}^*   & = (1-\gamma)    \phi^*_k+\gamma  f(x_{k })+{\eta  (1-\gamma)        \over 2\gamma} \|x_k-y_k\|^2   - { \gamma^2 \over 2\eta }  \| h(x_k) \|^2     +  (1-\gamma)   (x_k-y_k)^Th(x_k)\notag
\\& = (1-\gamma)    \phi^*_k+\gamma  f(x_{k })      - { \gamma^2 \over 2\eta }  \| h(x_k) \|^2     +  (1-\gamma)   (x_k-y_k)^T \nabla f(x_k) \notag
\\& \quad +{\eta  (1-\gamma)        \over 2\gamma} \|x_k-y_k\|^2 + (1-\gamma)   (x_k-y_k)^T   w_{k,N_k}
\\&  \geq (1-\gamma)    \phi^*_k+\gamma  f(x_{k })      - { \gamma^2 \over 2\eta }  \| h(x_k) \|^2   +  (1-\gamma)   (x_k-y_k)^T \nabla f(x_k)  -{(1-\gamma)\gamma \over 2\eta} \| w_{k,N_k}\|^2, \notag
\end{align}
where the last inequality follows by ${\eta  (1-\gamma)        \over 2\gamma} \|a\|^2+ (1-\gamma) a^Tb\geq -{(1-\gamma)\gamma \over 2\eta}\|b\|^2.$

We proceed to show that $ \phi_k^* \geq f(y_k)-p_k$ for any $k\geq 0$. By the definitions  \eqref{def-ph1} and \eqref{def-p},  we see  that $ \phi_0^* =f(x_0)=f(y_0)$ and  $p_0=0 $. Hence $ \phi_k^* \geq f(y_k)-p_k$  holds for $k=0.$ We inductively assume that $   f(y_k) \leq \phi_k^*   +p_k$,  and  aim to show that  $  f(y_{k+1})-  \phi_{k+1}^*  \leq p_{k+1}.$
Since $f$  is $L$-smooth, by using  $h(x_k)={x_k-y_{k+1} \over \alpha} =\nabla f(x_k)  + w_{k,N_k},$ we obtain that
\begin{align}\label{bd-fk}
f(y_{k+1} )&\leq f(x_k) +(y_{k+1}-x_k)^T \nabla f(x_k) +{ L \over 2} \|y_{k+1}-x_k\|^2 \notag
\\  &\leq f(x_k)  -\alpha h(x_k)^T(h(x_k)- w_{k,N_k})+{ L\alpha^2 \over 2} \|h(x_k) \|^2 \notag
\\& \leq    f(x_k)  +\left({ L\alpha^2 \over 2}-\alpha \right)  \| h(x_k)\|^2+   \alpha \|w_{k,N_k}\|^2 +  \alpha  w_{k,N_k} ^T \nabla f(x_k),
\end{align}
where the last inequality holds because   $h(x_k)^T w_{k,N_k}=    \|w_{k,N_k}\|^2 +   w_{k,N_k} ^T \nabla f(x_k).$
By the induction assumption  $   f(y_k) \leq \phi_k^*   +p_k$  and  the convexity of  $f$,   we obtain that
\begin{align*}
 f(x_k) & = (1-\gamma ) f(y_k) + (1-\gamma) \left( f(x_k) -f(y_k) \right)+\gamma   f(x_k)
\\&   \leq  (1-\gamma ) \phi_k^*   + (1-\gamma ) p_k+(1-\gamma) (x_k-y_k)^T \nabla  f(x_k)  + \gamma   f(x_k)  .
\end{align*}
The above bound combined with \eqref{bd-fk}  produces
\begin{align*}
f(y_{k+1} ) &   \leq   (1-\gamma ) \phi_k^*+ \gamma   f(x_k)    +(1-\gamma) (x_k-y_k)^T  \nabla  f(x_k) + (1-\gamma ) p_k \\&
  \quad + \left({ L\alpha^2 \over 2}-\alpha \right)   \| h(x_k)\|^2+   \alpha \|w_{k,N_k}\|^2 +  \alpha  w_{k,N_k} ^T \nabla f(x_k).
\end{align*}
It incorporated  with  \eqref{rec-phi2} leads to the following relation:
\begin{align}\label{rec-phi3}
   f(y_{k+1} )   -\phi_{k+1}^*   \leq & \left({ L\alpha^2 \over 2}-\alpha+{ \gamma^2 \over 2\eta }  \right)   \| h(x_k) \|^2+ (1-\gamma ) p_k +  \alpha  w_{k,N_k} ^T \nabla f(x_k) \notag
    \\&   + \left( \alpha +{(1-\gamma)\gamma \over 2\eta}\right) \|w_{k,N_k}\|^2    \leq  p_{k+1},
\end{align}
where the last inequality holds by the definition of $p_{k+1}$ in \eqref{def-p}
and ${ L\alpha^2 \over 2}-\alpha+{ \gamma^2 \over 2\eta }={ \alpha ( \alpha L -1)  \over 2} \leq 0 $
 from  $\alpha\leq 1/L$  and $\gamma^2=\alpha \eta <1$.
 Therefore, we conclude  that $f(y_k) \leq \phi_k^*+p_k$ for any $k\geq 0.$

Because $ f$ is $\eta$-strongly convex, from  \eqref{def-ph2} and $\mathbb{E}[ h(x_k) | \mathcal{F}_k]= \nabla f(x_k)$ it follows that
\begin{align*}
&\mathbb{E}[\phi_{k+1}(x) | \mathcal{F}_k]   =  (1-\gamma) \phi_k(x)   +\gamma \left(f(x_{k })   +(x-x_k)^T  \nabla f(x_k)+{\eta \over 2} \|x-x_k\|^2\right)  \\& \qquad \qquad \qquad   \leq  (1-\gamma) \phi_k(x)  +\gamma f(x),\quad \forall x \in \mathbb{R}^m.
\end{align*}
By taking unconditional expectations, we obtain that $\mathbb{E}[\phi_{k+1}(x) ]   \leq  (1-\gamma) \mathbb{E}[\phi_k(x)]  +\gamma f(x)$ for any $x \in \mathbb{R}^m.$
Therefore, by rearranging terms and setting $x=x^*$ in the above inequality, we have
\begin{align*}
\mathbb{E}[\phi_{k+1}(x^*)]-f(x^*)& \leq (1-\gamma) ( \mathbb{E}[\phi_k(x^*)]- f(x^*))
 \leq (1-\gamma)^{k+1}  ( \mathbb{E}[\phi_0(x^*)]- f(x^*)) .
\end{align*}
Then by the fact  that $f(y_{k }) \leq \phi_{k }^* +p_{k}$, there holds
\begin{align}\label{re-fy}
\mathbb{E}[f(y_{k}) ]-f(x^*) & \leq \mathbb{E}[\phi_{k}(x^*)]-f(x^*)  +\mathbb{E}[p_{k}]  \leq  (1-\gamma)^{k }  ( \mathbb{E}[\phi_0(x^*)]- f(x^*)) +\mathbb{E}[p_{k }] \notag
\\& \leq  {\eta+L \over 2} (1-\gamma)^{k }    \mathbb{E}[\| x_0-x^*\|^2] +\mathbb{E}[p_{k }],
\end{align}
where the last inequality holds because $\phi_0(x^*)-f(x^*)=f(x_0)-f(x^*)+{\eta \over 2} \|x^*-x_0\|^2  \leq {\eta+L \over 2}\|x^*-x_0\|^2$ by the $L$-smoothness of  $f$  and $ \nabla f(x^*)=0.$
Next, we derive a bound on  $\mathbb{E}[p_{k}]$. By taking expectations on both sides of \eqref{def-p},
 using   $p_0=0$ and \eqref{ass-bd-noise}, we obtain that
\begin{align*}
\mathbb{E}[p_{k}] & =  (1-\gamma )\mathbb{E}[p_{k -1}] + \left( \alpha +{(1-\gamma)\gamma \over 2\eta}\right)  \mathbb{E}[\|w_{ k-1 ,N_{k-1 }}\|^2 ]
\\&  =\left( \alpha +{(1-\gamma)\gamma \over 2\eta}\right)  \sum_{i=0}^{k-1}  (1-\gamma )^i \mathbb{E}[\|w_{ k-i-1 ,N_{k-i -1}}\|^2 ]
   \leq  \nu^2 \left( \alpha +{(1-\gamma)\gamma \over 2\eta}\right)  \sum_{i=0}^{k-1}  (1-\gamma )^i/N_{k-1-i}.
\end{align*}
This together with \eqref{re-fy}   proves Lemma \ref{lem-acc}. \hfill $\blacksquare$

\section{Proof of Lemma \ref{Lem-HB}} \label{App-Lem-HB}
The proof of Lemma \ref{Lem-HB} is motivated by \cite[Theorem 2]{ghadimi2015global}.

Define
\begin{align}\label{def-pk}
p_k={\beta \over 1-\beta} (x_k-x_{k-1}).
\end{align}
Then by \eqref{Alg3}, we have
\begin{align*}
x_{k+1}+p_{k+1}& =x_{k+1}+{\beta \over 1-\beta} (x_{k+1}-x_{k })={1\over 1-\beta} x_{k+1} -{\beta \over 1-\beta}  x_{k }  \notag
\\&={x_k -\alpha     (\nabla f(x_k)+w_{k,N_k}) +\beta(x_k-x_{k-1})\over 1-\beta} x_{k+1} -{\beta \over 1-\beta}  x_{k }
\\& =x_k+p_k-{\alpha \over 1-\beta}    (\nabla f(x_k)+w_{k,N_k}) \notag.
\end{align*}
This implies that
\begin{align}\label{rec-xp}
\| x_{k+1}+p_{k+1}-x^*\|^2\notag
&= \|x_k+p_k-x^*\|^2-{2\alpha \over 1-\beta}    (\nabla f(x_k)+w_{k,N_k})^T(x_k+p_k-x^*)
\\&+ {\alpha^2 \over (1-\beta)^2}  \| \nabla f(x_k)+w_{k,N_k}\|^2.
\end{align}

With the definition  of  $\mathcal{F}_k$ and the update rule \eqref{Alg3}, we see that $x_k$ is adapted to $\mathcal{F}_k$, and hence
$p_k$ is  adapted to $\mathcal{F}_k$. Then by taking conditional expectations on both sides of \eqref{rec-xp} on $\mathcal{F}_k$,
 we obtain
 \begin{align}\label{rec-xp2}
& \mathbb{E}[ \| x_{k+1}+p_{k+1}-x^*\|^2|\mathcal{F}_k ] = \|x_k+p_k-x^*\|^2-{2\alpha \over 1-\beta}    (\nabla f(x_k)+ \mathbb{E}[w_{k,N_k}|\mathcal{F}_k])^T(x_k+p_k-x^*)\notag
\\&
+ {\alpha^2 \over (1-\beta)^2} \left( \| \nabla f(x_k)\|^2 +2\nabla f(x_k) ^T  \mathbb{E}[w_{k,N_k} | \mathcal{F}_k ]+  \mathbb{E}[\|w_{k,N_k}\|^2| \mathcal{F}_k ]\right) \notag
\\ & \overset{  \eqref{ass-bd-noise}}{ \leq}  \|x_k+p_k-x^*\|^2-{2\alpha \over 1-\beta}    \nabla f(x_k)^T(x_k+p_k-x^*)
+ {\alpha^2 \over (1-\beta)^2} \left( \| \nabla f(x_k)\|^2  + \tfrac{\nu^2}{N_k}\right)
\\& \overset{\eqref{def-pk}}{=} \|x_k+p_k-x^*\|^2-\tfrac{2\alpha}{ 1-\beta}    \nabla f(x_k)^T(x_k -x^*)-\tfrac{2\alpha \beta }{ (1-\beta)^2}    \nabla f(x_k)^T(x_k-x_{k-1})
+ \tfrac{\alpha^2 }{(1-\beta)^2} \left( \| \nabla f(x_k)\|^2  + \tfrac{\nu^2}{N_k}\right). \notag
\end{align}

Since  $ \nabla f(x^*)=0, $ Assumption \ref{ass-fun}(i) and (ii) hold, we  recall from \cite[(2.1.24)]{nesterov2013introductory} that
\begin{align}\label{bd-Hb1}
  \nabla f(x_k)^T(x_k -x^*) \geq  \tfrac{\eta L}{ L+\eta} \|x_k -x^* \|^2+\tfrac{1 }{ L+\eta} \|  \nabla f(x_k)\|^2.
\end{align}
Since $f$  is $\eta$-strongly convex from Assumption \ref{ass-fun}(ii), we have
\begin{align}\label{bd-Hb2}
  & f(x_{k-1})- f(x_k) \geq   \nabla f(x_k)^T (x_{k-1}-x_k) +\tfrac{\eta}{ 2}  \|x_{k-1}-x_k\|^2. \notag
   \\& \Rightarrow \nabla f(x_k)^T(x_k-x_{k-1}) \geq  f(x_k) -f(x_{k-1}) +\tfrac{\eta}{ 2}   \|x_{k-1}-x_k\|^2.
\end{align}
 Then by substituting   \eqref{bd-Hb1}  and    \eqref{bd-Hb2}  into \eqref{rec-xp2}, we  obtain
 \begin{align*}
& \mathbb{E}[ \| x_{k+1}+p_{k+1}-x^*\|^2|\mathcal{F}_k ]  \leq  \|x_k+p_k-x^*\|^2+ \tfrac{\alpha^2 }{(1-\beta)^2} \left( \| \nabla f(x_k)\|^2  + \tfrac{\nu^2}{N_k}\right)
\\&-\tfrac{2\alpha}{ 1-\beta}  \left(\tfrac{\eta L}{ L+\eta} \|x_k -x^* \|^2+\tfrac{1 }{ L+\eta} \|  \nabla f(x_k)\|^2\right)-\tfrac{2\alpha \beta }{ (1-\beta)^2}   \left( f(x_k) -f(x_{k-1}) +\tfrac{\eta}{ 2}   \|x_{k-1}-x_k\|^2\right)
. \notag
\end{align*}
By taking unconditional expectations on both sides of the above equation, and rearranging the terms, we    get
 \begin{align}\label{rec-xp3}
& \mathbb{E}[ \| x_{k+1}+p_{k+1}-x^*\|^2|\mathcal{F}_k ] +\tfrac{2\alpha \beta }{ (1-\beta)^2} ( f(x_k)-f(x^*)) \leq  \|x_k+p_k-x^*\|^2+\tfrac{2\alpha \beta }{ (1-\beta)^2}( f(x_{k-1})-f(x^*)) \notag
\\&+ \tfrac{\alpha^2 \nu^2}{(1-\beta)^2N_k}
-\tfrac{2\alpha\eta L}{ (1-\beta)(L+\eta)} \|x_k -x^* \|^2 - \tfrac{\alpha \beta \eta}{ (1-\beta)^2}      \|x_{k-1}-x_k\|^2+
 \tfrac{\alpha  }{ 1-\beta } \left( \tfrac{\alpha  }{ 1-\beta }-\tfrac{2 }{ L+\eta}  \right)  \| \nabla f(x_k)\|^2 .
\end{align}

Since  $ \nabla f(x^*)=0  $  and Assumption \ref{ass-fun}(ii) holds, we  recall from \cite[(2.1.19)]{nesterov2013introductory} that
\begin{align}\label{bd-Hb3}
2\eta(  f(x_k) -f(x^*)) \leq   \|  \nabla f(x_k)\|^2.
\end{align}
By  recalling that $\beta\in (0,1)$ and
$\alpha\in (0, \tfrac{2( 1-\beta) }{ L+\eta})$, the last term on the right-hand side of \eqref{rec-xp3} becomes negative.
Thus, by substituting \eqref{bd-Hb3} into \eqref{rec-xp3}, we derive
 \begin{align}\label{rec-xp4}
& \mathbb{E}[ \| x_{k+1}+p_{k+1}-x^*\|^2|\mathcal{F}_k ] +\tfrac{2\alpha  }{ 1-\beta}\left(\tfrac{\beta-\alpha \eta }{ 1-\beta}+
 \tfrac{2\eta }{\eta+L} \right) ( f(x_k)-f(x^*)) \leq  \|x_k+p_k-x^*\|^2\notag
\\&+\tfrac{2\alpha \beta }{ (1-\beta)^2}( f(x_{k-1})-f(x^*))
-\tfrac{2\alpha\eta L}{ (1-\beta)(L+\eta)} \|x_k -x^* \|^2 - \tfrac{\alpha \beta \eta}{ (1-\beta)^2}      \|x_{k-1}-x_k\|^2 + \tfrac{\alpha^2 \nu^2}{(1-\beta)^2N_k} .
\end{align}

Define $\hat{z}_k=\begin{pmatrix}
  x_{k}-x^*
\\ x_k-x_{k-1} \end{pmatrix}.$ From \eqref{def-pk} it can be seen that
\[ \| x_{k}+p_{k}-x^*\|^2=\hat{z}_k^T \mathbf{M} \hat{z}_k {\rm~with~} \mathbf{M}\triangleq  \left(
                                                                                        \begin{array}{cc}
                                                                                          \mathbf{I}_m &  \tfrac{\beta}{1-\beta}    \mathbf{I}_m \\
                                                                                       \tfrac{\beta}{1-\beta} \mathbf{I}_m & \left(\tfrac{\beta}{(1-\beta)^2}\right)^2               \mathbf{I}_m \\
                                                                                        \end{array}
                                                                                      \right).
\]
Then \eqref{rec-xp4} can be rewritten as
 \begin{align}\label{rec-xp5}
& \mathbb{E}[ \hat{z}_{k+1}^T \mathbf{M} \hat{z}_{k+1}|\mathcal{F}_k ] +\hat{m}  ( f(x_k)-f(x^*)) \leq  \hat{z}_k^T \mathbf{N} \hat{z}_k+
\hat{n} ( f(x_{k-1})-f(x^*))+ \tfrac{\alpha^2 \nu^2}{(1-\beta)^2N_k} ,
\end{align}
where
\begin{align*}
\hat{m} \triangleq \tfrac{2\alpha  }{ 1-\beta}\left(\tfrac{\beta-\alpha \eta }{ 1-\beta}+
 \tfrac{2\eta }{\eta+L} \right), ~ \hat{n}\triangleq \tfrac{2\alpha \beta }{ (1-\beta)^2},
 {\rm ~and}~\mathbf{N}\triangleq  \left(      \begin{array}{cc}
                                                                                        \left(1-\tfrac{2\alpha\eta L}{ (1-\beta)(L+\eta)}\right)  \mathbf{I}_m &  \tfrac{\beta}{1-\beta}    \mathbf{I}_m \\
                                                                                       \tfrac{\beta}{1-\beta} \mathbf{I}_m & \tfrac{\beta(\beta-\alpha\eta )}{(1-\beta)^2}             \mathbf{I}_m \\
                                                                                        \end{array}
                                                                                      \right).
\end{align*}
With the definitions of $q_1$ and $q_2$ in \eqref{def-q1q2},
it has been shown in \cite[Appendix]{ghadimi2015global} that
$\hat{n}\leq q_1 \hat{m}$ and $q_2\mathbf{M}-\mathbf{N}$ is positive semidefinite.
Then from \eqref{rec-xp5} and $q=\max\{q_1,q_2\} \in(0,1)$, we obtain that
 \[ \mathbb{E}[ \hat{z}_{k+1}^T \mathbf{M} \hat{z}_{k+1}|\mathcal{F}_k ] +\hat{m}  ( f(x_k)-f(x^*)) \leq  q\left( \hat{z}_k^T \mathbf{M} \hat{z}_k+ \hat{m} ( f(x_{k-1})-f(x^*)) \right)+ \tfrac{\alpha^2 \nu^2}{(1-\beta)^2N_k}\]
By taking unconditional expectations on both sides of the above inequality, we achieve
 \begin{align}\label{rec-xp6}
& \mathbb{E}[ \hat{z}_{k+1}^T \mathbf{M} \hat{z}_{k+1} +\hat{m}  ( f(x_k)-f(x^*)) ]\leq  q\mathbb{E} \left [\hat{z}_k^T \mathbf{M} \hat{z}_k+ \hat{m} ( f(x_{k-1})-f(x^*)) \right]+ \tfrac{\alpha^2 \nu^2}{(1-\beta)^2N_k}  \notag
\\& \leq q^{k+1}\mathbb{E} \left [\hat{z}_0^T \mathbf{M} \hat{z}_0+ \hat{m} ( f(x_{-1})-f(x^*)) \right]+ \tfrac{\alpha^2 \nu^2}{(1-\beta)^2} \sum_{i=0}^k q^i/N_{k-i}\notag
\\& =q^{k+1}\mathbb{E} \left [\|  x_{0}-x^*\|^2+ \hat{m} ( f(x_0)-f(x^*)) \right]+ \tfrac{\alpha^2 \nu^2}{(1-\beta)^2} \sum_{i=0}^k q^i/N_{k-i} ,
\end{align}
where the last inequality holds since $x_{-1}=x_0 $ and $\hat{z}_0=\begin{pmatrix}
  x_{0}-x^*
\\0 \end{pmatrix}$.

Since  $ \nabla f(x^*)=0  $ and Assumption \ref{ass-fun}(ii) holds, we have $ f(x_k)-f(x^*)\geq 2\eta\|x_k-x^*\|^2.$
Then by using  \eqref{rec-xp6} and recalling  the positive  semidefiniteness  of  $\mathbf{M}$, we obtain
 \begin{align*}
&2 \hat{m} \eta \mathbb{E}[\|x_k-x^*\|^2]\leq  q^{k+1}\mathbb{E} \left [\|  x_{0}-x^*\|^2+ \hat{m} ( f(x_0)-f(x^*)) \right]+ \tfrac{\alpha^2 \nu^2}{(1-\beta)^2} \sum_{i=0}^k q^i/N_{k-i} .
\end{align*}
Thus the required result holds.
  \hfill $\blacksquare$

\section{Proof of Proposition \ref{prp2}}
By  using $\rho_2 \in (1-\gamma,1)$ and   substituting $N_{k }=\lceil \rho_2^{-(k+1)} \rceil $  into   \eqref{lem-fy}, we obtain that
\begin{align*}
\mathbb{E}[f(y_{k }) ]-f^* & \leq         \rho_2^{k }   \tfrac{\eta+L}{2}    \mathbb{E}[\| x_0-x^*\|^2]   +  \nu^2 \left( \alpha +\tfrac{(1-\gamma)\gamma}{2\eta}\right)  \rho_2^k \sum_{i=0}^{k-1}  \left(\tfrac{  1-\gamma}{\rho_2} \right)^i  ,
\end{align*}
which leads to \eqref{lem-fy2} by using  the bound $ \sum_{i=0}^{k-1}  \left({  1-\gamma   \over \rho_2} \right)^i \leq {1 \over 1-{  1-\gamma   \over \rho_2}}= {\rho_2   \over \rho_2-(1-\gamma)}$.

  Because  $f$  is $\eta$-strongly convex and $\nabla f(x^*)=0,$  we have $f(x)-f (x^*) \geq {\eta \over 2} \|x-x^*\|^2.$ Thus, $ y_k $ generated by Algorithm  \ref{Alg_2}  satisfies
$\|y_k-x^*\|^2 \leq {2 \over \eta} (f(y_{k })  -f^* ).$
Then from \eqref{lem-fy2} it follows that   $$\mathbb{E}[\|y_k-x^*\|^2] \leq c \rho_2^k  {~~\rm with ~~} c\triangleq \tfrac{2}{\eta}\left(\tfrac{\eta+L}{2}    \mathbb{E}[\| x_0-x^*\|^2] + \tfrac{\rho_2 \nu^2}{\rho_2-(1-\gamma)}\left( \alpha +\tfrac{(1-\gamma)\gamma}{2\eta}\right)\right)  .$$  Note  by \eqref{Alg22} that  $\|x_{k}-x^*\| \leq \| (1+\beta) (y_{k}-x^*) \| +\| \beta (y_{k-1}-x^*)\| $,
 and hence
\begin{align*}\mathbb{E}[\|x_{k}-x^*\|^2] \notag & \leq 2(1+\beta)^2 \mathbb{E}[\|y_k-x^*\|^2] +2\beta^2 \mathbb{E}[\|y_{k-1}-x^*\|^2]   \\
& \leq c\left( 2(1+\beta)^2+2\beta^2 \rho_2^{-1} \right)\rho_2^k .
\end{align*}
Therefore,  sequences $\{x_k\} $ and $ \{y_k\}$ converge to the optimal solution  $x^*$ at a geometric rate $\mathcal{O}(\rho_2^k)$ in the mean-squared sense.

 Suppose we set $\alpha\triangleq{1\over L}$. Then
$\gamma = \sqrt{   \eta  \over L } ={1\over \sqrt{\kappa}}$ and $\beta=  {\sqrt{\kappa}-1\over  \sqrt{\kappa}+1}$.
Select $\rho_2 \triangleq 1-{1\over 2 \sqrt{\kappa}} $ such that $\rho_2>1-\gamma$.
We can show that  the number of iterations   required  to obtain an  $\epsilon $-optimal solution in a mean-squared sense  is  $  \mathcal{O} \left (  {   \ln \left ({1/\epsilon} \right)  \over   \ln  \left ( 1/\rho_2 \right)} \right) = \mathcal{O} \left (   \sqrt{\kappa}    \ln \left ({1/\epsilon} \right)   \right)$  since $\ln \left(1 \over 1-{1/(2 \sqrt{\kappa})}  \right)\approx {1\over 2 \sqrt{\kappa}} $ for large  $\kappa,$ and hence the oracle complexity $
\sum_{k=0}^{K(\epsilon) -1 } N_k  =\mathcal{O}\left ({\sqrt{ \kappa } / \epsilon} \right) .  $
 \hfill  $\blacksquare$

 \section{Proof of Proposition \ref{prp-alg3-Q}}
 By substituting $   \alpha ={ 4 \over (\sqrt{\eta}+\sqrt{L})^2}$ into $\beta =  \max\{|1-\sqrt{\alpha \eta}|,|1-\sqrt{\alpha L}| \}^2,$   there holds     $\beta = \left(1-{2 \over \sqrt{\kappa} +1}\right)^2<1.$  Then   Lemma \ref{Lem-HB-Q} holds. Therefore, by using  \eqref{lem-hb-q},   $N_{k  }=\lceil \rho_3^{- (k+1)} \rceil $, $x_{-1}=x_0$,  {$\rho_3 \in (  \beta ,1)$,  and  $\iota \in (0,\sqrt{\rho_3}-\sqrt{\beta})$,
 there exists a constant $ c(\iota)$ such that
\begin{align*}
\mathbb{E}& \left [ \left \| \begin{pmatrix}
  x_{k+1}-x^*
&\\ x_k-x^* \end{pmatrix} \right\|^2    \right]   \leq    2  c(\iota )^2(\sqrt{\beta}+\iota )^{2(k+1) } \mathbb{E}[ \| x_0 - x^*\|^2]+\alpha^2 \nu^2  c(\iota )^2 \sum_{t=0}^{k}(\sqrt{\beta}+\iota )^{2(k-t)} \rho_3^{t+1} , \quad \forall k\geq 0.
\end{align*}
Since  $\iota \in (0,\sqrt{\rho_3}-\sqrt{\beta})$, we have  $\sqrt{\beta}+\iota \in (\sqrt{\beta}, \rho_3).$
This together with  $\sum_{t=0}^{k}(\sqrt{\beta}+\iota )^{2(k-t)} \rho_3^{t+1}=\sum_{t=0}^{k}  (\sqrt{\beta}+\iota )^{2t} \rho_3^{k+1-t} =\rho_3^{k+1}\sum_{t=0}^{k}  \left ( \tfrac{(\sqrt{\beta}+\iota )^{2}}{\rho_3}\right)^{t} \leq{\rho_3^{k+1}\over 1-(\sqrt{\beta}+\iota )^{2}/\rho_3}$,  proving  \eqref{prp-hb-q}.}

 By \eqref{prp-hb-q}, $ \mathbb{E} [\| x_k-x^* \|]\leq  c \rho_3^{k}  $ for some constant $c>0.$ Suppose $\rho_3 = \left(1-{1 \over \sqrt{\kappa} +1}\right)^2> \beta$. Akin  to the proof of   Proposition \ref{prp1}(ii), we can show that  the number of iterations   required  to obtain an  $\epsilon $-optimal solution satisfying   $\mathbb{E}[\|x_k-x^*\|^2] \leq \epsilon$   is   $  \mathcal{O} \left (  {   \ln \left ({1/\epsilon} \right)  \over   \ln  \left ( 1/\rho_3 \right)} \right) = \mathcal{O} \left (   \sqrt{\kappa}    \ln \left ({1/\epsilon} \right)   \right)$  since $\ln \left(1 \over 1-{1/(  \sqrt{\kappa}+1)}  \right)\approx {1\over   \sqrt{\kappa}+1} $ for large  $\kappa,$ and the oracle complexity $
\sum_{k=0}^{K(\epsilon) -1 } N_k  =\mathcal{O}\left ({\sqrt{ \kappa } / \epsilon} \right) .  $
\hfill $\blacksquare$

 \section{Proof of Proposition \ref{prp-alg3}}

By substituting  $N_{k  }=\lceil \rho_4^{- (k+1)} \rceil $  into \eqref {lem-hb}, we obtain
 \begin{align*}  \mathbb{E}[\|x_k-x^*\|^2]&\leq \tfrac{1}{2 \hat{m} \eta}\mathbb{E} \left [\|  x_{0}-x^*\|^2+ \hat{m} ( f(x_0)-f(x^*)) \right] q^{k+1}+ \tfrac{\alpha^2 \nu^2}{2(1-\beta)^2 \hat{m} \eta} \rho_4^{k+1} \sum_{i=0}^k  (q/\rho_4)^i
 \\&\leq \tfrac{1}{2 \hat{m} \eta}\mathbb{E} \left [\|  x_{0}-x^*\|^2+ \hat{m} ( f(x_0)-f(x^*)) \right]  \rho_4^{k+1}+ \tfrac{\alpha^2 \nu^2}{2(1-\beta)^2 \hat{m} \eta(1-q/\rho_4)} \rho_4^{k+1}  ,
\end{align*}
where the last inequality holds by $\rho_4\in (q,1). $  Then \eqref{prp-hb} holds. \hfill $\blacksquare$

\section{Proof of Proposition \ref{prp3}}
 By substituting   $N_k   \triangleq \lceil (k+1)^v\rceil $ into \eqref{recursion1}, using   $q= 1-{2 \alpha \eta L \over \eta+L} $ and Lemma \ref{lem-recur}(i), one obtains
 \begin{align*}
 \mathbb{E}[\| x_{k}-x^*\|^2  ]  & \leq  q \mathbb{E}[\| x_{k-1} - x^*\|^2]+\alpha^2\nu^2 k^{-v} = q^{k} \mathbb{E}[ \| x_0 - x^*\|^2] +\alpha^2\nu^2 \sum_{m=1}^{k} q^{k-m} m^{-v} \\ & =  q^k \mathbb{E}[ \| x_0 - x^*\|^2]  +\alpha^2\nu^2  \left(q^{k } \tfrac{e^{2v}q^{-1}-1}{1-q}+\tfrac{2   k^{-v}}{q \ln(1/q) }\right) ,
\end{align*}
which together with  Lemma \ref{lem-recur}(ii) proves  \eqref{rate-pol}.
 Then  $  \mathbb{E}[ \|x_{k}-x^*\|^2  ] \leq \epsilon$ for any $  k\geq  K(\epsilon)\triangleq  \left( { C_v  \over \epsilon }\right)^{1/v}$.    By noting  that $C_v=\mathcal{O}(e^v v^v)$,   the iteration complexity is  $\mathcal{O} (  {v} (1/\epsilon)^{1/v})$.
Therefore, the  number of sampled gradients required to obtain    an  $\epsilon-$NE  is   bounded by
\begin{equation*}
\begin{array}{l} \sum_{k=0}^{ K(\epsilon) -1 }   \lceil  (k+1)^v \rceil
  \leq  K(\epsilon) + (K(\epsilon))^v + \sum_{k=1}^{ K(\epsilon) -1 }    k^v
\leq  \left( { C_v  \over \epsilon }\right)^{1/v}+\tfrac{ C_v}{\epsilon } +  \int_{1}^{ K(\epsilon)  } t^v dt
\\ \qquad \qquad \qquad ~\quad= \left( { C_v  \over \epsilon }\right)^{1/v}+\tfrac{ C_v}{\epsilon }+\tfrac{t^{v+1} }{ v+1}\Big |_{1}^{ K(\epsilon) }
 =\left( { C_v  \over \epsilon }\right)^{1/v}+\tfrac{ C_v}{\epsilon } + (v+1)^{-1} \left(\tfrac{ C_v}{\epsilon }\right)^{1+{1\over v}} - (v+1)^{-1} .
  \end{array}
  \end{equation*}
Therefore, the oracle complexity is $\mathcal{O} \left(  {e^v v^v}\left({ 1/\epsilon }\right)^{1+{1\over v}}  \right).$
\hfill $\blacksquare$

\section{Proof of Proposition \ref{prp4}}
 It is noticed   by $t = k-i$ and Lemma \ref{lem-recur}(i),
\begin{align}\label{bd-rk}
  \sum_{i=0}^{k-1}  (1-\gamma )^i (k-i)^{-v} &   =    \sum_{t=1}^{k } (1-\gamma )^{k-t} t^{-v} \leq      \tfrac{(1-\gamma )^{k  }(e^{2v} (1-\gamma )^{-1}-1)}{\gamma }+\tfrac{2  k^{-v}}{(1-\gamma ) \ln(1/ (1-\gamma )) } .
\end{align}
By substituting $N_{k }=\lceil (k+1)^v \rceil $  into  \eqref{lem-fy},    we obtain that
\begin{align*}
\mathbb{E}[f(y_{k }) ]-f^* & \leq          \tfrac{(1-\gamma)^{k }(\eta+L)\mathbb{E}[\| x_0-x^*\|^2]}{2}       +  \nu^2  \left( \alpha +\tfrac{(1-\gamma)\gamma}{2\eta}\right)  \sum_{i=0}^{k-1}  (1-\gamma )^i (k-i)^{-v} .
\end{align*}
This  together   with   Lemma \ref{lem-recur}(ii) and \eqref{bd-rk}   proves the  required result.\hfill $\blacksquare$

\section{Proof of Proposition \ref{prp2-alg3-Q}}
By substituting $N_k=\lceil (k+1)^v\rceil $  into \eqref{lem-hb-q},  we obtain that
{for any $\iota \in (0,1-\sqrt{\beta})$,
 there exists a constant $ c(\iota )$ such that for any $k\geq 0,$
\begin{align*}
\mathbb{E} \left [ \left \| \begin{pmatrix}
  x_{k+1}-x^*
\\ x_k-x^* \end{pmatrix} \right\|^2    \right]  &   \leq 2   (c(\iota ))^2 (\sqrt{\beta}+\iota )^{2(k+1) } \mathbb{E}[ \| x_0 - x^*\|^2]+\alpha^2 \nu^2 (c(\iota ))^2  \sum_{t=1}^{k+1} (\sqrt{\beta}+\iota )^{2(k+1-t)} t^{-v}
 \\& \leq 2  (c(\iota ))^2 (\sqrt{\beta}+\iota )^{2(k+1) } \mathbb{E}[ \| x_0 - x^*\|^2]
 \\&+\alpha^2 \nu^2 (c(\iota ))^2   \left((\sqrt{\beta}+\iota )^{2( k+1)} {e^{2 v}(\sqrt{\beta}+\iota )^{-2}-1 \over 1-(\sqrt{\beta}+\iota )^2}+{    2(k+1)^{-v}\over 2 (\sqrt{\beta}+\iota )^2  \ln(1/(\sqrt{\beta}+\iota )) }\right),
\end{align*}
where the last inequality holds by   using  Lemma \ref{lem-recur}(i) and recalling  that $\sqrt{\beta}+\iota \in (\sqrt{\beta},1)$.
This  together with  Lemma \ref{lem-recur}(ii) yields the result.}
\hfill $\blacksquare$

\section{Proof of Proposition \ref{prp2-alg3}}
By substituting $N_k=\lceil (k+1)^v\rceil $  into \eqref{lem-hb},  we obtain that {
\begin{align*}
&\mathbb{E}[\|x_k-x^*\|^2]\leq \tfrac{1}{2 \hat{m} \eta}\mathbb{E} \left [\|  x_{0}-x^*\|^2+ \hat{m} ( f(x_0)-f(x^*)) \right] q^{k+1}+ \tfrac{\alpha^2 \nu^2}{2(1-\beta)^2 \hat{m} \eta} \sum_{i=0}^k q^i(k+1-i)^{-v}, \forall k\geq 0 .
\end{align*}
By using    Lemma \ref{lem-recur}(i), we derive $\sum_{i=0}^k q^i(k+1-i)^{-v}=\sum_{t=1}^{k+1} q^{k+1-t}  t^{-v} \leq q^{k+1  }
\tfrac{e^{2v}q^{-1}-1 }{ 1-q}+\tfrac{2  (k+1)^{-v}}{q \ln(1/q) }$.
Thus,
\begin{align*}
&\mathbb{E}[\|x_k-x^*\|^2]\leq \tfrac{1}{2 \hat{m} \eta}\mathbb{E} \left [\|  x_{0}-x^*\|^2+ \hat{m} ( f(x_0)-f(x^*)) \right] q^{k+1}+ \tfrac{\alpha^2 \nu^2}{2(1-\beta)^2 \hat{m} \eta} \left( q^{k+1  }
\tfrac{e^{2v}q^{-1}-1 }{ 1-q}+\tfrac{2  (k+1)^{-v}}{q \ln(1/q) }\right)  .
\end{align*} This  together with  Lemma \ref{lem-recur}(ii) yields the result.}
\hfill $\blacksquare$

\end{document}